\documentclass{amsproc}
\usepackage{color}
\newtheorem{theorem}{Theorem}[section]
\newtheorem{lemma}[theorem]{Lemma}
\newtheorem{proposition}[theorem]{Proposition}
\newtheorem{corollary}[theorem]{Corollary}

\theoremstyle{definition}
\newtheorem{definition}[theorem]{Definition}

\newtheorem{assumption}[theorem]{Assumption}

\theoremstyle{remark}

\numberwithin{equation}{section}

\newcommand{\abs}[1]{\lvert#1\rvert}

\newcommand{\bd}{\noindent {\sc Proof}.\ \ }

\begin{document}

\title[Nontrivial attractor-repellor maps]
{Nontrivial attractor-repellor maps of $S^2$ and rotation numbers}

\author{Shigenori Matsumoto}
\address{Department of Mathematics, College of
Science and Technology, Nihon University, 1-8-14 Kanda, Surugadai,
Chiyoda-ku, Tokyo, 101-8308 Japan}
\email{matsumo@math.cst.nihon-u.ac.jp}

\thanks{The author is partially supported by Grant-in-Aid for
Scientific Research (C) No.\ 20540096.}

\date{\today}

\newcommand{\AAA}{{\mathbb A}}
\newcommand{\BBB}{{\mathbb B}}
\newcommand{\LL}{{\mathcal L}}
\newcommand{\MCG}{{\rm MCG}}
\newcommand{\PSL}{{\rm PSL}}
\newcommand{\R}{{\mathbb R}}
\newcommand{\Z}{{\mathbb Z}}
\newcommand{\XX}{{\mathcal X}}
\newcommand{\per}{{\rm per}}
\newcommand{\N}{{\mathbb N}}

\newcommand{\PP}{{\mathcal P}}
\newcommand{\GG}{{\mathbb G}}
\newcommand{\FF}{{\mathcal F}}
\newcommand{\EE}{{\mathbb E}}
\newcommand{\BB}{{\mathbb B}}
\newcommand{\CC}{{\mathcal C}}
\newcommand{\HH}{{\mathcal H}}
\newcommand{\UU}{{\mathcal U}}
\newcommand{\oboundary}{{\mathbb S}^1_\infty}
\newcommand{\Q}{{\mathbb Q}}
\newcommand{\DD}{{\mathcal D}}
\newcommand{\rot}{{\rm rot}}
\newcommand{\Cl}{{\rm Cl}}
\newcommand{\Index}{{\rm Index}}
\newcommand{\Int}{{\rm Int}}
\newcommand{\Fr}{{\rm Fr}}

\date{\today }

\maketitle

\begin{abstract}
We consider an orientation preserving homeomorphism $h$ of $S^2$
which admits a repellor denoted $\infty$ and an attractor $-\infty$
such that $h$ is not a North-South map 
and that the basins of $\infty$ and $-\infty$ intersect.
We study various aspects of the rotation number of
$h:S^2\setminus\{\pm\infty\}\to S^2\setminus\{\pm\infty\}$,
especially its relationship with the existence of periodic orbits.
\end{abstract}

\section{Introduction}
Let $h$ be a homeomorphism of the 2-sphere $S^2$. A fixed
point $a$ of $h$ is called an {\em attractor} if there is an open disk $V$
containing $a$ such that $h({\rm Cl}(V))\subset V$ and $ 
\bigcap_{i\in\N}h^i(V)=\{a\}$. For such $V$, the set $ 
W_a=\bigcup_{i\in\N}h^{-i}(V)$ is called the {\em basin} of $a$.
A point $x$ of $W_a$ is characterised by the property: 
$  \lim_{i\to\infty}h^{i}(x)=a$.
An attractor $b$ of the inverse $h^{-1}$ is called a {\em repellor}
of $h$, and its basin $W_b$ is defined likewise.
The basins are invariant by $h$ and homeomorphic to open disks.

Let $\infty$ and $-\infty$ be distinct points of $S^2$.

\begin{definition} \label{HH}
A homeomorphism $h$ of $S^2$
which satisfy the following conditions 
is called a {\em nontrivial attractor-repellor map}.
\begin{enumerate}
\item $h$ is orientation preserving.

\item $-\infty$ is an attractor of $h$ with basin $W_{-\infty}$ and
$\infty$ a repellor with basin $W_\infty$.

\item $Z=S^2\setminus(W_{-\infty}\cup W_\infty)$ is nonempty.

\item  $W_{-\infty}\cap W_\infty\neq\emptyset$.

\end{enumerate}
\end{definition}

Condition (3) is equivalent to saying that $h$ is
{\em not} a North-South map.
Condition (4) is equivalent to saying that there is {\em no} $h$-invariant
continuum separating $-\infty$ and $\infty$. 
Denote by $\HH$ the set of nontrivial attractor-repellor maps.

Let $W_{\pm\infty}^*=W_{\pm\infty}\cup\partial W^*_{\pm\infty}$ be
the prime end compactification of $W_{\pm\infty}$, 
where
$\partial W^*_{\pm\infty}$ is
the set of prime ends of $W_{\pm\infty}$.
The homeomorphism $h\in\HH$ induces a homeomorphism $h_{\pm\infty}^*$ of 
$W_{\pm\infty}^*$. See Section 2 for more details.

The open annulus $\AAA=S^1\times\R$ is identified with
$S^2\setminus\{\pm\infty\}$ in such a way that the end
$S^1\times\{\pm\infty\}$ is identified with the deleted point $\pm\infty$.
Then $h$ induces an orientation and end preserving homeomorphism of
$\AAA$, which we still denoted by $h$. The set
$U_{\pm\infty}=W_{\pm\infty}\setminus\{\pm\infty \}$ is
 considered to be a subset of $\AAA$. Denote 
$U_{\pm\infty}^*=W_{\pm\infty}^*\setminus \{\pm\infty\}$.

The universal covering space of
$U_{\pm\infty}$ is
defined as the set of the homotopy classes of paths
from the base point, and is considerd 
to be simultaneously a subspace of $\tilde\AAA$, the universal covering space of $\AAA$
and of $\tilde U^*_{\pm\infty}$, the universal covering space of
$U^*_{\pm\infty}$.

Denote the both covering maps
by $\pi:\tilde\AAA\to\AAA$ and
$\pi:\tilde U_{\pm\infty}^*\to U^*_{\pm\infty}$. 
The inverse image  $\pi^{-1}(U_{\pm\infty})$ is 
simultaneously considered to be a subspace of
$\tilde\AAA$ and of $\tilde U_{\pm\infty}^*$.

Fix once for
all a lift $\tilde h:\tilde\AAA\to\tilde\AAA$ of $h$. 
Corresponding to $\tilde h$, a lift $\tilde
h^*_{\pm\infty}:\tilde U_{\pm\infty}^*\to\tilde U_{\pm\infty}^*$ of $h^*_{\pm\infty}
$ is specified in such a way that
they coincide on 
$\pi^{-1}(U_{\pm\infty})$ under the above identification.

The rotation number (taking value in $\R$)
of
the restriction of $\tilde h_{\pm\infty}^*$
to the boundary $\partial \tilde U^*_{\pm\infty}=\pi^{-1}(\partial U^*_{\pm\infty})$
is called 
the {\em prime end rotation number of $\tilde h$ at $\pm\infty$} 
and is denoted by $\rot(\tilde h,\pm\infty)$.

In \cite{OR} it is shown that if 
one of the prime end rotation numbers, say
$\rot(\tilde h,\infty)$ of $h\in\HH$ is rational, then
there are periodic points in $Z$. 
In \cite{HOR} a 
partial
converse is shown:
if $\rot(\tilde h,\infty)$ is irrational and if the point $-\infty$
is {\em accessible} from $W_{\infty}$, then there is no periodic
points
in $Z$. The second condition means that there is a path $\gamma:[0,1]\to
S^2$ such that $\gamma([0,1))\subset W_{\infty}$ and $\gamma(1)=-\infty$.
Our first result shows that the accessibility condition
is actually necessary, contrary to a conjecture therein. 

\begin{theorem} \label{main}
Given any real numbers $\alpha$ and $\beta$, there is a homeomorphism
$h\in\HH$ with its lift $\tilde h$ 
such that $\rot(\tilde h,\infty)=\alpha$ and $(\tilde h,-\infty)=\beta$.
\end{theorem}

The accessibility condition is necessary since if we choose
$\alpha$ to be rational and $\beta$ irrational, there is
a periodic point (\cite{OR}) and thus the irrationality of one prime end rotation number
does dot imply the nonexistence of periodic point.

\medskip
A nontrivial attractor repellor map $h$ has a structure similar to
a gradient flow. Most relevant to this structure is the chain recurrent
set $C$ of $h$. Except $\pm\infty$, $C$ is contained in $Z$, and
partitioned into the union of chain transitive classes. Each chain
transitive class is closed and $h$-invariant. See Section 3 for
a review of these concepts.

The example in Theorem \ref{main} constructed in
Section 2 shows that Poincar\'e-Birkhoff type theorem does not
hold for $h\in\HH$. But when restricted to a single chain transitive class,
we get a variant of it.

We consider $h$ to be a homeomorphism
of the annulus $\AAA$, and
fix a lift $\tilde h:\tilde\AAA\to\tilde\AAA$ of $h$. 
Then for any $h$-invariant
probability measure $\mu$ of $Z$, the rotation number $\rot(\tilde h,\mu)$
is defined as follows. Denote by $\Pi_1:\tilde\AAA\to\R$ the projection
onto the first factor. Then the function $\Pi_1\circ\tilde h-\Pi_1$ is
invariant under the covering transformations, and hence defines a
function on $\AAA$. We set 
$$
\rot(\tilde h,\mu)=\langle \mu,\Pi_1\circ\tilde h-\Pi_1\rangle.
$$
For a periodic point $x$ of $h$,
we denote by $\rot(\tilde h,x)$ the rotation number $\rot(\tilde h,\mu)$
for $\mu$ the average of the point masses along the orbit of $x$.

\begin{theorem} \label{main2}
Suppose $x_1$ and $x_2$ are periodic points of $h$ belonging
to the same chain transitive class $C_0$ such that
$\rot(\tilde h,x_\nu)=\alpha_\nu$ ($\nu=1,2$).
Then for any rational number
$\alpha\in[\alpha_1,\alpha_2]$ there
is a periodic point $x$ in $C_0$ such that $\rot(\tilde h,x)=\alpha$.
\end{theorem}

Let us define the {\em rotation set}\, $\rot(\tilde h,C_0)$ of
a chain transitive class $C_0$ as
the set of the values $\rot(\tilde h,\mu)$, where $\mu$
runs over the space of $h$-invariant probability measures
supported on  $C_0$.
The rotation set  is a closed interval or a singleton.

\begin{corollary} \label{c1}
Suppose $C_0$ is a chain transitive class with
$\rot(\tilde h,C_0)=[\alpha_1,\alpha_2]$, where $\alpha_\nu$
 are distinct rational numbers. 
Then for any rational number
$\alpha\in[\alpha_1,\alpha_2]$ there
is a periodic point $x$ in $C_0$ such that $\rot(\tilde h,x)=\alpha$.
\end{corollary}

The proof of Theorem \ref{main2} and Corollary \ref{c1}, as well as an example of $h\in\HH$
which shows that Theorem \ref{main2} is nonvoid is given in Section 3.
The author cannot improve Corollary \ref{c1} so as to include the case
where
$\alpha_\nu$ is irrational.

\medskip
Next we study an influence of the prime end rotation number 
$\rot(\tilde h,\infty)$ on the dynamics of $h$ on $Z$.

\begin{theorem}\label{main3}
If $\rot(\tilde h,\infty)=p/q$ {\em ($(p,q)=1$)}, there is a periodic point 
$x\in \AAA$
of period $q$  such that $\rot(\tilde h,x)=p/q$.
\end{theorem}

This is already known (\cite{OR}) except for the last statement about 
the rotation number.
Section 4 is devoted to the proof of Theorem \ref{main3}.

\medskip
Our last theorem is concerned about the case where $-\infty$
is accessible from $U_\infty$. 
Then the dynamics of $h$ on $Z$
is shown to be quite simple in the view point of rotation numbers. 
This
is a refinement of a result in \cite{HOR} cited above.
The proof is given in Section 5.

\begin{theorem}\label{main4}
Assume that $-\infty$ is accessible from $U_\infty$ and let
$\alpha=\rot(\tilde h,\infty)$. Then
\begin{enumerate}
\item $\rot(\tilde h,\mu)=\alpha$ for any $h$-invariant probability
      measure supported on $Z$.
\item $\rot(\tilde h,-\infty)=\alpha$.
\end{enumerate}
\end{theorem}

\bigskip
{\sc Acknowledgement:} Hearty thanks are due to the referee,
for valuable comments.

\section{Prime end rotation numbers}

\noindent
{\bf 2.1.} First of all, we
recall fundamental facts about the prime end compactification of
$W_{\pm\infty}$. See \cite{E,M,P,MN} for an detailed exposition. 

A properly embedded copy of the real line $c$ in $W_{\pm\infty}$ which does not
pass through ${\pm\infty}$ is called a {\em cross cut} of
$W_{\pm\infty}$. 
The word ``proper'' means that the inverse image of any compact set
is compact. The
connected component of the complement of a cross cut $c$ which does not contain the point
${\pm\infty} $ is denoted by $V(c)$. A sequence $\{c_i\}_{i\in\N}$ of cross
cuts is called
a {\em topological chain} (\cite{M}) if the following conditions are satisfied.
\begin{enumerate}
\item $c_{i+1}\subset V(c_i)$, $\forall i\in\N$.
\item ${\rm Cl}(c_i)\cap {\rm Cl}(c_j)=\emptyset$ if $i\neq j$, where
${\rm Cl}(\cdot)$ denotes the closure in $S^2$.
\item ${\rm diam}(c_i)\to0$ as $i\to\infty$, where the diameter is taken
      with respect to the spherical metric of $S^2$.
\end{enumerate}
 Two topological chains $\{c_i\}$ and $\{c'_i\}$ are said to be {\em
 equivalent} if for any $i$, there is $j$ such that $c_j\subset V(c'_i)$
and $c_j'\subset V(c_i)$. 

An equivalence class of topological chains is
 called a {\em prime end} of $W_{\pm\infty}$. The set of prime ends is
 denoted by $\partial W^*_{\pm\infty}$. The set
 $W_{\pm\infty}^*=W_{\pm\infty}\cup\partial W^*_{\pm\infty}$ is called the {\em
 prime end compactification of $W_{\pm\infty}$}.
It is topologized as follows. A
neighbourhood system in $W_{\pm\infty}^*$ of a point in $W_{\pm\infty}$
is  the same as a given system for $W_{\pm\infty}$.
Choose a point $\xi$ in $\partial W^*_{\pm\infty}$ represented by a topological chain
 $\{c_i\}$. The set of points in $V(c_i)$, together with the prime ends
 represented by topological chains contained in $V(c_i)$, for each $i$,
 forms a fundamental neighbourhood system of $\xi$. It is a classical fact due
 to Carath\'eodory that $W_{\pm\infty}^*$ is homeomorphic to a closed disk.

It is clear by the topological nature of
the definition that the homeomorphism $h$ of $S^2$ induces a
homeomorphism $h^*_{\pm\infty}$ of $W_{\pm\infty}^*$.

\medskip
\noindent {\bf 2.2.}
Now let us embark upon the costruction of the homeomorphism
$h\in\HH$ in Theorem \ref{main}.
We shall construct it as a homeomorphism of the annulus $\AAA$.
Roughly speaking, on the subannulus $S^1\times[5,\infty)$,
$h$ is of the form
$$
h(\theta,t)=(f_\alpha(\theta), t-g(\theta,t)),$$
where $f_\alpha$ is a rigid rotation of $S^1$ if $\alpha$ is rational,
and a Denjoy hemeomorphism if irrational. By choosing the $[0,1]$-valued
function
$g$ appropriately, one can form the homeomorphism $h$ which has
a unique minimal set on the level $t=10$. Also
$h$ satisfies
$$
h(S^1\times[5,\infty))=S^1\times[4,\infty).$$

Likewise we define $h$ on $(-\infty,-5]$ using a homeomorphism $f_\beta$
of $S^1$ of rotation number $\beta$. 
It has a unique minimal set on the level $t=-10$.
Finally on $S^1\times[-5,5]$,
we define $h$ as
$$h(\theta,t)=(\varphi_t(\theta),t-1),$$
by using an isotopy  $\varphi_t$ ($t\in[-5,5]$) joining $f_\beta$ and
$f_\alpha$. 

\medskip

Let us start a concrete construction.
Given $\alpha\in\R$, let 
us define a homeomorphism $f_\alpha$ of $S^1$ 
and its lift $\tilde f_\alpha:\R\to\R$ with rotation number 
$\rot(\tilde
f_\alpha)=\alpha$
as follows. For $\alpha$ rational, let $\tilde f_\alpha$
be the translation by $\alpha$. Thus $f_\alpha$
is the rigid rotation of $S^1$.
For  $\alpha$ irrational,
let $f_\alpha$ be
a Denjoy homeomorphism and $\tilde f_\alpha$ the lift of $f_\alpha$
such that $\rot(\tilde f)=\alpha$. 
Let $C_\alpha\subset S^1$ be a minimal set 
of $f_\alpha$. Thus $C_\alpha$ is a single periodic orbit if $\alpha$
is rational, and a Cantor set if $\alpha$ is irrational.
For $\alpha$ irrational, we assume furthermore that 
the complement of $C_\alpha$ consists of the orbit of
a single wandering interval. That is, there
is a connected component $I_\alpha$ of $S^1\setminus C_\alpha$ such that
$  \bigcup_{i\in\Z}f^i_\alpha(I_\alpha)=S^1\setminus C_\alpha$.

Define a continuous function $g_\alpha:S^1\to[0,1]$ such that

(a) $g_\alpha^{-1}(0)=C_\alpha$, and

(b) for any $\theta\in S^1\setminus C_\alpha$, 
$  \sum_{i\geq0}g_\alpha(f_\alpha^i(\theta))=\infty$ and
$  \sum_{i\leq0}g_\alpha(f_\alpha^i(\theta))=\infty$.

The existence of such $g_\alpha$ is clear for $\alpha$ rational.
For $\alpha$ irrational, first define $g_\alpha$ on the
interval ${\rm Cl}(I_\alpha)$ so that 
$g_\alpha^{-1}(0)=\partial I_\alpha$. 
For any $i\in\Z\setminus\{0\}$, define $g_\alpha$ on
$f_\alpha^i(I_\alpha)$ by
$g_\alpha(f_\alpha^i(\theta))=\abs{i}^{-1}g_\alpha(\theta)$.
Finally set $g_\alpha=0$ on $C_\alpha$.
Then $g_\alpha$ is continuous and satisfies (a) and (b).

For $\beta\in\R$, we define $f_\beta$, $\tilde f_\beta$, $C_\beta$ 
and $g_\beta$
likewise. 

Define a continuous function $g:S^1\times\R\to[0,1]$,
differentiable along the $\R$-direction, such that

(c) $g^{-1}(0)=C_\alpha\times\{10\}\cup C_\beta\times\{-10\}$,

(d) {\color{red}for $t\in[9,11]$, $g(\theta,t)\geq g_\alpha(\theta)$ and 
for $t\in[-11,-9]$, $g(\theta,t)\geq g_\beta(\theta)$, with the equality
only for $t=\pm10$,}

(e)  $g=1$ on $S^1\times((-\infty, 15]\cup[-5,5]\cup[15,\infty))$ and

(f) $\partial g/\partial t<1$.

Choose a continuous family $\varphi_t$ ($t\in\R$) of
homeomorphisms
of $S^1$ and its continuous lift $\tilde \varphi_t$
such that 

(g) $\tilde\varphi_t=\tilde f_\alpha$ for $t\in[5,\infty)$
and $\tilde \varphi_t=\tilde f_\beta$ for $t\in(-\infty,-5]$.

Finally define a homeomorphism $h:S^1\times\R\to S^1\times\R$ by
$$
h(\theta,t)=(\varphi_t(\theta),t-g(\theta,t)).$$

\medskip
\noindent {\bf 2.3.}
We shall show that $h$ satisfies the conditions of Theorem \ref{main}.
First let us verify that $h$ is a homeomorphism of $\AAA$.
Clearly $h$ is continuous and by (e) maps the circle $S^1\times\{15\}$
(resp.\ $S^1\times\{-15\}$) onto $S^1\times\{14\}$ (resp.\
$S^1\times\{-16\}$). This shows that $h$ is surjective.
To show that $h$ is injective, assume $h(\theta_1,t_1)=h(\theta_2,t_2)$.
Then by (e) $h$ maps $S^1\times[5,\infty)$, $S^1\times(-\infty,-5]$
and $S^1\times[-5,5]$ respecively onto $S^1\times[4,\infty)$, $S^1\times(-\infty,-6]$
and $S^1\times[-6,4]$. Therefore the two points $(\theta_1,t_1)$ and
$(\theta_2,t_2)$ must simultaneously belong to either one of the
subannuli $S^1\times[5,\infty)$, $S^1\times(-\infty,-5]$
and $S^1\times[-5,5]$. In the first case we have
$$
h(\theta_i,t_i)=(f_\alpha(\theta_i),t_i-g(\theta_i,t_i)),$$
and thus $\theta_1=\theta_2$. On the other hand by (f), 
$h\vert_{\{\theta\}\times\R}$ is injective, showing that $t_1=t_2$.
The second case can be dealt with similarly.

In the last case, we have 
$$h(\theta_i,t_i)=(\varphi_{t_i}(\theta_i),t_i-1).$$
Thus $t_1=t_2$, which implies $\theta_1=\theta_2$.

Next let us show that $h\in\HH$. Conditions (1) $\sim$ (3)
of Definition \ref{HH} are clear. Let us show (4).
Consider the basin
$W_\infty$ of the repellor $\infty$ (corresponding to the end
$S^1\times\{\infty\}$ of the cylinder $\AAA$). Recall the notation
$U_\infty=W_\infty\setminus\{\infty\}\subset\AAA$.
We shall show
\begin{equation} \label{e1}
U_\infty\cap (S^1\times[5,\infty))=(S^1\times[5,\infty))\setminus
(C_\alpha\times [5,10]).
\end{equation}

{\color{red}To show this, first notice that the minimum value of $g$ on 
$S^1\times[10+\varepsilon,\infty)$ is positive for any $\varepsilon>0$.
This shows that $S^1\times(10,\infty)\subset U_\infty$. Next by (d) and (b),
any point in $(S^1\setminus C_\alpha)\times(10,11)$ can be moved below the
level $t=10$ by an iterate of $h$.  Since 
$$U_\infty=\bigcup_{i\in\N}h^i(S^1\times(10,\infty)),$$ 
we have 
$$(S^1\times[5,\infty))\setminus
(C_\alpha\times [5,10])\subset U_\infty.
$$
On the other hand, the opposite inclusion of (\ref{e1}) is easy to show.}

The basin $U_\infty$ is obtained as the increasing union of the images of
the set in (\ref{e1}) by the positive iterates of $h$. Therefore it is clear
that $U_\infty\cap(S^1\times(-5,5))$ is open and dense in $S^1\times(-5,5)$.
Likewise we can prove that $U_{-\infty}\cap(S^1\times(-5,5))$ is open and dense
in $S^1\times(-5,5)$. This shows that $U_\infty\cap U_{-\infty}\neq\emptyset$,
as is required.

What is left is to show that $\rot(\tilde h,\infty)=\alpha$, 
the other assertion
$\rot(\tilde h,-\infty)=\beta$ being proven similarly.

Now for any point $\theta\in S^1$, define the ray
$r_\theta:(0,\infty)\to U_\infty$ by
$$r_\theta(t)=(\theta, t^{-1}+10).$$ 
For $\theta\not\in C_\alpha$, the end point
$r_\theta(\infty)=(\theta,10)$
is a point in $U_\infty$. For $\theta\in C_\alpha$, the end point 
$r_\theta(\infty)$ is
defined as a prime end, i.\ e.\
a point of $\partial U^*_\infty(=\partial W^*_\infty)$ as
follows. 

For any $i\in\N$, let $S_i$ be
the circle centered at $(\theta,10)$ and of radius $i^{-1}$, and $c_i$
the cross cut of $U_\infty$ obtained as the connected component of
$S_i\cap U_\infty$ that intersects the ray $r_\theta$.
Clearly $\{c_i\}$ is a topological chain. Denote the prime end it
determines by
$r_\theta(\infty)$. 

Define a map $\gamma:S^1\to U^*_\infty$ by
$\gamma(\theta)=r_\theta(\infty)$.
The map  $\gamma$ is clearly injective.
It  is also continuous according to the definition of the topology
of $U_\infty^*$.
The intersection $C^*$ of  the curve $\gamma$ with the set of prime ends
$\partial U^*_\infty$
is either a finite set or a Cantor set, and $\gamma$
maps $C_\alpha$ homeomorphically onto
$C^*$
in a way to preserve the cyclic order  and conjugates
$f_\alpha\vert_{C_\alpha}$ to $h_\infty^*\vert_{C^*}$.
Moreover there is a lift of $\gamma$ defined on $\R$ 
taking values on $\tilde U_\infty^*$ which maps
$\pi^{-1}(C_\alpha)$ homeomorphically onto $\pi^{-1}(C^*)$ in an order preserving
way and conjugates
$\tilde f_\alpha\vert{\pi^{-1}(C_\alpha)}$ to 
$\tilde h_\infty^*\vert{\pi^{-1}(C^*)}$.
Since the lift $\tilde h$ of $h$ is determined by (g), we have
$\rot(\tilde h,\infty)=\alpha$, completing the proof of Theorem
\ref{main}.

\section{The rotation set of a chain transitive class}

\noindent
{\bf 3.1.} Fix $h\in\HH$. For $\varepsilon>0$ and $x,\,y\in S^2$, a sequence $\{x=x_0,x_2,\ldots,x_r=y\}$
of points of $S^2$ is called an  {\em $\varepsilon$-chain of $h$
of length $r$
from $x$ to $y$} if for any $0\leq i\leq n-1$, $d(h(x_{i}), x_{i+1})<\varepsilon$,
and an  {\em $\varepsilon$-cycle at} $x$ if furthermore $x=y$.
A point $x\in S^2$ is called {\em chain recurrent} if for any
$\varepsilon>0$,
there
is an $\varepsilon$-cycle at $x$.
The set $C$ of the chain recurrent points is called the {\em chain recurrent
set}. It is a closed set invariant by $h$. 

Two points $x$ and $y$ of $C$
are said to be {\em chain transitive}, denoted $x\sim y$,
if for any $\varepsilon>0$, there are an $\varepsilon$-chain from $x$ to $y$
and another from $y$ to $x$. An equivalence class of $\sim$ is called
a {\em chain transitive class}.
Again it is  closed and invariant by $h$.

\begin{definition} \label{d2}
A continuous function $H:S^2\to\R$ is called a {\em complete Lyapunov function}
of $h$ if it satisfies the following conditions.
\begin{enumerate}
\item If $x\not\in C$, then $H(h(x))< H(x)$.
\item If $x,\,y\in C$, $H(x)=H(y)$ if and only if $x\sim y$.
\item The set of values $H(C)$ is closed and Lebesgue null in $\R$.
\end{enumerate}
A value in $\R\setminus H(C)$ is called a {\em
dynamically regular value} of $H$. If $a$ is dynamically regular, then
$H^{-1}(a)$ is mapped by $h$ into $H^{-1}((-\infty,a))$
\end{definition}

The existence of a complete Lyapunov function for any homeomorphism of a
compact metric space is shown in \cite{F}. For our purpose, the following proposition
is more convenient. The proof can be found in Appendix.

\begin{proposition} \label{p1}
For any $h\in\HH$, there is a $C^\infty$ complete Lyapunov
function $H$ of $h$.
\end{proposition}

Condition (3) above and the Sard theorem say that dynamically regular
and regular (in the usual sense) values are Lebesgue full.
For such a value $a$, $H^{-1}(a)$ is a 1 dimensional $C^\infty$ 
submanifold and $h$ maps $H^{-1}(a)$ into $H^{-1}(-\infty,a)$.

\medskip
\noindent
{\bf 3.2.} Let us construct an example of $C^\infty$ diffeomorphism
$h\in\HH$ which admits a chain transitive class $C_0$ such that the
rotation set $\rot(\tilde h,C_0)$ is a nontrivial interval.
We construct $h$ as an area preserving diffeomorphism of the annulus
$\AAA$,
which is so to call a ``winding horseshoe map''. See Figure 1.
Let us denote by $m$ the (infinite) measure on $\AAA$ given by
the area form $d\theta\wedge dt$. 

\begin{figure}
\unitlength 0.1in
\begin{picture}( 37.9000, 25.2000)( 14.0000,-29.1000)
%
{\color[named]{Black}{%
\special{pn 8}%
\special{pa 1400 390}%
\special{pa 1400 2910}%
\special{fp}%
}}%
%
{\color[named]{Black}{%
\special{pn 8}%
\special{pa 5190 390}%
\special{pa 5190 390}%
\special{fp}%
\special{pa 5100 390}%
\special{pa 5100 2890}%
\special{fp}%
}}%
%
{\color[named]{Black}{%
\special{pn 8}%
\special{pa 2410 670}%
\special{pa 4200 670}%
\special{pa 4200 2460}%
\special{pa 2410 2460}%
\special{pa 2410 670}%
\special{pa 4200 670}%
\special{fp}%
}}%
%
{\color[named]{Black}{%
\special{pn 8}%
\special{pa 2410 1050}%
\special{pa 4200 1050}%
\special{fp}%
}}%
%
{\color[named]{Black}{%
\special{pn 8}%
\special{pa 2420 2100}%
\special{pa 4190 2100}%
\special{fp}%
}}%
%
{\color[named]{Black}{%
\special{pn 8}%
\special{pa 2810 670}%
\special{pa 2810 2460}%
\special{dt 0.045}%
}}%
%
{\color[named]{Black}{%
\special{pn 8}%
\special{pa 3820 670}%
\special{pa 3820 2450}%
\special{dt 0.045}%
}}%
%
{\color[named]{Black}{%
\special{pn 8}%
\special{pa 4210 670}%
\special{pa 4242 674}%
\special{pa 4274 676}%
\special{pa 4370 688}%
\special{pa 4400 694}%
\special{pa 4432 702}%
\special{pa 4462 708}%
\special{pa 4582 748}%
\special{pa 4612 760}%
\special{pa 4640 774}%
\special{pa 4670 786}%
\special{pa 4698 800}%
\special{pa 4726 816}%
\special{pa 4754 830}%
\special{pa 4784 846}%
\special{pa 4840 878}%
\special{pa 4866 896}%
\special{pa 4894 914}%
\special{pa 4922 930}%
\special{pa 4950 948}%
\special{pa 4978 968}%
\special{pa 5004 986}%
\special{pa 5060 1022}%
\special{pa 5086 1042}%
\special{pa 5100 1050}%
\special{fp}%
}}%
%
{\color[named]{Black}{%
\special{pn 8}%
\special{pa 4210 1050}%
\special{pa 4274 1054}%
\special{pa 4370 1066}%
\special{pa 4400 1072}%
\special{pa 4460 1088}%
\special{pa 4520 1108}%
\special{pa 4550 1120}%
\special{pa 4578 1132}%
\special{pa 4606 1146}%
\special{pa 4636 1160}%
\special{pa 4664 1174}%
\special{pa 4748 1222}%
\special{pa 4774 1240}%
\special{pa 4830 1276}%
\special{pa 4856 1294}%
\special{pa 4884 1314}%
\special{pa 4936 1354}%
\special{pa 4964 1374}%
\special{pa 5016 1414}%
\special{pa 5044 1436}%
\special{pa 5096 1476}%
\special{pa 5100 1480}%
\special{fp}%
}}%
%
{\color[named]{Black}{%
\special{pn 8}%
\special{pa 2420 2090}%
\special{pa 2356 2086}%
\special{pa 2324 2082}%
\special{pa 2294 2078}%
\special{pa 2198 2066}%
\special{pa 2166 2060}%
\special{pa 2136 2054}%
\special{pa 2104 2048}%
\special{pa 2074 2040}%
\special{pa 2042 2034}%
\special{pa 2010 2026}%
\special{pa 1950 2010}%
\special{pa 1918 2002}%
\special{pa 1888 1994}%
\special{pa 1856 1984}%
\special{pa 1826 1976}%
\special{pa 1796 1966}%
\special{pa 1764 1956}%
\special{pa 1704 1936}%
\special{pa 1672 1926}%
\special{pa 1582 1896}%
\special{pa 1550 1884}%
\special{pa 1520 1874}%
\special{pa 1490 1862}%
\special{pa 1460 1852}%
\special{pa 1430 1840}%
\special{pa 1400 1830}%
\special{fp}%
}}%
%
{\color[named]{Black}{%
\special{pn 8}%
\special{pa 2400 2460}%
\special{pa 2336 2452}%
\special{pa 2306 2448}%
\special{pa 2274 2444}%
\special{pa 2242 2438}%
\special{pa 2210 2434}%
\special{pa 2180 2428}%
\special{pa 2116 2416}%
\special{pa 2086 2410}%
\special{pa 2054 2402}%
\special{pa 2022 2396}%
\special{pa 1992 2388}%
\special{pa 1960 2382}%
\special{pa 1930 2374}%
\special{pa 1898 2366}%
\special{pa 1868 2358}%
\special{pa 1836 2350}%
\special{pa 1806 2342}%
\special{pa 1774 2332}%
\special{pa 1744 2324}%
\special{pa 1712 2314}%
\special{pa 1682 2306}%
\special{pa 1652 2296}%
\special{pa 1620 2288}%
\special{pa 1590 2278}%
\special{pa 1558 2270}%
\special{pa 1498 2250}%
\special{pa 1466 2240}%
\special{pa 1436 2232}%
\special{pa 1400 2220}%
\special{fp}%
}}%
\put(22.8000,-5.6000){\makebox(0,0){$a$}}%
\put(43.8000,-25.8000){\makebox(0,0){$b$}}%
\put(33.1000,-15.5000){\makebox(0,0){$R$}}%
\put(46.5000,-9.8000){\makebox(0,0){$h(R)$}}%
\put(26.3000,-15.7000){\makebox(0,0){$R_0$}}%
\put(40.2000,-15.7000){\makebox(0,0){$R_1$}}%
\end{picture}%

\caption{}
\end{figure}

Choose a rectangle 
$R=[0,4^{-1}]\times [-4^{-1},0]$ in $\AAA=(\R/\Z)\times\R$. 
Strech $R$ horizontally
by 5 and contract vertically by $5^{-1}$, and embed the resultant
long and thin rectangle into $\AAA$ in a way to wind
the annulus $\AAA$. The precise conditions for a map $h:R\to\AAA$
is the following.

\begin{enumerate} 
\item $h$ is an $m$-preserving $C^\infty$ embedding.
\item Restricted to the subrectangle $R_0=[0,20^{-1}]\times[-4^{-1},0]$, \
$$h(\theta,t)=(5\theta,5^{-1}t).$$
\item Restricted to the subrectangle
      $R_1=[5^{-1},4^{-1}]\times[-4^{-1},0]$, 
$$h(\theta,t)=(5\theta-1,5^{-1}t-5^{-1}).$$
\item $h^{-1}(R)=R_0\cup R_1$.
\item $R\cup h(R)$ separates both ends of $\AAA$.
\end{enumerate}

Notice that the points $a=(0,0)$ and $b=(4^{-1},-4^{-1})$ are the (only)
fixed points of $h$. 
Extend $h$ to a $C^\infty$ diffeomorphism $h_0$ of $\AAA$ so as to
satisfy the following conditions.

\smallskip
(6) On $S^1\times ((-\infty,-10]\cup[10,\infty))$,
    $h_0(\theta,t)=(\theta,t-1)$.

\smallskip
The  measure $(h_0)_*m$ coincides with $m$ on $h_0(R)$, since
$h_0$ is $m$-preserving on $R$, and likewise on
$h_0(S^1\times((-\infty,-10]\cup[10,\infty)))$. 
Now by Moser's lemma (\cite{HZ},
p.16),
there is a $C^\infty$ diffeomorphsim $h_1$ on $\AAA$
such that $(h_1)_*((h_0)_*(m))=m$ which is the identity
on $h_0(R\cup(S^1\times(-\infty,-10]\cup[10,\infty)))$.

Now the composite $h=h_1\circ h_0$ is $m$-preserving. Let us show that
$h$ satisfies the condition raised in the beginning of {\bf 3.2}.
First of all clearly $h$ satisfies condition (1) $\sim$ (3) of
Definition
\ref{HH}. Moreover since $h$ is $m$-preserving, it cannot admit a
invariant continuum separating both ends of
$\AAA$. Therefore it satisfies (4) also.

Choose a lift $\tilde h$ of $h$ so that each point of $\pi^{-1}(a)$ is 
fixed by $\tilde h$. Then we have $\rot(\tilde h,a)=0$ and $\rot(\tilde
h,b)=1$.

Finally we have $a\sim b$, since the stable manifold of $a$ intersects
the unstable manifold of $b$, and the unstable manifold of $a$
intersects the stable manifold of $b$. 
Therefore the chain transitive class $C_0$ of $a$ and $b$ satisfies
$[0,1]\subset\rot(\tilde h,C_0)$.

\medskip
\noindent
{\bf 3.3.} Here we shall show Theorem \ref{main2} by a rather
lengthy argument. We consider
$h\in\HH$ to be a homeomorphism of the annulus $\AAA$.
Denote the generator of the covering transormations by
$T:\tilde\AAA\to\tilde\AAA$: \ $T(\theta,t)=(\theta+1,t)$.
First of all we have the following fundamental lemma.

\begin{lemma} \label{l3.3}
Suppose $h^q(z)=z$  ($q\in\N$, $z\in\AAA$) and let $p\in\Z$.
(We do not assume $(p,q)=1$.)
Then the following
conditions are equivalent.
\begin{enumerate}
\item
$\rot(\tilde h,z)=p/q$.
\item
$\tilde h^q(\tilde z)=T^p(\tilde z)$ for a lift $\tilde z$ of $z$.
\end{enumerate}
\end{lemma}

Notice that condition (2) is independent of the choice of the lift
$\tilde z$. This is because $\tilde h$ commutes with $T$.

\bd
Recall that $\pi:\tilde\AAA\to\AAA$ is the universal covering map
and $\Pi_1:\tilde \AAA\to\R$ is the canonical projection onto the first
factor.
Define a function $\varphi:\AAA\to\R$ by
$$\pi\circ\varphi=\Pi_1\circ \tilde h-\Pi_1.$$
Denote the average of the Dirac masses along the orbit of $z$ by $\mu$,
that is,
$$\mu=q^{-1}(\delta_{z}+\delta_{h(z)}+\cdots+\delta_{h^{q-1}(z)}).$$
Notice that for any $i\in\N$,
\begin{equation}\label{el3.3.1}
\langle \delta_{h^i}(z),\varphi\rangle=
\langle\pi_*\delta_{\tilde h^i(\tilde z)},\varphi\rangle
=\langle \delta_{\tilde h^i(\tilde z)},\varphi\circ\pi\rangle
\end{equation}
$$
=\langle\delta_{\tilde h^i(\tilde z)},\Pi_1\circ \tilde h-\Pi_1\rangle
=\Pi_1(\tilde h^{i+1}(\tilde z))-\Pi_1(\tilde h^{i}(\tilde z)).$$

Now assume (1): 
$
\rot(\tilde h, z)=\langle \mu,\,\,\varphi\rangle=p/q
$.
Then
we have by (\ref{el3.3.1})
$$
\Pi_1(h^q(\tilde z))-\Pi_1(\tilde z)=p.
$$
That is,
\begin{equation}\label{el3.3.2}
\Pi_1(h^q(\tilde z))=\Pi_1(T^p(\tilde z)).
\end{equation} 
On the other hand, since $h^q(z)=z$, we have
\begin{equation}\label{el3.3.3}
\tilde h^q(\tilde z)=T^j(\tilde z)
\end{equation}
for some $j\in\Z$.
Now (\ref{el3.3.2}) and (\ref{el3.3.3}) imply that $j=p$.
We obtain condition (2). The converse can be shown by a reversed
argument. \qed

\medskip
Let us begin the proof of Theorem \ref{main2}.
Let $h\in\HH$ and $C_0$ a chain transitive class of $h$.
Assume $x_\nu\in C_0$ are periodic points  such that
$\rot(\tilde h,x_\nu)=\alpha_\nu$ ($\nu=1,2$) and let $\alpha$ be
a rational number in $[\alpha_1,\alpha_2]$. If $\alpha_1=\alpha_2$,
there is nothing to prove. So assume $\alpha_1<\alpha<\alpha_2$.
Then
it is possible to choose a number $q\in\N$ such that

(a) the rational numbers $\alpha_\nu$ and $\alpha$ are written as
$$
\alpha_\nu=p_\nu/q,\ \  (\nu=1,2),\ \ \alpha=p/q, \ \
p_1<p<p_2,\ \ \mbox{ and}$$

(b) the periodic points $x_\nu$  satisfies
$h^q(x_\nu)=x_\nu$.

By Lemma \ref{l3.3}, lifts $\tilde x_\nu$ of $x_\nu$ satisfy
$\tilde h^q(\tilde x_\nu)=T^{p_\nu}(\tilde x_\nu)$.
Our purpose is to show the existence of a periodic point $ x\in C_0$
of period $q$ such that $\rot(\tilde h,x)=p/q$, that is, whose lifts
 $\tilde x$ satisfy 
$\tilde h^q(\tilde x)=T^{p}(\tilde x)$.

However a simultaneous proof for all $p\in(p_1,p_2)$ has an elementary
number theoretic difficulty.
We shall avoid it by employing
an induction on $p-p_1$. 
Namely we first show only for $p=p_1+1$. Then the newly obtained periodic
points
can serve as an assumption for the next step $p=p_1+2$.
This way,
Theorem \ref{main2} reduces to the following.

\begin{proposition} \label{p2}
Let $q>0$ and $p_1+1<p_2$ and let $C_0$ be a chain transitive class of
$h\in\HH$. Assume 
there are points $x_\nu\in C_0$ with a lift $\tilde x_\nu$ such that 
$\tilde h^q(\tilde x_\nu)=T^{p_\nu}(\tilde x_\nu)$ ($\nu=1,2$). 
Then there is 
a point $x\in C_0$ with a lift $\tilde x$
 such that $\tilde h^q(\tilde x)=T^{p_1+1}(\tilde x)$.
\end{proposition}

Now Proposition \ref{p2} itself reduces to the following.

\begin{proposition}\label{p2a}
Let $q>0$ and $p_1+1<p_2$ and let $C_0$ be a chain transitive class of
$h\in\HH$. Assume 
there are points $x_\nu\in C_0$ with a lift $\tilde x_\nu$ such that 
$\tilde h^q(\tilde x_\nu)=T^{p_\nu}(\tilde x_\nu)$ ($\nu=1,2$). 
Let $H$ be a $C^\infty$ complete Lyapunov function such that $H(C_0)=0$,
and let $-a'<0$ and $a>0$ be  dynamically regular 
and regular values of $H$.
Then there is 
a point $x$ in the subsurface $H^{-1}([-a',a])$ with a lift $\tilde x$
 such that $\tilde h^q(\tilde x)=T^{p_1+1}(\tilde x)$.
\end{proposition}

Postponing the proof, we shall show the reduction.
Thanks to the Sard theorem and condition (3) of Definition
\ref{d2}, one can find dynamically regular 
and regular values $-a'<0<a$ 
of
$H$ as close to 0 as we want. 
Let
$$F_0=\pi({\rm Fix}(\tilde h^q\circ T^{-p_1-1})).$$
Then Proposition \ref{p2a} says that 
$F_0\cap H^{-1}([-a',a])\neq\emptyset$
for any such values. 
By the compactness of $F_0$, this implis that
$F_0\cap H^{-1}(0)\neq\emptyset$.
On the other hand, since $H$ is a complete Lyapunov function,
 $C\cap H^{-1}(0)=C_0$, showing $F_0\cap C_0\neq\emptyset$,
as is required in Proposition \ref{p2}.

\medskip

The rest of this paragraph is devoted to the proof of Proposition
\ref{p2a}.
The subsurface $H^{-1}([-a',a])$ admits
a single distinguished connected component $X$ which is
homotopically nontrivial in $\AAA$. In fact, if there were more than
one such components, then in the complement, one could find a 
forward invariant compact subannulus. The intersection of its forward
images would be a $h$-invariant continuum separating $U_\infty$ and
$U_{-\infty}$, contradicting condition (4) of Definition \ref{HH}.

Let us consider the upper boundary $H^{-1}(a)\cap X$ of $X$.
It has a unique homotopically nontrivial component $\partial A^+$.
The curve $\partial A^+$ bounds an infinite annulus $A^+$
on the opposite side of $X$. The intersection of $\Int(A^+)$ with the level
$H^{-1}(a)$
consists of finitely many circles $\partial D_i^+$. They bound
discs $D_i^+$ in $A^+$. See Figure 2.

The components of $H^{-1}(a)\cap X$ other than $\partial A^+$
are denoted by $\partial E_k^+$. They are finite
in number and bound discs $E_k^+$ in $\AAA$.

Likewise we define an annulus $A^-$, discs $D^-_j$ and $E^-_l$ by
considering the lower boundary $H^{-1}(-a')\cap X$ of $X$.
Then we have
$$
\AAA=X\cup A^-\cup A^+\cup\bigcup_k E_k^+\cup\bigcup_l E_l^-.$$

\begin{figure}
\unitlength 0.1in
\begin{picture}( 32.1900, 34.2200)( 14.0200,-37.9900)
%
{\color[named]{Black}{%
\special{pn 8}%
\special{pa 2460 430}%
\special{pa 2462 462}%
\special{pa 2462 526}%
\special{pa 2464 558}%
\special{pa 2464 622}%
\special{pa 2466 654}%
\special{pa 2466 718}%
\special{pa 2468 750}%
\special{pa 2468 846}%
\special{pa 2470 878}%
\special{pa 2470 1198}%
\special{pa 2468 1230}%
\special{pa 2468 1294}%
\special{pa 2466 1326}%
\special{pa 2466 1358}%
\special{pa 2464 1390}%
\special{pa 2464 1422}%
\special{pa 2460 1486}%
\special{pa 2460 1518}%
\special{pa 2456 1582}%
\special{pa 2452 1614}%
\special{pa 2450 1646}%
\special{pa 2446 1678}%
\special{pa 2440 1710}%
\special{pa 2434 1740}%
\special{pa 2428 1772}%
\special{pa 2420 1804}%
\special{pa 2410 1834}%
\special{pa 2400 1866}%
\special{pa 2388 1896}%
\special{pa 2374 1926}%
\special{pa 2358 1956}%
\special{pa 2340 1982}%
\special{pa 2318 2006}%
\special{pa 2292 2024}%
\special{pa 2262 2038}%
\special{pa 2228 2048}%
\special{pa 2194 2050}%
\special{pa 2160 2048}%
\special{pa 2128 2040}%
\special{pa 2098 2026}%
\special{pa 2076 2006}%
\special{pa 2058 1982}%
\special{pa 2044 1954}%
\special{pa 2034 1922}%
\special{pa 2024 1888}%
\special{pa 2016 1854}%
\special{pa 2008 1818}%
\special{pa 1996 1750}%
\special{pa 1994 1716}%
\special{pa 1992 1684}%
\special{pa 1994 1652}%
\special{pa 2000 1622}%
\special{pa 2008 1594}%
\special{pa 2022 1568}%
\special{pa 2040 1544}%
\special{pa 2062 1520}%
\special{pa 2088 1500}%
\special{pa 2112 1480}%
\special{pa 2140 1458}%
\special{pa 2164 1438}%
\special{pa 2188 1414}%
\special{pa 2208 1390}%
\special{pa 2224 1364}%
\special{pa 2238 1338}%
\special{pa 2248 1308}%
\special{pa 2258 1276}%
\special{pa 2264 1244}%
\special{pa 2268 1210}%
\special{pa 2272 1138}%
\special{pa 2272 1104}%
\special{pa 2268 1070}%
\special{pa 2260 1040}%
\special{pa 2244 1014}%
\special{pa 2222 998}%
\special{pa 2196 986}%
\special{pa 2162 978}%
\special{pa 2128 974}%
\special{pa 2094 972}%
\special{pa 2058 968}%
\special{pa 2026 966}%
\special{pa 1994 966}%
\special{pa 1964 968}%
\special{pa 1932 976}%
\special{pa 1902 988}%
\special{pa 1874 1006}%
\special{pa 1850 1028}%
\special{pa 1832 1054}%
\special{pa 1818 1082}%
\special{pa 1810 1112}%
\special{pa 1806 1146}%
\special{pa 1802 1178}%
\special{pa 1800 1212}%
\special{pa 1800 1244}%
\special{pa 1798 1278}%
\special{pa 1796 1310}%
\special{pa 1794 1340}%
\special{pa 1792 1372}%
\special{pa 1792 1404}%
\special{pa 1788 1468}%
\special{pa 1788 1500}%
\special{pa 1790 1532}%
\special{pa 1790 1564}%
\special{pa 1798 1628}%
\special{pa 1808 1788}%
\special{pa 1812 1818}%
\special{pa 1816 1850}%
\special{pa 1828 1914}%
\special{pa 1836 1944}%
\special{pa 1846 1976}%
\special{pa 1866 2036}%
\special{pa 1902 2126}%
\special{pa 1916 2154}%
\special{pa 1932 2182}%
\special{pa 1946 2210}%
\special{pa 1962 2238}%
\special{pa 1980 2266}%
\special{pa 1994 2294}%
\special{pa 2026 2350}%
\special{pa 2062 2440}%
\special{pa 2070 2470}%
\special{pa 2078 2502}%
\special{pa 2082 2532}%
\special{pa 2088 2628}%
\special{pa 2092 2660}%
\special{pa 2096 2694}%
\special{pa 2102 2726}%
\special{pa 2106 2758}%
\special{pa 2106 2790}%
\special{pa 2104 2820}%
\special{pa 2098 2852}%
\special{pa 2084 2882}%
\special{pa 2068 2912}%
\special{pa 2046 2938}%
\special{pa 2020 2960}%
\special{pa 1992 2978}%
\special{pa 1962 2990}%
\special{pa 1930 2996}%
\special{pa 1898 2996}%
\special{pa 1864 2992}%
\special{pa 1834 2982}%
\special{pa 1804 2968}%
\special{pa 1778 2948}%
\special{pa 1754 2926}%
\special{pa 1734 2902}%
\special{pa 1718 2876}%
\special{pa 1704 2846}%
\special{pa 1692 2816}%
\special{pa 1682 2786}%
\special{pa 1658 2690}%
\special{pa 1650 2660}%
\special{pa 1642 2628}%
\special{pa 1636 2596}%
\special{pa 1628 2566}%
\special{pa 1624 2534}%
\special{pa 1618 2502}%
\special{pa 1614 2472}%
\special{pa 1610 2440}%
\special{pa 1606 2376}%
\special{pa 1606 2248}%
\special{pa 1608 2216}%
\special{pa 1608 2184}%
\special{pa 1610 2152}%
\special{pa 1610 1992}%
\special{pa 1608 1960}%
\special{pa 1608 1928}%
\special{pa 1602 1832}%
\special{pa 1602 1800}%
\special{pa 1594 1672}%
\special{pa 1594 1640}%
\special{pa 1592 1608}%
\special{pa 1592 1576}%
\special{pa 1590 1544}%
\special{pa 1590 1416}%
\special{pa 1600 1256}%
\special{pa 1608 1192}%
\special{pa 1610 1160}%
\special{pa 1610 1032}%
\special{pa 1614 1000}%
\special{pa 1620 970}%
\special{pa 1634 940}%
\special{pa 1650 914}%
\special{pa 1670 888}%
\special{pa 1694 864}%
\special{pa 1720 844}%
\special{pa 1748 826}%
\special{pa 1776 810}%
\special{pa 1806 794}%
\special{pa 1834 780}%
\special{pa 1862 764}%
\special{pa 1888 748}%
\special{pa 1914 728}%
\special{pa 1936 704}%
\special{pa 1954 676}%
\special{pa 1970 646}%
\special{pa 1978 614}%
\special{pa 1982 582}%
\special{pa 1976 550}%
\special{pa 1964 520}%
\special{pa 1948 492}%
\special{pa 1924 468}%
\special{pa 1898 444}%
\special{pa 1870 426}%
\special{pa 1840 410}%
\special{pa 1810 400}%
\special{pa 1780 392}%
\special{pa 1750 388}%
\special{pa 1718 384}%
\special{pa 1650 380}%
\special{pa 1614 378}%
\special{pa 1580 378}%
\special{pa 1548 386}%
\special{pa 1524 398}%
\special{pa 1504 420}%
\special{pa 1490 448}%
\special{pa 1480 480}%
\special{pa 1472 514}%
\special{pa 1464 546}%
\special{pa 1454 578}%
\special{pa 1434 638}%
\special{pa 1426 668}%
\special{pa 1414 728}%
\special{pa 1410 758}%
\special{pa 1406 790}%
\special{pa 1404 820}%
\special{pa 1402 852}%
\special{pa 1402 916}%
\special{pa 1404 948}%
\special{pa 1404 980}%
\special{pa 1406 1014}%
\special{pa 1408 1046}%
\special{pa 1412 1078}%
\special{pa 1414 1112}%
\special{pa 1416 1144}%
\special{pa 1420 1176}%
\special{pa 1422 1210}%
\special{pa 1424 1242}%
\special{pa 1428 1276}%
\special{pa 1430 1308}%
\special{pa 1430 1340}%
\special{pa 1432 1372}%
\special{pa 1432 1470}%
\special{pa 1430 1502}%
\special{pa 1430 1534}%
\special{pa 1428 1566}%
\special{pa 1428 1598}%
\special{pa 1424 1662}%
\special{pa 1424 1692}%
\special{pa 1422 1724}%
\special{pa 1422 1756}%
\special{pa 1420 1788}%
\special{pa 1420 2012}%
\special{pa 1422 2044}%
\special{pa 1422 2076}%
\special{pa 1424 2108}%
\special{pa 1424 2140}%
\special{pa 1426 2172}%
\special{pa 1426 2204}%
\special{pa 1428 2236}%
\special{pa 1428 2300}%
\special{pa 1430 2332}%
\special{pa 1430 2396}%
\special{pa 1432 2428}%
\special{pa 1432 2750}%
\special{pa 1430 2780}%
\special{pa 1430 2972}%
\special{pa 1434 3036}%
\special{pa 1438 3068}%
\special{pa 1444 3100}%
\special{pa 1454 3132}%
\special{pa 1466 3162}%
\special{pa 1482 3192}%
\special{pa 1502 3220}%
\special{pa 1524 3240}%
\special{pa 1552 3256}%
\special{pa 1582 3264}%
\special{pa 1614 3268}%
\special{pa 1648 3270}%
\special{pa 1712 3274}%
\special{pa 1746 3276}%
\special{pa 1778 3280}%
\special{pa 1810 3280}%
\special{pa 1840 3282}%
\special{pa 1904 3282}%
\special{pa 1936 3280}%
\special{pa 1966 3278}%
\special{pa 1998 3274}%
\special{pa 2030 3272}%
\special{pa 2064 3268}%
\special{pa 2096 3268}%
\special{pa 2128 3270}%
\special{pa 2160 3274}%
\special{pa 2192 3282}%
\special{pa 2220 3296}%
\special{pa 2248 3314}%
\special{pa 2272 3334}%
\special{pa 2296 3356}%
\special{pa 2336 3400}%
\special{pa 2352 3424}%
\special{pa 2366 3450}%
\special{pa 2378 3480}%
\special{pa 2386 3516}%
\special{pa 2390 3558}%
\special{pa 2392 3606}%
\special{pa 2388 3706}%
\special{pa 2384 3750}%
\special{pa 2382 3782}%
\special{pa 2384 3798}%
\special{pa 2388 3796}%
\special{pa 2398 3770}%
\special{pa 2400 3760}%
\special{fp}%
}}%
%
{\color[named]{Black}{%
\special{pn 8}%
\special{pa 3380 420}%
\special{pa 3378 452}%
\special{pa 3374 484}%
\special{pa 3372 518}%
\special{pa 3370 550}%
\special{pa 3370 614}%
\special{pa 3374 678}%
\special{pa 3378 708}%
\special{pa 3394 772}%
\special{pa 3406 804}%
\special{pa 3420 832}%
\special{pa 3438 858}%
\special{pa 3462 880}%
\special{pa 3490 896}%
\special{pa 3518 910}%
\special{pa 3548 920}%
\special{pa 3578 928}%
\special{pa 3642 940}%
\special{pa 3678 944}%
\special{pa 3714 952}%
\special{pa 3748 960}%
\special{pa 3780 972}%
\special{pa 3806 988}%
\special{pa 3826 1010}%
\special{pa 3836 1038}%
\special{pa 3840 1070}%
\special{pa 3842 1104}%
\special{pa 3846 1136}%
\special{pa 3848 1166}%
\special{pa 3854 1198}%
\special{pa 3858 1230}%
\special{pa 3864 1262}%
\special{pa 3872 1326}%
\special{pa 3876 1356}%
\special{pa 3880 1388}%
\special{pa 3886 1420}%
\special{pa 3890 1452}%
\special{pa 3890 1516}%
\special{pa 3892 1548}%
\special{pa 3898 1578}%
\special{pa 3916 1674}%
\special{pa 3926 1704}%
\special{pa 3942 1732}%
\special{pa 3948 1764}%
\special{pa 3958 1794}%
\special{pa 3978 1820}%
\special{pa 4012 1838}%
\special{pa 4048 1848}%
\special{pa 4080 1846}%
\special{pa 4100 1830}%
\special{pa 4106 1804}%
\special{pa 4102 1770}%
\special{pa 4094 1732}%
\special{pa 4088 1692}%
\special{pa 4082 1656}%
\special{pa 4078 1622}%
\special{pa 4076 1590}%
\special{pa 4074 1560}%
\special{pa 4074 1468}%
\special{pa 4070 1408}%
\special{pa 4068 1376}%
\special{pa 4060 1312}%
\special{pa 4054 1216}%
\special{pa 4048 1186}%
\special{pa 4040 1154}%
\special{pa 4032 1124}%
\special{pa 4028 1092}%
\special{pa 4028 1060}%
\special{pa 4030 1028}%
\special{pa 4030 964}%
\special{pa 4026 900}%
\special{pa 4028 868}%
\special{pa 4034 836}%
\special{pa 4042 804}%
\special{pa 4056 772}%
\special{pa 4074 744}%
\special{pa 4096 720}%
\special{pa 4120 700}%
\special{pa 4150 686}%
\special{pa 4180 678}%
\special{pa 4214 674}%
\special{pa 4248 676}%
\special{pa 4282 682}%
\special{pa 4316 692}%
\special{pa 4346 706}%
\special{pa 4376 724}%
\special{pa 4400 744}%
\special{pa 4422 768}%
\special{pa 4440 794}%
\special{pa 4456 820}%
\special{pa 4470 850}%
\special{pa 4482 882}%
\special{pa 4492 914}%
\special{pa 4500 946}%
\special{pa 4518 1048}%
\special{pa 4520 1080}%
\special{pa 4520 1112}%
\special{pa 4516 1142}%
\special{pa 4506 1172}%
\special{pa 4490 1198}%
\special{pa 4472 1226}%
\special{pa 4432 1278}%
\special{pa 4414 1304}%
\special{pa 4398 1332}%
\special{pa 4384 1360}%
\special{pa 4352 1416}%
\special{pa 4338 1444}%
\special{pa 4320 1472}%
\special{pa 4292 1528}%
\special{pa 4282 1558}%
\special{pa 4272 1590}%
\special{pa 4264 1622}%
\special{pa 4260 1652}%
\special{pa 4254 1684}%
\special{pa 4246 1812}%
\special{pa 4242 1844}%
\special{pa 4236 1876}%
\special{pa 4226 1908}%
\special{pa 4212 1936}%
\special{pa 4190 1960}%
\special{pa 4164 1978}%
\special{pa 4134 1992}%
\special{pa 4102 2000}%
\special{pa 4068 2002}%
\special{pa 4038 2000}%
\special{pa 4006 2000}%
\special{pa 3972 2002}%
\special{pa 3940 2004}%
\special{pa 3910 2000}%
\special{pa 3878 1990}%
\special{pa 3850 1974}%
\special{pa 3826 1952}%
\special{pa 3808 1926}%
\special{pa 3796 1898}%
\special{pa 3786 1868}%
\special{pa 3780 1836}%
\special{pa 3776 1802}%
\special{pa 3770 1770}%
\special{pa 3762 1740}%
\special{pa 3754 1708}%
\special{pa 3742 1678}%
\special{pa 3712 1588}%
\special{pa 3704 1556}%
\special{pa 3700 1524}%
\special{pa 3700 1492}%
\special{pa 3702 1458}%
\special{pa 3706 1424}%
\special{pa 3706 1392}%
\special{pa 3690 1368}%
\special{pa 3664 1348}%
\special{pa 3634 1332}%
\special{pa 3604 1322}%
\special{pa 3572 1320}%
\special{pa 3538 1328}%
\special{pa 3510 1346}%
\special{pa 3498 1374}%
\special{pa 3498 1406}%
\special{pa 3502 1440}%
\special{pa 3500 1472}%
\special{pa 3500 1536}%
\special{pa 3504 1600}%
\special{pa 3508 1632}%
\special{pa 3510 1664}%
\special{pa 3514 1694}%
\special{pa 3518 1726}%
\special{pa 3520 1758}%
\special{pa 3520 1822}%
\special{pa 3516 1886}%
\special{pa 3514 1920}%
\special{pa 3514 1952}%
\special{pa 3516 1982}%
\special{pa 3528 2046}%
\special{pa 3530 2078}%
\special{pa 3526 2110}%
\special{pa 3520 2142}%
\special{pa 3524 2172}%
\special{pa 3534 2202}%
\special{pa 3546 2232}%
\special{pa 3560 2262}%
\special{pa 3572 2292}%
\special{pa 3592 2352}%
\special{pa 3600 2384}%
\special{pa 3606 2414}%
\special{pa 3614 2448}%
\special{pa 3622 2480}%
\special{pa 3632 2510}%
\special{pa 3646 2540}%
\special{pa 3664 2564}%
\special{pa 3688 2584}%
\special{pa 3718 2600}%
\special{pa 3750 2610}%
\special{pa 3784 2614}%
\special{pa 3818 2612}%
\special{pa 3850 2606}%
\special{pa 3880 2592}%
\special{pa 3906 2576}%
\special{pa 3930 2554}%
\special{pa 3952 2530}%
\special{pa 3972 2502}%
\special{pa 3988 2472}%
\special{pa 4004 2440}%
\special{pa 4018 2408}%
\special{pa 4030 2376}%
\special{pa 4040 2344}%
\special{pa 4048 2312}%
\special{pa 4058 2280}%
\special{pa 4070 2250}%
\special{pa 4082 2222}%
\special{pa 4096 2194}%
\special{pa 4114 2168}%
\special{pa 4136 2146}%
\special{pa 4162 2124}%
\special{pa 4192 2106}%
\special{pa 4222 2092}%
\special{pa 4254 2082}%
\special{pa 4288 2078}%
\special{pa 4318 2078}%
\special{pa 4348 2086}%
\special{pa 4376 2100}%
\special{pa 4400 2118}%
\special{pa 4424 2142}%
\special{pa 4444 2168}%
\special{pa 4464 2198}%
\special{pa 4482 2226}%
\special{pa 4500 2256}%
\special{pa 4516 2286}%
\special{pa 4544 2346}%
\special{pa 4556 2376}%
\special{pa 4568 2404}%
\special{pa 4588 2464}%
\special{pa 4596 2494}%
\special{pa 4602 2524}%
\special{pa 4614 2588}%
\special{pa 4616 2620}%
\special{pa 4620 2652}%
\special{pa 4620 2684}%
\special{pa 4622 2718}%
\special{pa 4620 2754}%
\special{pa 4620 2788}%
\special{pa 4612 2860}%
\special{pa 4606 2896}%
\special{pa 4598 2930}%
\special{pa 4588 2962}%
\special{pa 4574 2992}%
\special{pa 4560 3020}%
\special{pa 4542 3042}%
\special{pa 4522 3062}%
\special{pa 4498 3078}%
\special{pa 4472 3088}%
\special{pa 4442 3094}%
\special{pa 4410 3098}%
\special{pa 4378 3098}%
\special{pa 4342 3096}%
\special{pa 4308 3092}%
\special{pa 4272 3086}%
\special{pa 4204 3074}%
\special{pa 4170 3066}%
\special{pa 4138 3060}%
\special{pa 4074 3044}%
\special{pa 4042 3034}%
\special{pa 4012 3026}%
\special{pa 3982 3016}%
\special{pa 3954 3004}%
\special{pa 3924 2992}%
\special{pa 3868 2968}%
\special{pa 3840 2952}%
\special{pa 3812 2938}%
\special{pa 3786 2922}%
\special{pa 3734 2886}%
\special{pa 3682 2846}%
\special{pa 3658 2826}%
\special{pa 3610 2782}%
\special{pa 3566 2734}%
\special{pa 3546 2708}%
\special{pa 3526 2684}%
\special{pa 3508 2656}%
\special{pa 3490 2630}%
\special{pa 3476 2602}%
\special{pa 3462 2572}%
\special{pa 3450 2544}%
\special{pa 3440 2514}%
\special{pa 3432 2482}%
\special{pa 3420 2418}%
\special{pa 3416 2386}%
\special{pa 3414 2352}%
\special{pa 3408 2320}%
\special{pa 3400 2290}%
\special{pa 3386 2262}%
\special{pa 3370 2236}%
\special{pa 3354 2206}%
\special{pa 3342 2172}%
\special{pa 3328 2144}%
\special{pa 3304 2128}%
\special{pa 3270 2128}%
\special{pa 3240 2142}%
\special{pa 3222 2166}%
\special{pa 3212 2198}%
\special{pa 3210 2234}%
\special{pa 3210 2360}%
\special{pa 3212 2392}%
\special{pa 3212 2424}%
\special{pa 3210 2456}%
\special{pa 3210 2488}%
\special{pa 3208 2522}%
\special{pa 3208 2586}%
\special{pa 3212 2616}%
\special{pa 3216 2648}%
\special{pa 3224 2680}%
\special{pa 3234 2710}%
\special{pa 3258 2770}%
\special{pa 3290 2826}%
\special{pa 3308 2854}%
\special{pa 3344 2906}%
\special{pa 3364 2932}%
\special{pa 3418 3010}%
\special{pa 3436 3038}%
\special{pa 3452 3064}%
\special{pa 3466 3092}%
\special{pa 3480 3122}%
\special{pa 3494 3150}%
\special{pa 3504 3180}%
\special{pa 3516 3210}%
\special{pa 3524 3240}%
\special{pa 3534 3270}%
\special{pa 3542 3300}%
\special{pa 3554 3364}%
\special{pa 3560 3394}%
\special{pa 3568 3458}%
\special{pa 3572 3492}%
\special{pa 3576 3524}%
\special{pa 3578 3556}%
\special{pa 3582 3590}%
\special{pa 3586 3654}%
\special{pa 3590 3722}%
\special{pa 3590 3740}%
\special{fp}%
}}%
%
{\color[named]{Black}{%
\special{pn 8}%
\special{pa 1800 1670}%
\special{pa 2010 1670}%
\special{fp}%
}}%
%
{\color[named]{Black}{%
\special{pn 8}%
\special{pa 2440 1670}%
\special{pa 3510 1670}%
\special{fp}%
}}%
%
{\color[named]{Black}{%
\special{pn 8}%
\special{pa 3740 1670}%
\special{pa 3920 1670}%
\special{fp}%
}}%
%
{\color[named]{Black}{%
\special{pn 8}%
\special{pa 1420 2380}%
\special{pa 1620 2380}%
\special{fp}%
}}%
%
{\color[named]{Black}{%
\special{pn 8}%
\special{pa 2030 2380}%
\special{pa 3210 2380}%
\special{fp}%
}}%
%
{\color[named]{Black}{%
\special{pn 8}%
\special{pa 3420 2380}%
\special{pa 3600 2380}%
\special{fp}%
}}%
\put(29.7000,-12.2000){\makebox(0,0){$A^+$}}%
\put(27.1000,-20.6000){\makebox(0,0){$X$}}%
\put(42.4000,-11.0000){\makebox(0,0){$D_i^+$}}%
\put(42.6000,-26.7000){\makebox(0,0){$E_l^-$}}%
\put(16.4000,-6.2000){\makebox(0,0){$D_j^-$}}%
\put(19.8000,-12.8000){\makebox(0,0){$E^+_k$}}%
\put(26.1000,-28.6000){\makebox(0,0){$A^-$}}%
\end{picture}%

\caption{}
\end{figure}

Let us study how family of the discs $\DD^+=\{D_i^+\}$  are mapped by $h$,
and show the following. Denote $\abs{\DD^+}=\bigcup_iD_i^+$.

\begin{proposition} \label{l6}
The chain transitive class $C_0$ 
is disjoint from $A^+\cup A^-$.
\end{proposition}

This is not obvious since discs $D^+_i\subset A^+$ may intersect $H^{-1}(0)$.
Our overall strategy after having shown Proposition \ref{l6} is
to replace $h$ by a homeomorphism which has no periodic points in
$A^+\cup A^-$, and seek for  periodic points
in the rest of $\AAA$.
The argument will be divided into two cases. In the first case
we employ a topological method, while in the second a dynamical. 

To establish Proposition \ref{l6}, we must prepare some lemmas.

\begin{lemma} \label{l5}
For any small $\varepsilon>0$, an $\varepsilon$-chain joining two points in
$C_0$ is contained in $H^{-1}((-a',a))$.
\end{lemma}

\bd Notice that $-a'$ is
a dynamically regular value (Definiton \ref{d2}) and therefore $h$ maps
the level $H^{-1}(-a')$ below itself. 
Therefore there is $\varepsilon_0>0$ such that the $\varepsilon_0$-neighbourhood of
any point in $H^{-1}((-\infty,-a'])$ is mapped by $h$
into $H^{-1}((-\infty,-a'))$.
If we choose $\varepsilon<\varepsilon_0$
and if the $\varepsilon$-chain joining two points of $C_0$
falls into $H^{-1}((-\infty,-a'])$,
then the rest of the chain cannot escape $H^{-1}((-\infty,-a'))$
forever.  A contradiction.
The opposite case of falling into $H^{-1}([a,\infty))$
can be dealt with similarly by considering  $h^{-1}$
and the reversed chain. \qed

\medskip
Let $B^+=A^+\setminus \abs{\DD^+}$. Then we have 
$h^{-1}(B^+)\subset B^+$. In fact, a point $z\in B^+$ is
characterized by the existence of a  path in
$H^{-1}([a,\infty))$
starting at
$z$ and ending at a point in $\partial A^+$ without passing $H^{-1}(a)$
in the middle. This property is inherited to $h^{-1}(z)$ since 
there is a path from $h^{-1}(\partial A^+)$ to $\partial A^+$ 
which does not
pass $H^{-1}(a)$ in the middle.

The above inclusion
 implies that any disc $D^+_i\in\DD^+$ is mapped by $h$ to the
complement of $B^+$,
either into $\Int(D^+_{i'})$ for some $D^+_{i'}\in\DD^+$ or into
$\AAA\setminus A^+$. 
Notice that $\AAA\setminus A^+$ is forward
invariant by $h$. 

Let us call a sequence in $\DD^+$
$$\PP=\{D^+_{i_0},D^+_{i_1},\ldots,D^+_{i_n}\}$$
a {\em cycle of discs}, if $h(D^+_{i_{j-1}})\subset \Int(D^+_{i_{j}})$
($0< j\leq n$) and $D^+_{i_n}=D^+_{i_0}$.
Denote $\abs{\PP}=\bigcup_{j=0}^{n-1}D^+_{i_{j}}$.

\begin{lemma}\label{lnew}
If $C_0\cap A^+\neq\emptyset$, then there is a cycle of
discs $\PP$ in $\DD^+$ such that $C_0\cap\abs{\PP}\neq\emptyset$.
\end{lemma}

\bd
One can show as in the proof of Lemma \ref{l5} that for any small
$\varepsilon>0$, there is no $\varepsilon$-chain from a point
in $\AAA\setminus A^+$ to a point in $A^+$, since 
$h(\Cl(\AAA\setminus A^+))\subset \AAA\setminus A^+$.

Notice that $C_0$ is $h$-invariant.
Now if $D^+_i\cap C_0\neq\emptyset$ for some $D^+_i\in\DD^+$, $D^+_i$ cannot
be mapped into $\AAA\setminus A^+$ by a positive iterate of $h$.
Then the iterated images of $\DD^+_i$ are eventually periodic. 
That is, there are $m>0$ and a cycle of discs $\PP$
such that $h^m(D_i^+)\subset\abs{\PP}$. Since $C_0$
is $h$-invariant, this shows the lemma.
\qed

\begin{lemma}\label{lnew2}
If $C_0\cap\abs{\PP}\neq\emptyset$ for
some cycle of discs $\PP$ of $\DD^+$, then $C_0\subset\abs{\PP}$.
\end{lemma}

\bd The set $H^{-1}([-a',a])\setminus \abs{\PP}$ is compact,
as well as $\abs{\PP}$.
Thus there is $\varepsilon_0>0$ such that  any 
$z\in H^{-1}([-a',a])\setminus \abs{\PP}$ and $w\in\abs{\PP}$ satisfy
$d(z,w)>\varepsilon_0$.

Let $x\in C_0\cap\abs{\PP}$ and 
let  $y$ be an arbitrary point in $C_0$. Choose $\varepsilon$ small enough
so that $\varepsilon<\varepsilon_0$ and that any $\varepsilon$-chain
$x_0,x_1,\ldots,x_r$ from $x$ to $y$
is contained in $H^{-1}((-a',a))$ (Lemma \ref{l5}).
We shall show inductively that $x_i\in\abs{\PP}$.
This is true for $i=0$. Assume $x_{i-1}\in \abs{\PP}$.
Then $h(x_{i-1})\in\abs{\PP}$ since $\abs{\PP}$ is forward invariant.
On the other hand,
 $d(x_i,h(x_{i-1}))< \varepsilon_0$ and $x_i\in H^{-1}([-a',a])$.
By the definition of $\varepsilon_0$, this implies $x_i\in\abs{\PP}$.
Inductively we have $y\in\abs{\PP}$, as is required.
\qed

\begin{lemma} \label{l4}
If two periodic points $z_\nu$ ($\nu=1,2$) are contained in $\abs{\PP}$,
where $\PP$ is a cycle of discs in $\DD^+$, then we have
$\rot(\tilde h,z_1)=\rot(\tilde h,z_2)$.
\end{lemma}

\bd 
By replacing $z_\nu$ by their iterate,
one may assume both $z_\nu$ belong to the disc $D_{i_0}^+$.
Choose the lift $\tilde z_\nu$ of $z_\nu$ from the same lift
of $D_{i_0}^+$. Then for any $j\in\N$, their
images $\tilde h^j(\tilde z_\nu)$ must belong to the same lift of
the same disc $D_{i_j}^+$, showing the lemma. \qed

\medskip

\noindent
{\sc Proof of Proposition \ref{l6}}. Assume on the contrary
that $C_0\cap A^+\neq\emptyset$.
Then by Lemmas \ref{lnew} and \ref{lnew2}, $C_0\subset\abs{\PP}$
for a cycle of discs $\PP$ in $\DD^+$. But then Lemma \ref{l4}
contradicts the assumption of $C_0$ (the existence of two periodic
points of different rotation number).  The case 
$C_0\cap A^-\neq\emptyset$ can be dealt with similarly.
\qed

\medskip
Now let us deform the homeomorphism $h$ in $A^-\cup A^+$
so that it has no
periodic points in $A^-\cup A^+$. 
Namely we replace $h$ with a map in $\HH$ with very simple dynamics
in $A^-\cup A^+$. 
Notice that for any small $\varepsilon$,
any $\varepsilon$-chain starting and ending at $C_0$ never falls
into $A^+\cup A^-$.
Proposition \ref{l6}, together with this fact, shows 
that the chain transitive class $C_0$ of the
old $h$ is unchanged for the new $h$.
Especially the points $x_\nu\in C_0$ 
in the assumption of Proposition \ref{p2a} are
still the periodic points of the new $h$. 
Moreover if we find
a periodic point of the new $h$ in $C_0$, it is
a periodic point of the old $h$ in $C_0$ of the same rotation number. 
Therefore in the proof of 
Proposition \ref{p2a}, it is no loss of generality to assume the following.

\begin{assumption} \label{a1}
There is $\beta>0$ such that for any $z\in A^-\cup A^+$,
we have $d(z, h^q(z))>\beta$.
\end{assumption}

\medskip
The rest of the proof is divided into two cases according to
whether $  F_0\cap(\bigcup_k E_k^+\cup\bigcup_lE_l^-)=\emptyset$ or not, where
 $F_0=\pi({\rm Fix}(\tilde h^q\circ T^{-p_1-1}))$.

\medskip
\noindent
{\sc Case 1.} $  F_0\cap(\bigcup_k E_k^+\cup\bigcup_l
E_l^-)\neq\emptyset$.

The argument in this case is based upon the Nielsen fixed
point theory (\cite{J}), which is a refinement of
the Lefschetz index theorem. 
Let us give a brief summary of the theory for the special case
of a continuous map $f$ of the closed annulus $\BBB$.
Let us denote by $\pi:\tilde \BBB\to\BBB$ the universal covering map.
Let $\{\tilde f_i\}_{i\in I}$ be the family of the lifts of $f$
to $\tilde \BBB$. Then  $F_i=\pi({\rm Fix}(\tilde f_i))$ is
a closed subset of ${\rm Fix}(f)$,
called a {\em Nielsen class} of ${\rm Fix}(f)$. It is empty but for finitely many
lifts $\tilde f_i$,
and  ${\rm Fix}(f)$ is partitioned into a finite disjoint union
of nonempty Nielsen classes. To each Nielsen class $F_i$, an integer
$\Index(f,F_i)$, called the {\em index of $F_i$},
is assigned so that the sum of indices is equal to the Lefschetz number
of $f$. 

The most important feature of the index is the following. Suppose
that two maps $f$ and $f'$ are homotopic, and 
a lift $\tilde f_i$ of $f$ is joined with a lift $\tilde f'_i$
of $f'$
by a lift of the homotopy. Then the corresponding Nielsen classes
$F_i=\pi({\rm Fix}(\tilde f'_i))$ and $F_i'=\pi({\rm Fix}(\tilde f'_i))$ 
have the same index: \ $\Index(f,F_i)=\Index(f',F_i')$.

In particular if $f:\BBB\to\BBB$ is homotopic to the identity, then for any Nielsen
class $F_i$, we have $\Index(f,F_i)=0$, since $f$ is homotopic to a
fixed point free homeomorphism.

The index $\Index(f,F_i)$ is computed as follows.
If a Nielsen class $F_i$ is partitioned into a finite disjoint union
of closed subsets: $F_i=\bigcup_jG_j$, then
$$
\Index(f,F_i)=\sum_j\Index(f,G_j).$$
Assume there is a closed disc $D$ such that
\begin{equation} \label{e3.8}
G_j={\rm Fix}(f)\cap D\subset\Int(D).
\end{equation}
Let $\tilde f_i$ be the lift correspoinding to the Nielsen class $F_i$
that contains $G_j$
and $\tilde D$ any lift of $D$. Consider an inclusion
$\tilde\BBB\subset\R^2$. Then $\Index(f,G_j)$ is the mapping degree
of the map
$$
{\rm Id}-\tilde f_i:\partial \tilde D\to\R^2\setminus\{0\}.$$
This is independent of the choice of the disc $D$ satisfying (\ref{e3.8}).
In particular 
if $G_j$ is nonempty and if $f^{-1}(D)\subset \Int(D)$, then $\Index(f,G_j)=1$.

Now let us start the proof of Proposition \ref{p2a} in Case 1.
We apply the Nielsen fixed point theory to the map
$h^q$. 
For this purpose, the homeomophism $h^q:\AAA\to\AAA$ must be deformed in
the exterior of a compact subannulus and extended to a homeomorphism of
$\BBB$ in such a way that the fixed point of the new extended
$h^q$ is the same as
the original $h^q$. But this can easily be done.
In the sequal, we forget about this change, and just consider
the original $h^q$.

We are interested in the particular lift 
$\tilde h^q\circ
T^{-p_1-1}$
and the corresponding Nielsen class
$F_0=\pi({\rm Fix}(\tilde h^q\circ T^{-p_1-1}))$.
Our purpose is to show that 
$F_0\cap H^{-1}([-a',a])\neq\emptyset$. 
We have
\begin{equation}\label{e3.9}
\Index(h^q,F_0)=0.
\end{equation}

Assume $F_0\cap E_k^+\neq\emptyset$ for some $k$. Then we have
$h^{-q}(E_k^+)\subset\Int(E_k^+)$.  
Condition (\ref{e3.8}) above is satisfied
for $f=h^q$, $D=E_k^+$ and $G_j=F_0\cap E_k^+$. Thus we have 
$\Index(h^q,F_0\cap E_k^+)=1$. Likewise if $F_0\cap E_l^-\neq\emptyset$,
then $\Index(h^q,F_0\cap E_l^-)=1$. 

On the other hand by Assumption \ref{a1}, $F_0\cap(A^-\cup A^+)=\emptyset$.
By (\ref{e3.9}), this implies that $\Index(h^q,F_0\cap X)<0$, showing
that $F_0\cap X\neq\emptyset$, and hence
$F_0\cap H^{-1}([-a',a])\neq\emptyset$, as is required.

\medskip
\noindent
{\sc Case 2.} $  F_0\cap(\bigcup_k E_k^+\cup\bigcup_l
E_l^-)=\emptyset$.

In this case, $F_0$, if nonempty, must be contained in $X\subset
H^{-1}([-a',a])$.
Therefore we only need to show
that $F_0$ is nonempty in $\AAA$. The proof is by absurdity.
Assume throughout Case 2 
that the map $\tilde h^q\circ T^{-p_1-1}$ is fixed point free.
This, together with Assumption \ref{a1}, 
implies that there is $\alpha>0$ with the following property.
\begin{enumerate} 
\item For any $\tilde z\in\tilde\AAA$,
$d(\tilde h^q(\tilde z),T^{p_1+1}(\tilde z))>2\alpha$.
\end{enumerate}
Here $d$ denotes the distance function given by the lift
of the standard Riemannian metric $d\theta^2+dt^2$ of $\AAA$. Thus the covering
transformation $T$ is an isometry for $d$.

There is $\delta>0$ such that for a lift $\tilde\varphi$ of
a homeomorphism $\varphi$ of $\AAA$, the following holds. We denote by
$\Vert\cdot\Vert_0$ the supremum norm.

\begin{enumerate} \addtocounter{enumi}{1}
\item
If $\Vert\tilde\varphi-{\rm Id}\Vert_0<2\delta$, then $\Vert (\tilde \varphi\circ\tilde h)^q-\tilde h^q\Vert_0<\alpha$.
\end{enumerate}

Conditions (1) and (2) implies in particular that
for any $\tilde z\in\tilde\AAA$, we have
$$d((\tilde\varphi\circ\tilde h)^q(\tilde z),T^{p_1+1}(\tilde z))>\alpha,$$
and therefore we have the following.
\begin{enumerate} \addtocounter{enumi}{2}
\item There is no fixed point of $(\tilde\varphi\circ\tilde h)^q\circ T^{-p_1-1}$.
\end{enumerate}

Fix once and for all the number $\delta>0$ that satisfies (2).

Recall the periodic points $x_\nu$ and their lift $\tilde x_\nu$
in the assumption of Proposition \ref{p2a}.
Consider a $\delta$-chain $\gamma=(z_0,z_1,\ldots,z_i)$ of length $i$
from $x_\nu$ to $x_{\nu'}$ ($\nu,\,\nu'=1,2$). 
Let 
$\tilde \gamma=(\tilde z_0,\tilde z_1,\ldots,\tilde z_i)$ be 
a lift of $\gamma$
starting at
$\tilde x_\nu$ which is a $\delta$-chain for $\tilde h$.
Assume that $\tilde\gamma$ ends at $T^j(\tilde x_{\nu'})$ for some $j\in\Z$. Then
the pair $(i,j)$ is called the {\em dynamical index} of $\gamma$.
We have the following lemma, which is a variant of the method for
finding periodic points invented in
\cite{F2}. 

\begin{lemma} \label{l1}
There is no $\delta$-cycle at $x_1$ of dynamical index $(\xi q,\xi (p_1+1))$ 
for any $\xi\in\N$.
\end{lemma}

\bd Assume for contradiction that there is a $\delta$-cycle
$\gamma=(z_0,z_1\ldots, z_r)$ at $x_1$ of dynamical index $\xi(q,p_1+1)$ for some
$\xi>0$, Thus $r=\xi q$ and $z_0=z_r=x_1$. 
Then there is a homeomorphism $\varphi$ of $\AAA$ such that
$\varphi(h(z_i))=z_{i+1}$ ($0\leq i<\xi q$) and that 
$\Vert \varphi-{\rm Id}\Vert_0 <2\delta$. 

To show this, consider
the product $\AAA\times[0,1]$ and the line segments 
joining $(h(z_i),0)$ to $(z_{i+1},1)$. A general position argument
shows  that
the line segments can be moved slightly so that
they are mutually disjoint. Define a vector field $X$ pointing upwards, tangent to
the segments. With an appropriate choice of $X$, the holonomy map of $X$
from $\AAA\times\{0\}$ to $\AAA\times\{1\}$
yields a desired homeomorphism $\varphi$. See Figure 3.

\begin{figure}
\unitlength 0.1in
\begin{picture}( 38.5000, 28.3000)( 16.4000,-33.3000)
%
{\color[named]{Black}{%
\special{pn 8}%
\special{pa 4875 1236}%
\special{pa 4868 1238}%
\special{fp}%
\special{pa 4831 1245}%
\special{pa 4829 1246}%
\special{pa 4824 1247}%
\special{fp}%
\special{pa 4787 1254}%
\special{pa 4785 1255}%
\special{pa 4780 1256}%
\special{pa 4780 1256}%
\special{fp}%
\special{pa 4743 1263}%
\special{pa 4739 1264}%
\special{pa 4736 1264}%
\special{fp}%
\special{pa 4699 1270}%
\special{pa 4692 1272}%
\special{fp}%
\special{pa 4655 1278}%
\special{pa 4654 1278}%
\special{pa 4648 1280}%
\special{fp}%
\special{pa 4611 1285}%
\special{pa 4610 1285}%
\special{pa 4606 1286}%
\special{pa 4603 1287}%
\special{fp}%
\special{pa 4566 1292}%
\special{pa 4566 1292}%
\special{pa 4561 1293}%
\special{pa 4559 1293}%
\special{fp}%
\special{pa 4522 1298}%
\special{pa 4520 1298}%
\special{pa 4515 1299}%
\special{pa 4515 1299}%
\special{fp}%
\special{pa 4478 1304}%
\special{pa 4474 1305}%
\special{pa 4470 1305}%
\special{fp}%
\special{pa 4433 1310}%
\special{pa 4431 1310}%
\special{pa 4426 1310}%
\special{pa 4426 1310}%
\special{fp}%
\special{pa 4389 1315}%
\special{pa 4388 1315}%
\special{pa 4383 1315}%
\special{pa 4381 1315}%
\special{fp}%
\special{pa 4344 1320}%
\special{pa 4340 1320}%
\special{pa 4337 1321}%
\special{fp}%
\special{pa 4300 1324}%
\special{pa 4296 1325}%
\special{pa 4292 1325}%
\special{fp}%
\special{pa 4255 1328}%
\special{pa 4251 1329}%
\special{pa 4248 1329}%
\special{fp}%
\special{pa 4210 1332}%
\special{pa 4206 1333}%
\special{pa 4203 1333}%
\special{fp}%
\special{pa 4166 1336}%
\special{pa 4161 1336}%
\special{pa 4158 1337}%
\special{fp}%
\special{pa 4121 1339}%
\special{pa 4120 1339}%
\special{pa 4115 1340}%
\special{pa 4114 1340}%
\special{fp}%
\special{pa 4076 1343}%
\special{pa 4074 1343}%
\special{pa 4069 1343}%
\special{fp}%
\special{pa 4032 1345}%
\special{pa 4027 1346}%
\special{pa 4024 1346}%
\special{fp}%
\special{pa 3987 1348}%
\special{pa 3981 1348}%
\special{pa 3979 1348}%
\special{fp}%
\special{pa 3942 1350}%
\special{pa 3939 1350}%
\special{pa 3934 1351}%
\special{fp}%
\special{pa 3897 1352}%
\special{pa 3897 1352}%
\special{pa 3891 1353}%
\special{pa 3890 1353}%
\special{fp}%
\special{pa 3852 1354}%
\special{pa 3845 1354}%
\special{fp}%
\special{pa 3807 1356}%
\special{pa 3800 1356}%
\special{fp}%
\special{pa 3762 1357}%
\special{pa 3755 1357}%
\special{fp}%
\special{pa 3718 1358}%
\special{pa 3710 1358}%
\special{fp}%
\special{pa 3673 1359}%
\special{pa 3665 1359}%
\special{fp}%
\special{pa 3628 1360}%
\special{pa 3620 1360}%
\special{fp}%
\special{pa 3583 1360}%
\special{pa 3575 1360}%
\special{fp}%
\special{pa 3538 1360}%
\special{pa 3530 1360}%
\special{fp}%
\special{pa 3493 1360}%
\special{pa 3485 1360}%
\special{fp}%
\special{pa 3448 1359}%
\special{pa 3440 1359}%
\special{fp}%
\special{pa 3403 1359}%
\special{pa 3395 1359}%
\special{fp}%
\special{pa 3358 1358}%
\special{pa 3350 1358}%
\special{fp}%
\special{pa 3313 1357}%
\special{pa 3305 1357}%
\special{fp}%
\special{pa 3268 1356}%
\special{pa 3266 1356}%
\special{pa 3261 1355}%
\special{pa 3260 1355}%
\special{fp}%
\special{pa 3223 1354}%
\special{pa 3215 1354}%
\special{fp}%
\special{pa 3178 1352}%
\special{pa 3170 1352}%
\special{fp}%
\special{pa 3133 1350}%
\special{pa 3126 1350}%
\special{fp}%
\special{pa 3088 1348}%
\special{pa 3087 1348}%
\special{pa 3082 1347}%
\special{pa 3081 1347}%
\special{fp}%
\special{pa 3044 1345}%
\special{pa 3036 1345}%
\special{fp}%
\special{pa 2999 1342}%
\special{pa 2991 1342}%
\special{fp}%
\special{pa 2954 1339}%
\special{pa 2948 1339}%
\special{pa 2947 1339}%
\special{fp}%
\special{pa 2909 1335}%
\special{pa 2907 1335}%
\special{pa 2902 1335}%
\special{fp}%
\special{pa 2865 1332}%
\special{pa 2862 1332}%
\special{pa 2857 1331}%
\special{fp}%
\special{pa 2820 1328}%
\special{pa 2817 1328}%
\special{pa 2813 1327}%
\special{fp}%
\special{pa 2776 1324}%
\special{pa 2773 1323}%
\special{pa 2768 1323}%
\special{pa 2768 1323}%
\special{fp}%
\special{pa 2731 1319}%
\special{pa 2729 1319}%
\special{pa 2724 1318}%
\special{pa 2723 1318}%
\special{fp}%
\special{pa 2686 1314}%
\special{pa 2685 1314}%
\special{pa 2681 1313}%
\special{pa 2679 1313}%
\special{fp}%
\special{pa 2642 1309}%
\special{pa 2642 1309}%
\special{pa 2638 1308}%
\special{pa 2634 1308}%
\special{fp}%
\special{pa 2598 1303}%
\special{pa 2596 1303}%
\special{pa 2591 1303}%
\special{pa 2590 1303}%
\special{fp}%
\special{pa 2553 1298}%
\special{pa 2549 1297}%
\special{pa 2546 1296}%
\special{fp}%
\special{pa 2509 1291}%
\special{pa 2508 1291}%
\special{pa 2504 1290}%
\special{pa 2501 1290}%
\special{fp}%
\special{pa 2465 1284}%
\special{pa 2464 1284}%
\special{pa 2459 1284}%
\special{pa 2457 1284}%
\special{fp}%
\special{pa 2420 1277}%
\special{pa 2420 1277}%
\special{pa 2416 1277}%
\special{pa 2413 1276}%
\special{fp}%
\special{pa 2376 1270}%
\special{pa 2373 1269}%
\special{pa 2369 1269}%
\special{pa 2369 1269}%
\special{fp}%
\special{pa 2332 1262}%
\special{pa 2331 1262}%
\special{pa 2325 1260}%
\special{fp}%
\special{pa 2288 1254}%
\special{pa 2281 1252}%
\special{fp}%
\special{pa 2244 1245}%
\special{pa 2237 1243}%
\special{fp}%
\special{pa 2200 1235}%
\special{pa 2193 1233}%
\special{fp}%
\special{pa 2157 1224}%
\special{pa 2151 1223}%
\special{pa 2149 1223}%
\special{fp}%
\special{pa 2113 1214}%
\special{pa 2112 1214}%
\special{pa 2108 1213}%
\special{pa 2106 1212}%
\special{fp}%
\special{pa 2070 1202}%
\special{pa 2067 1201}%
\special{pa 2063 1201}%
\special{pa 2062 1201}%
\special{fp}%
\special{pa 2027 1190}%
\special{pa 2024 1189}%
\special{pa 2020 1188}%
\special{pa 2019 1188}%
\special{fp}%
\special{pa 1984 1176}%
\special{pa 1983 1176}%
\special{pa 1979 1175}%
\special{pa 1976 1174}%
\special{fp}%
\special{pa 1941 1162}%
\special{pa 1934 1160}%
\special{fp}%
\special{pa 1899 1147}%
\special{pa 1893 1145}%
\special{pa 1892 1144}%
\special{fp}%
\special{pa 1857 1130}%
\special{pa 1857 1130}%
\special{pa 1855 1129}%
\special{pa 1852 1127}%
\special{pa 1851 1126}%
\special{fp}%
\special{pa 1817 1111}%
\special{pa 1814 1110}%
\special{pa 1810 1108}%
\special{fp}%
\special{pa 1777 1090}%
\special{pa 1770 1087}%
\special{fp}%
\special{pa 1739 1067}%
\special{pa 1737 1066}%
\special{pa 1736 1065}%
\special{pa 1733 1063}%
\special{fp}%
\special{pa 1704 1040}%
\special{pa 1704 1040}%
\special{pa 1702 1039}%
\special{pa 1700 1037}%
\special{pa 1698 1036}%
\special{fp}%
\special{pa 1673 1009}%
\special{pa 1672 1008}%
\special{pa 1668 1004}%
\special{pa 1668 1003}%
\special{fp}%
\special{pa 1649 973}%
\special{pa 1649 972}%
\special{pa 1649 971}%
\special{pa 1648 970}%
\special{pa 1648 969}%
\special{pa 1647 967}%
\special{pa 1647 966}%
\special{pa 1647 966}%
\special{fp}%
\special{pa 1640 931}%
\special{pa 1640 924}%
\special{fp}%
\special{pa 1648 890}%
\special{pa 1649 889}%
\special{pa 1649 887}%
\special{pa 1650 886}%
\special{pa 1650 885}%
\special{pa 1652 883}%
\special{fp}%
\special{pa 1671 853}%
\special{pa 1671 853}%
\special{pa 1672 851}%
\special{pa 1676 847}%
\special{fp}%
\special{pa 1702 822}%
\special{pa 1703 821}%
\special{pa 1704 819}%
\special{pa 1705 818}%
\special{pa 1707 817}%
\special{pa 1707 817}%
\special{fp}%
\special{pa 1737 795}%
\special{pa 1740 793}%
\special{pa 1741 791}%
\special{pa 1742 790}%
\special{fp}%
\special{pa 1774 771}%
\special{pa 1781 768}%
\special{pa 1781 768}%
\special{fp}%
\special{pa 1814 751}%
\special{pa 1819 748}%
\special{pa 1821 747}%
\special{fp}%
\special{pa 1854 731}%
\special{pa 1855 731}%
\special{pa 1858 730}%
\special{pa 1860 729}%
\special{pa 1861 729}%
\special{fp}%
\special{pa 1896 714}%
\special{pa 1897 714}%
\special{pa 1899 713}%
\special{pa 1903 712}%
\special{fp}%
\special{pa 1938 699}%
\special{pa 1938 699}%
\special{pa 1945 697}%
\special{fp}%
\special{pa 1980 684}%
\special{pa 1983 684}%
\special{pa 1988 682}%
\special{fp}%
\special{pa 2023 671}%
\special{pa 2027 670}%
\special{pa 2031 669}%
\special{fp}%
\special{pa 2066 658}%
\special{pa 2067 658}%
\special{pa 2071 657}%
\special{pa 2074 656}%
\special{fp}%
\special{pa 2110 647}%
\special{pa 2112 646}%
\special{pa 2116 645}%
\special{pa 2117 645}%
\special{fp}%
\special{pa 2153 636}%
\special{pa 2156 635}%
\special{pa 2159 635}%
\special{pa 2161 635}%
\special{fp}%
\special{pa 2197 626}%
\special{pa 2201 625}%
\special{pa 2204 624}%
\special{pa 2204 624}%
\special{fp}%
\special{pa 2241 617}%
\special{pa 2248 615}%
\special{fp}%
\special{pa 2285 608}%
\special{pa 2292 606}%
\special{fp}%
\special{pa 2329 599}%
\special{pa 2336 597}%
\special{pa 2336 597}%
\special{fp}%
\special{pa 2373 591}%
\special{pa 2374 591}%
\special{pa 2380 589}%
\special{fp}%
\special{pa 2417 583}%
\special{pa 2417 583}%
\special{pa 2421 583}%
\special{pa 2424 582}%
\special{fp}%
\special{pa 2461 576}%
\special{pa 2465 575}%
\special{pa 2469 575}%
\special{fp}%
\special{pa 2505 569}%
\special{pa 2509 569}%
\special{pa 2513 568}%
\special{fp}%
\special{pa 2550 563}%
\special{pa 2550 563}%
\special{pa 2555 562}%
\special{pa 2557 562}%
\special{fp}%
\special{pa 2594 557}%
\special{pa 2597 557}%
\special{pa 2601 556}%
\special{pa 2602 556}%
\special{fp}%
\special{pa 2639 551}%
\special{pa 2639 551}%
\special{pa 2643 551}%
\special{pa 2646 550}%
\special{fp}%
\special{pa 2683 546}%
\special{pa 2686 546}%
\special{pa 2691 545}%
\special{fp}%
\special{pa 2728 541}%
\special{pa 2730 541}%
\special{pa 2735 541}%
\special{pa 2735 541}%
\special{fp}%
\special{pa 2772 537}%
\special{pa 2774 537}%
\special{pa 2779 536}%
\special{pa 2780 536}%
\special{fp}%
\special{pa 2817 532}%
\special{pa 2818 532}%
\special{pa 2823 532}%
\special{pa 2824 532}%
\special{fp}%
\special{pa 2861 528}%
\special{pa 2863 528}%
\special{pa 2868 528}%
\special{pa 2869 528}%
\special{fp}%
\special{pa 2906 524}%
\special{pa 2908 524}%
\special{pa 2914 524}%
\special{fp}%
\special{pa 2951 521}%
\special{pa 2958 521}%
\special{fp}%
\special{pa 2996 518}%
\special{pa 3000 518}%
\special{pa 3003 517}%
\special{fp}%
\special{pa 3040 515}%
\special{pa 3047 515}%
\special{pa 3048 515}%
\special{fp}%
\special{pa 3085 513}%
\special{pa 3088 512}%
\special{pa 3093 512}%
\special{fp}%
\special{pa 3130 510}%
\special{pa 3137 510}%
\special{fp}%
\special{pa 3175 508}%
\special{pa 3182 508}%
\special{fp}%
\special{pa 3220 506}%
\special{pa 3227 506}%
\special{fp}%
\special{pa 3264 505}%
\special{pa 3267 504}%
\special{pa 3272 504}%
\special{fp}%
\special{pa 3309 503}%
\special{pa 3317 503}%
\special{fp}%
\special{pa 3354 502}%
\special{pa 3362 502}%
\special{fp}%
\special{pa 3399 501}%
\special{pa 3407 501}%
\special{fp}%
\special{pa 3444 501}%
\special{pa 3449 500}%
\special{pa 3452 500}%
\special{fp}%
\special{pa 3489 500}%
\special{pa 3497 500}%
\special{fp}%
\special{pa 3534 500}%
\special{pa 3542 500}%
\special{fp}%
\special{pa 3579 500}%
\special{pa 3587 500}%
\special{fp}%
\special{pa 3624 500}%
\special{pa 3631 500}%
\special{pa 3632 500}%
\special{fp}%
\special{pa 3669 501}%
\special{pa 3677 501}%
\special{fp}%
\special{pa 3714 502}%
\special{pa 3722 502}%
\special{fp}%
\special{pa 3759 503}%
\special{pa 3767 503}%
\special{fp}%
\special{pa 3804 504}%
\special{pa 3812 504}%
\special{fp}%
\special{pa 3849 506}%
\special{pa 3856 506}%
\special{fp}%
\special{pa 3894 507}%
\special{pa 3898 508}%
\special{pa 3901 508}%
\special{fp}%
\special{pa 3939 510}%
\special{pa 3940 510}%
\special{pa 3946 510}%
\special{fp}%
\special{pa 3983 512}%
\special{pa 3991 512}%
\special{fp}%
\special{pa 4028 514}%
\special{pa 4034 515}%
\special{pa 4036 515}%
\special{fp}%
\special{pa 4073 517}%
\special{pa 4075 517}%
\special{pa 4080 518}%
\special{pa 4081 518}%
\special{fp}%
\special{pa 4118 520}%
\special{pa 4121 521}%
\special{pa 4125 521}%
\special{fp}%
\special{pa 4162 524}%
\special{pa 4170 524}%
\special{fp}%
\special{pa 4207 527}%
\special{pa 4212 528}%
\special{pa 4215 528}%
\special{fp}%
\special{pa 4252 531}%
\special{pa 4252 531}%
\special{pa 4257 532}%
\special{pa 4259 532}%
\special{fp}%
\special{pa 4296 536}%
\special{pa 4297 536}%
\special{pa 4302 536}%
\special{pa 4304 536}%
\special{fp}%
\special{pa 4341 540}%
\special{pa 4341 540}%
\special{pa 4346 541}%
\special{pa 4348 541}%
\special{fp}%
\special{pa 4385 545}%
\special{pa 4389 545}%
\special{pa 4393 546}%
\special{fp}%
\special{pa 4430 550}%
\special{pa 4432 550}%
\special{pa 4437 551}%
\special{fp}%
\special{pa 4474 556}%
\special{pa 4475 556}%
\special{pa 4479 556}%
\special{pa 4482 557}%
\special{fp}%
\special{pa 4519 562}%
\special{pa 4521 562}%
\special{pa 4526 562}%
\special{pa 4526 562}%
\special{fp}%
\special{pa 4563 568}%
\special{pa 4566 568}%
\special{pa 4571 569}%
\special{fp}%
\special{pa 4607 574}%
\special{pa 4611 575}%
\special{pa 4615 576}%
\special{fp}%
\special{pa 4652 581}%
\special{pa 4659 583}%
\special{pa 4659 583}%
\special{fp}%
\special{pa 4696 588}%
\special{pa 4702 590}%
\special{pa 4703 590}%
\special{fp}%
\special{pa 4740 597}%
\special{pa 4740 597}%
\special{pa 4744 597}%
\special{pa 4747 598}%
\special{fp}%
\special{pa 4784 605}%
\special{pa 4791 607}%
\special{fp}%
\special{pa 4828 614}%
\special{pa 4835 616}%
\special{fp}%
\special{pa 4872 623}%
\special{pa 4872 623}%
\special{pa 4879 625}%
\special{fp}%
\special{pa 4915 633}%
\special{pa 4917 634}%
\special{pa 4923 635}%
\special{fp}%
\special{pa 4959 644}%
\special{pa 4961 645}%
\special{pa 4964 645}%
\special{pa 4966 646}%
\special{fp}%
\special{pa 5002 656}%
\special{pa 5003 656}%
\special{pa 5006 657}%
\special{pa 5010 657}%
\special{fp}%
\special{pa 5046 668}%
\special{pa 5046 668}%
\special{pa 5050 669}%
\special{pa 5053 670}%
\special{fp}%
\special{pa 5089 681}%
\special{pa 5096 684}%
\special{fp}%
\special{pa 5131 695}%
\special{pa 5139 698}%
\special{fp}%
\special{pa 5174 710}%
\special{pa 5176 711}%
\special{pa 5178 712}%
\special{pa 5181 713}%
\special{fp}%
\special{pa 5215 727}%
\special{pa 5218 728}%
\special{pa 5220 729}%
\special{pa 5222 730}%
\special{fp}%
\special{pa 5256 745}%
\special{pa 5258 747}%
\special{pa 5261 748}%
\special{pa 5263 749}%
\special{fp}%
\special{pa 5296 766}%
\special{pa 5297 767}%
\special{pa 5302 770}%
\special{fp}%
\special{pa 5334 788}%
\special{pa 5335 789}%
\special{pa 5337 790}%
\special{pa 5339 792}%
\special{pa 5340 793}%
\special{fp}%
\special{pa 5369 814}%
\special{pa 5371 815}%
\special{pa 5373 817}%
\special{pa 5375 818}%
\special{pa 5376 819}%
\special{fp}%
\special{pa 5402 844}%
\special{pa 5403 846}%
\special{pa 5407 850}%
\special{fp}%
\special{pa 5427 880}%
\special{pa 5427 880}%
\special{pa 5428 881}%
\special{pa 5428 883}%
\special{pa 5430 885}%
\special{pa 5430 886}%
\special{pa 5430 886}%
\special{fp}%
\special{pa 5440 920}%
\special{pa 5440 920}%
\special{pa 5440 927}%
\special{fp}%
\special{pa 5434 963}%
\special{pa 5434 963}%
\special{pa 5434 965}%
\special{pa 5433 966}%
\special{pa 5433 968}%
\special{pa 5432 969}%
\special{pa 5432 969}%
\special{fp}%
\special{pa 5415 1000}%
\special{pa 5414 1001}%
\special{pa 5413 1003}%
\special{pa 5410 1006}%
\special{fp}%
\special{pa 5384 1033}%
\special{pa 5384 1033}%
\special{pa 5383 1035}%
\special{pa 5382 1036}%
\special{pa 5380 1037}%
\special{pa 5379 1038}%
\special{fp}%
\special{pa 5350 1060}%
\special{pa 5349 1061}%
\special{pa 5346 1064}%
\special{pa 5344 1065}%
\special{fp}%
\special{pa 5313 1084}%
\special{pa 5311 1085}%
\special{pa 5308 1088}%
\special{pa 5307 1089}%
\special{fp}%
\special{pa 5274 1106}%
\special{pa 5268 1109}%
\special{pa 5267 1109}%
\special{fp}%
\special{pa 5233 1125}%
\special{pa 5233 1125}%
\special{pa 5230 1126}%
\special{pa 5228 1128}%
\special{pa 5227 1128}%
\special{fp}%
\special{pa 5193 1142}%
\special{pa 5192 1143}%
\special{pa 5186 1145}%
\special{pa 5185 1145}%
\special{fp}%
\special{pa 5150 1158}%
\special{pa 5143 1161}%
\special{fp}%
\special{pa 5108 1172}%
\special{pa 5101 1175}%
\special{fp}%
\special{pa 5065 1186}%
\special{pa 5063 1187}%
\special{pa 5059 1188}%
\special{pa 5058 1188}%
\special{fp}%
\special{pa 5022 1199}%
\special{pa 5017 1201}%
\special{pa 5015 1202}%
\special{fp}%
\special{pa 4979 1211}%
\special{pa 4975 1212}%
\special{pa 4972 1213}%
\special{pa 4971 1213}%
\special{fp}%
\special{pa 4935 1222}%
\special{pa 4928 1224}%
\special{pa 4928 1224}%
\special{fp}%
\special{pa 4891 1232}%
\special{pa 4887 1233}%
\special{pa 4884 1234}%
\special{fp}%
}}%
%
{\color[named]{Black}{%
\special{pn 8}%
\special{pa 5490 2870}%
\special{pa 5490 2878}%
\special{fp}%
\special{pa 5482 2912}%
\special{pa 5482 2913}%
\special{pa 5481 2914}%
\special{pa 5481 2915}%
\special{pa 5480 2917}%
\special{pa 5480 2918}%
\special{pa 5479 2919}%
\special{fp}%
\special{pa 5461 2950}%
\special{pa 5461 2950}%
\special{pa 5461 2952}%
\special{pa 5458 2955}%
\special{pa 5457 2956}%
\special{fp}%
\special{pa 5433 2983}%
\special{pa 5432 2984}%
\special{pa 5427 2989}%
\special{fp}%
\special{pa 5400 3012}%
\special{pa 5398 3013}%
\special{pa 5397 3014}%
\special{pa 5395 3015}%
\special{pa 5394 3017}%
\special{fp}%
\special{pa 5364 3037}%
\special{pa 5362 3038}%
\special{pa 5361 3039}%
\special{pa 5357 3041}%
\special{pa 5357 3041}%
\special{fp}%
\special{pa 5325 3060}%
\special{pa 5321 3062}%
\special{pa 5318 3063}%
\special{fp}%
\special{pa 5285 3080}%
\special{pa 5285 3080}%
\special{pa 5283 3082}%
\special{pa 5280 3083}%
\special{pa 5279 3084}%
\special{fp}%
\special{pa 5245 3099}%
\special{pa 5244 3099}%
\special{pa 5241 3100}%
\special{pa 5239 3102}%
\special{pa 5238 3102}%
\special{fp}%
\special{pa 5203 3116}%
\special{pa 5202 3117}%
\special{pa 5196 3119}%
\special{fp}%
\special{pa 5161 3131}%
\special{pa 5156 3133}%
\special{pa 5154 3135}%
\special{fp}%
\special{pa 5119 3147}%
\special{pa 5112 3149}%
\special{pa 5111 3149}%
\special{fp}%
\special{pa 5076 3160}%
\special{pa 5071 3162}%
\special{pa 5068 3163}%
\special{fp}%
\special{pa 5032 3173}%
\special{pa 5028 3174}%
\special{pa 5025 3175}%
\special{fp}%
\special{pa 4989 3184}%
\special{pa 4987 3185}%
\special{pa 4983 3186}%
\special{pa 4982 3186}%
\special{fp}%
\special{pa 4946 3196}%
\special{pa 4944 3196}%
\special{pa 4941 3197}%
\special{pa 4938 3198}%
\special{fp}%
\special{pa 4902 3207}%
\special{pa 4896 3208}%
\special{pa 4894 3208}%
\special{fp}%
\special{pa 4858 3216}%
\special{pa 4851 3218}%
\special{fp}%
\special{pa 4814 3225}%
\special{pa 4808 3227}%
\special{pa 4807 3227}%
\special{fp}%
\special{pa 4770 3235}%
\special{pa 4768 3235}%
\special{pa 4764 3235}%
\special{pa 4762 3235}%
\special{fp}%
\special{pa 4726 3243}%
\special{pa 4724 3243}%
\special{pa 4719 3244}%
\special{pa 4718 3244}%
\special{fp}%
\special{pa 4682 3250}%
\special{pa 4674 3252}%
\special{fp}%
\special{pa 4637 3258}%
\special{pa 4635 3258}%
\special{pa 4631 3258}%
\special{pa 4630 3258}%
\special{fp}%
\special{pa 4593 3264}%
\special{pa 4592 3264}%
\special{pa 4588 3265}%
\special{pa 4585 3266}%
\special{fp}%
\special{pa 4549 3271}%
\special{pa 4548 3271}%
\special{pa 4544 3271}%
\special{pa 4541 3272}%
\special{fp}%
\special{pa 4504 3277}%
\special{pa 4504 3277}%
\special{pa 4499 3277}%
\special{pa 4496 3278}%
\special{fp}%
\special{pa 4460 3283}%
\special{pa 4458 3283}%
\special{pa 4454 3283}%
\special{pa 4452 3283}%
\special{fp}%
\special{pa 4415 3287}%
\special{pa 4412 3288}%
\special{pa 4408 3289}%
\special{pa 4407 3289}%
\special{fp}%
\special{pa 4370 3293}%
\special{pa 4366 3293}%
\special{pa 4363 3294}%
\special{fp}%
\special{pa 4326 3297}%
\special{pa 4324 3297}%
\special{pa 4319 3298}%
\special{pa 4318 3298}%
\special{fp}%
\special{pa 4281 3302}%
\special{pa 4281 3302}%
\special{pa 4273 3302}%
\special{fp}%
\special{pa 4236 3305}%
\special{pa 4233 3306}%
\special{pa 4229 3306}%
\special{fp}%
\special{pa 4191 3309}%
\special{pa 4189 3309}%
\special{pa 4185 3310}%
\special{pa 4184 3310}%
\special{fp}%
\special{pa 4147 3312}%
\special{pa 4146 3312}%
\special{pa 4141 3313}%
\special{pa 4139 3313}%
\special{fp}%
\special{pa 4102 3315}%
\special{pa 4101 3315}%
\special{pa 4097 3316}%
\special{pa 4094 3316}%
\special{fp}%
\special{pa 4057 3318}%
\special{pa 4052 3318}%
\special{pa 4049 3319}%
\special{fp}%
\special{pa 4012 3321}%
\special{pa 4012 3321}%
\special{pa 4005 3321}%
\special{fp}%
\special{pa 3967 3323}%
\special{pa 3960 3323}%
\special{fp}%
\special{pa 3922 3325}%
\special{pa 3922 3325}%
\special{pa 3915 3325}%
\special{fp}%
\special{pa 3877 3326}%
\special{pa 3872 3326}%
\special{pa 3870 3326}%
\special{fp}%
\special{pa 3832 3328}%
\special{pa 3832 3328}%
\special{pa 3825 3328}%
\special{fp}%
\special{pa 3787 3329}%
\special{pa 3786 3329}%
\special{pa 3780 3329}%
\special{fp}%
\special{pa 3742 3329}%
\special{pa 3734 3329}%
\special{fp}%
\special{pa 3697 3330}%
\special{pa 3689 3330}%
\special{fp}%
\special{pa 3652 3330}%
\special{pa 3644 3330}%
\special{fp}%
\special{pa 3607 3330}%
\special{pa 3599 3330}%
\special{fp}%
\special{pa 3562 3330}%
\special{pa 3558 3330}%
\special{pa 3554 3329}%
\special{fp}%
\special{pa 3517 3329}%
\special{pa 3509 3329}%
\special{fp}%
\special{pa 3472 3328}%
\special{pa 3464 3328}%
\special{fp}%
\special{pa 3427 3327}%
\special{pa 3419 3327}%
\special{fp}%
\special{pa 3382 3325}%
\special{pa 3380 3325}%
\special{pa 3374 3325}%
\special{fp}%
\special{pa 3337 3324}%
\special{pa 3330 3324}%
\special{pa 3329 3324}%
\special{fp}%
\special{pa 3292 3322}%
\special{pa 3290 3322}%
\special{pa 3285 3321}%
\special{pa 3284 3321}%
\special{fp}%
\special{pa 3247 3319}%
\special{pa 3245 3319}%
\special{pa 3239 3319}%
\special{fp}%
\special{pa 3202 3317}%
\special{pa 3196 3317}%
\special{pa 3194 3317}%
\special{fp}%
\special{pa 3157 3314}%
\special{pa 3152 3314}%
\special{pa 3149 3313}%
\special{fp}%
\special{pa 3112 3311}%
\special{pa 3108 3311}%
\special{pa 3105 3310}%
\special{fp}%
\special{pa 3067 3308}%
\special{pa 3064 3307}%
\special{pa 3060 3307}%
\special{fp}%
\special{pa 3023 3303}%
\special{pa 3021 3303}%
\special{pa 3015 3303}%
\special{fp}%
\special{pa 2978 3300}%
\special{pa 2973 3299}%
\special{pa 2970 3299}%
\special{fp}%
\special{pa 2933 3295}%
\special{pa 2931 3295}%
\special{pa 2926 3294}%
\special{pa 2926 3294}%
\special{fp}%
\special{pa 2889 3290}%
\special{pa 2884 3290}%
\special{pa 2881 3289}%
\special{fp}%
\special{pa 2844 3285}%
\special{pa 2843 3285}%
\special{pa 2838 3285}%
\special{pa 2836 3285}%
\special{fp}%
\special{pa 2800 3280}%
\special{pa 2797 3279}%
\special{pa 2792 3279}%
\special{pa 2792 3279}%
\special{fp}%
\special{pa 2755 3274}%
\special{pa 2752 3274}%
\special{pa 2748 3273}%
\special{pa 2747 3273}%
\special{fp}%
\special{pa 2711 3268}%
\special{pa 2708 3267}%
\special{pa 2703 3267}%
\special{pa 2703 3267}%
\special{fp}%
\special{pa 2666 3261}%
\special{pa 2664 3261}%
\special{pa 2659 3260}%
\special{fp}%
\special{pa 2622 3254}%
\special{pa 2621 3254}%
\special{pa 2617 3253}%
\special{pa 2614 3253}%
\special{fp}%
\special{pa 2578 3247}%
\special{pa 2571 3245}%
\special{pa 2570 3245}%
\special{fp}%
\special{pa 2533 3239}%
\special{pa 2526 3237}%
\special{fp}%
\special{pa 2489 3230}%
\special{pa 2486 3229}%
\special{pa 2482 3229}%
\special{pa 2482 3229}%
\special{fp}%
\special{pa 2445 3221}%
\special{pa 2443 3221}%
\special{pa 2438 3220}%
\special{fp}%
\special{pa 2401 3212}%
\special{pa 2401 3212}%
\special{pa 2394 3210}%
\special{fp}%
\special{pa 2358 3201}%
\special{pa 2352 3200}%
\special{pa 2350 3199}%
\special{fp}%
\special{pa 2314 3190}%
\special{pa 2313 3190}%
\special{pa 2309 3190}%
\special{pa 2307 3189}%
\special{fp}%
\special{pa 2271 3179}%
\special{pa 2268 3178}%
\special{pa 2264 3178}%
\special{pa 2263 3178}%
\special{fp}%
\special{pa 2227 3167}%
\special{pa 2227 3167}%
\special{pa 2220 3165}%
\special{fp}%
\special{pa 2184 3153}%
\special{pa 2183 3153}%
\special{pa 2179 3152}%
\special{pa 2177 3151}%
\special{fp}%
\special{pa 2141 3139}%
\special{pa 2134 3137}%
\special{fp}%
\special{pa 2099 3125}%
\special{pa 2097 3124}%
\special{pa 2094 3123}%
\special{pa 2092 3122}%
\special{fp}%
\special{pa 2057 3108}%
\special{pa 2056 3108}%
\special{pa 2053 3106}%
\special{pa 2051 3105}%
\special{pa 2050 3105}%
\special{fp}%
\special{pa 2016 3090}%
\special{pa 2015 3090}%
\special{pa 2011 3088}%
\special{pa 2009 3087}%
\special{fp}%
\special{pa 1975 3071}%
\special{pa 1974 3070}%
\special{pa 1971 3069}%
\special{pa 1969 3068}%
\special{pa 1968 3067}%
\special{fp}%
\special{pa 1936 3049}%
\special{pa 1930 3046}%
\special{pa 1930 3046}%
\special{fp}%
\special{pa 1899 3026}%
\special{pa 1896 3023}%
\special{pa 1892 3021}%
\special{fp}%
\special{pa 1864 2999}%
\special{pa 1859 2994}%
\special{pa 1858 2994}%
\special{fp}%
\special{pa 1832 2967}%
\special{pa 1829 2964}%
\special{pa 1828 2962}%
\special{pa 1827 2961}%
\special{fp}%
\special{pa 1806 2931}%
\special{pa 1806 2931}%
\special{pa 1806 2930}%
\special{pa 1805 2929}%
\special{pa 1804 2927}%
\special{pa 1804 2926}%
\special{pa 1803 2925}%
\special{fp}%
\special{pa 1792 2892}%
\special{pa 1792 2890}%
\special{pa 1791 2888}%
\special{pa 1791 2884}%
\special{fp}%
\special{pa 1792 2848}%
\special{pa 1792 2847}%
\special{pa 1793 2845}%
\special{pa 1793 2843}%
\special{pa 1794 2842}%
\special{pa 1794 2841}%
\special{fp}%
\special{pa 1807 2808}%
\special{pa 1808 2806}%
\special{pa 1808 2805}%
\special{pa 1810 2803}%
\special{pa 1811 2802}%
\special{fp}%
\special{pa 1833 2772}%
\special{pa 1833 2772}%
\special{pa 1834 2770}%
\special{pa 1837 2767}%
\special{pa 1837 2766}%
\special{fp}%
\special{pa 1864 2741}%
\special{pa 1864 2741}%
\special{pa 1866 2740}%
\special{pa 1870 2736}%
\special{fp}%
\special{pa 1899 2715}%
\special{pa 1902 2712}%
\special{pa 1905 2710}%
\special{fp}%
\special{pa 1937 2691}%
\special{pa 1942 2688}%
\special{pa 1943 2687}%
\special{fp}%
\special{pa 1976 2669}%
\special{pa 1976 2669}%
\special{pa 1982 2666}%
\special{pa 1983 2666}%
\special{fp}%
\special{pa 2016 2650}%
\special{pa 2018 2649}%
\special{pa 2020 2648}%
\special{pa 2023 2647}%
\special{pa 2023 2647}%
\special{fp}%
\special{pa 2057 2632}%
\special{pa 2058 2631}%
\special{pa 2064 2629}%
\special{pa 2064 2629}%
\special{fp}%
\special{pa 2099 2615}%
\special{pa 2106 2613}%
\special{fp}%
\special{pa 2141 2601}%
\special{pa 2149 2598}%
\special{fp}%
\special{pa 2184 2586}%
\special{pa 2185 2586}%
\special{pa 2189 2585}%
\special{pa 2192 2584}%
\special{fp}%
\special{pa 2227 2573}%
\special{pa 2230 2572}%
\special{pa 2234 2571}%
\special{pa 2235 2571}%
\special{fp}%
\special{pa 2271 2561}%
\special{pa 2271 2561}%
\special{pa 2274 2560}%
\special{pa 2278 2559}%
\special{pa 2278 2559}%
\special{fp}%
\special{pa 2314 2550}%
\special{pa 2316 2549}%
\special{pa 2320 2548}%
\special{pa 2322 2547}%
\special{fp}%
\special{pa 2358 2538}%
\special{pa 2359 2538}%
\special{pa 2365 2536}%
\special{fp}%
\special{pa 2402 2529}%
\special{pa 2409 2527}%
\special{fp}%
\special{pa 2445 2519}%
\special{pa 2446 2519}%
\special{pa 2453 2517}%
\special{fp}%
\special{pa 2489 2510}%
\special{pa 2497 2508}%
\special{fp}%
\special{pa 2534 2501}%
\special{pa 2534 2501}%
\special{pa 2538 2501}%
\special{pa 2541 2500}%
\special{fp}%
\special{pa 2578 2493}%
\special{pa 2583 2492}%
\special{pa 2585 2492}%
\special{fp}%
\special{pa 2622 2486}%
\special{pa 2625 2485}%
\special{pa 2630 2485}%
\special{fp}%
\special{pa 2666 2479}%
\special{pa 2668 2479}%
\special{pa 2673 2478}%
\special{pa 2674 2478}%
\special{fp}%
\special{pa 2711 2472}%
\special{pa 2712 2472}%
\special{pa 2716 2471}%
\special{pa 2718 2471}%
\special{fp}%
\special{pa 2755 2466}%
\special{pa 2756 2466}%
\special{pa 2761 2465}%
\special{pa 2763 2465}%
\special{fp}%
\special{pa 2800 2460}%
\special{pa 2801 2460}%
\special{pa 2806 2459}%
\special{pa 2807 2459}%
\special{fp}%
\special{pa 2844 2455}%
\special{pa 2847 2454}%
\special{pa 2851 2454}%
\special{pa 2852 2454}%
\special{fp}%
\special{pa 2889 2450}%
\special{pa 2893 2449}%
\special{pa 2896 2449}%
\special{fp}%
\special{pa 2933 2445}%
\special{pa 2935 2445}%
\special{pa 2940 2444}%
\special{pa 2941 2444}%
\special{fp}%
\special{pa 2978 2441}%
\special{pa 2982 2440}%
\special{pa 2986 2440}%
\special{fp}%
\special{pa 3023 2436}%
\special{pa 3025 2436}%
\special{pa 3030 2436}%
\special{pa 3030 2436}%
\special{fp}%
\special{pa 3068 2432}%
\special{pa 3069 2432}%
\special{pa 3075 2432}%
\special{fp}%
\special{pa 3112 2429}%
\special{pa 3117 2429}%
\special{pa 3120 2428}%
\special{fp}%
\special{pa 3157 2426}%
\special{pa 3161 2426}%
\special{pa 3165 2425}%
\special{fp}%
\special{pa 3202 2423}%
\special{pa 3210 2423}%
\special{fp}%
\special{pa 3247 2421}%
\special{pa 3250 2420}%
\special{pa 3255 2420}%
\special{fp}%
\special{pa 3292 2418}%
\special{pa 3300 2418}%
\special{fp}%
\special{pa 3337 2416}%
\special{pa 3344 2416}%
\special{fp}%
\special{pa 3382 2415}%
\special{pa 3385 2414}%
\special{pa 3389 2414}%
\special{fp}%
\special{pa 3427 2413}%
\special{pa 3434 2413}%
\special{fp}%
\special{pa 3472 2412}%
\special{pa 3479 2412}%
\special{fp}%
\special{pa 3517 2411}%
\special{pa 3525 2411}%
\special{fp}%
\special{pa 3562 2410}%
\special{pa 3570 2410}%
\special{fp}%
\special{pa 3607 2410}%
\special{pa 3615 2410}%
\special{fp}%
\special{pa 3652 2410}%
\special{pa 3660 2410}%
\special{fp}%
\special{pa 3697 2410}%
\special{pa 3705 2410}%
\special{fp}%
\special{pa 3742 2411}%
\special{pa 3750 2411}%
\special{fp}%
\special{pa 3787 2411}%
\special{pa 3791 2412}%
\special{pa 3795 2412}%
\special{fp}%
\special{pa 3832 2412}%
\special{pa 3836 2413}%
\special{pa 3840 2413}%
\special{fp}%
\special{pa 3877 2414}%
\special{pa 3885 2414}%
\special{fp}%
\special{pa 3922 2415}%
\special{pa 3927 2416}%
\special{pa 3930 2416}%
\special{fp}%
\special{pa 3967 2417}%
\special{pa 3972 2417}%
\special{pa 3975 2418}%
\special{fp}%
\special{pa 4012 2419}%
\special{pa 4017 2420}%
\special{pa 4020 2420}%
\special{fp}%
\special{pa 4057 2422}%
\special{pa 4065 2422}%
\special{fp}%
\special{pa 4102 2425}%
\special{pa 4110 2425}%
\special{fp}%
\special{pa 4147 2427}%
\special{pa 4150 2428}%
\special{pa 4155 2428}%
\special{fp}%
\special{pa 4192 2431}%
\special{pa 4199 2431}%
\special{pa 4199 2431}%
\special{fp}%
\special{pa 4236 2435}%
\special{pa 4237 2435}%
\special{pa 4244 2435}%
\special{fp}%
\special{pa 4281 2438}%
\special{pa 4285 2439}%
\special{pa 4289 2439}%
\special{fp}%
\special{pa 4326 2443}%
\special{pa 4333 2443}%
\special{pa 4333 2443}%
\special{fp}%
\special{pa 4371 2447}%
\special{pa 4375 2448}%
\special{pa 4378 2448}%
\special{fp}%
\special{pa 4415 2452}%
\special{pa 4417 2452}%
\special{pa 4421 2453}%
\special{pa 4423 2453}%
\special{fp}%
\special{pa 4460 2457}%
\special{pa 4462 2458}%
\special{pa 4467 2459}%
\special{pa 4467 2459}%
\special{fp}%
\special{pa 4504 2463}%
\special{pa 4508 2464}%
\special{pa 4512 2464}%
\special{fp}%
\special{pa 4549 2469}%
\special{pa 4552 2470}%
\special{pa 4556 2470}%
\special{fp}%
\special{pa 4593 2476}%
\special{pa 4596 2476}%
\special{pa 4600 2477}%
\special{pa 4601 2477}%
\special{fp}%
\special{pa 4637 2483}%
\special{pa 4643 2484}%
\special{pa 4645 2484}%
\special{fp}%
\special{pa 4682 2490}%
\special{pa 4682 2490}%
\special{pa 4686 2491}%
\special{pa 4689 2491}%
\special{fp}%
\special{pa 4726 2498}%
\special{pa 4731 2499}%
\special{pa 4733 2499}%
\special{fp}%
\special{pa 4770 2505}%
\special{pa 4777 2507}%
\special{fp}%
\special{pa 4814 2514}%
\special{pa 4820 2516}%
\special{pa 4821 2516}%
\special{fp}%
\special{pa 4858 2524}%
\special{pa 4858 2524}%
\special{pa 4865 2526}%
\special{fp}%
\special{pa 4902 2533}%
\special{pa 4904 2534}%
\special{pa 4907 2535}%
\special{pa 4909 2536}%
\special{fp}%
\special{pa 4945 2544}%
\special{pa 4947 2545}%
\special{pa 4951 2545}%
\special{pa 4953 2545}%
\special{fp}%
\special{pa 4989 2555}%
\special{pa 4990 2555}%
\special{pa 4993 2556}%
\special{pa 4996 2557}%
\special{fp}%
\special{pa 5032 2567}%
\special{pa 5034 2568}%
\special{pa 5038 2569}%
\special{pa 5040 2570}%
\special{fp}%
\special{pa 5076 2580}%
\special{pa 5083 2582}%
\special{fp}%
\special{pa 5119 2594}%
\special{pa 5120 2594}%
\special{pa 5124 2595}%
\special{pa 5126 2596}%
\special{fp}%
\special{pa 5161 2608}%
\special{pa 5162 2609}%
\special{pa 5168 2611}%
\special{fp}%
\special{pa 5203 2624}%
\special{pa 5207 2625}%
\special{pa 5209 2626}%
\special{pa 5210 2627}%
\special{fp}%
\special{pa 5245 2642}%
\special{pa 5246 2642}%
\special{pa 5248 2643}%
\special{pa 5251 2644}%
\special{pa 5252 2644}%
\special{fp}%
\special{pa 5286 2660}%
\special{pa 5287 2661}%
\special{pa 5289 2662}%
\special{pa 5292 2663}%
\special{pa 5292 2663}%
\special{fp}%
\special{pa 5325 2680}%
\special{pa 5331 2683}%
\special{pa 5332 2684}%
\special{fp}%
\special{pa 5363 2703}%
\special{pa 5366 2704}%
\special{pa 5368 2706}%
\special{pa 5370 2707}%
\special{fp}%
\special{pa 5400 2728}%
\special{pa 5400 2728}%
\special{pa 5403 2731}%
\special{pa 5405 2732}%
\special{pa 5406 2733}%
\special{fp}%
\special{pa 5433 2757}%
\special{pa 5434 2757}%
\special{pa 5435 2759}%
\special{pa 5437 2761}%
\special{pa 5439 2762}%
\special{fp}%
\special{pa 5462 2790}%
\special{pa 5462 2791}%
\special{pa 5465 2794}%
\special{pa 5466 2796}%
\special{pa 5466 2796}%
\special{fp}%
\special{pa 5482 2828}%
\special{pa 5482 2828}%
\special{pa 5483 2830}%
\special{pa 5483 2831}%
\special{pa 5484 2832}%
\special{pa 5484 2835}%
\special{pa 5484 2835}%
\special{fp}%
}}%
%
{\color[named]{Black}{%
\special{pn 20}%
\special{pa 3180 2760}%
\special{pa 3710 1030}%
\special{fp}%
\special{sh 1}%
\special{pa 3710 1030}%
\special{pa 3672 1088}%
\special{pa 3694 1082}%
\special{pa 3710 1100}%
\special{pa 3710 1030}%
\special{fp}%
}}%
%
{\color[named]{Black}{%
\special{pn 20}%
\special{pa 3710 2720}%
\special{pa 3490 2070}%
\special{fp}%
}}%
%
{\color[named]{Black}{%
\special{pn 20}%
\special{pa 3390 1790}%
\special{pa 3140 1080}%
\special{fp}%
\special{sh 1}%
\special{pa 3140 1080}%
\special{pa 3144 1150}%
\special{pa 3158 1130}%
\special{pa 3182 1136}%
\special{pa 3140 1080}%
\special{fp}%
}}%
\put(39.6000,-29.5000){\makebox(0,0){$(h(z_i),0)$}}%
\put(29.6000,-30.1000){\makebox(0,0){$(h(z_j),0)$}}%
\put(39.7000,-7.9000){\makebox(0,0){$(z_{j+1},1)$}}%
\put(28.2000,-8.2000){\makebox(0,0){$(z_{i+1},1)$}}%
%
{\color[named]{Black}{%
\special{pn 8}%
\special{pa 2160 2850}%
\special{pa 2160 900}%
\special{fp}%
\special{sh 1}%
\special{pa 2160 900}%
\special{pa 2140 968}%
\special{pa 2160 954}%
\special{pa 2180 968}%
\special{pa 2160 900}%
\special{fp}%
}}%
%
{\color[named]{Black}{%
\special{pn 8}%
\special{pa 2500 2900}%
\special{pa 2530 1080}%
\special{fp}%
\special{sh 1}%
\special{pa 2530 1080}%
\special{pa 2510 1146}%
\special{pa 2530 1134}%
\special{pa 2550 1148}%
\special{pa 2530 1080}%
\special{fp}%
}}%
%
{\color[named]{Black}{%
\special{pn 8}%
\special{pa 4530 2860}%
\special{pa 4530 970}%
\special{fp}%
\special{sh 1}%
\special{pa 4530 970}%
\special{pa 4510 1038}%
\special{pa 4530 1024}%
\special{pa 4550 1038}%
\special{pa 4530 970}%
\special{fp}%
}}%
%
{\color[named]{Black}{%
\special{pn 8}%
\special{pa 4960 2970}%
\special{pa 4980 1010}%
\special{fp}%
\special{sh 1}%
\special{pa 4980 1010}%
\special{pa 4960 1076}%
\special{pa 4980 1064}%
\special{pa 5000 1078}%
\special{pa 4980 1010}%
\special{fp}%
}}%
%
{\color[named]{Black}{%
\special{pn 8}%
\special{pa 4140 2720}%
\special{pa 4128 2690}%
\special{pa 4100 2630}%
\special{pa 4088 2598}%
\special{pa 4074 2568}%
\special{pa 4062 2538}%
\special{pa 4048 2508}%
\special{pa 4036 2476}%
\special{pa 4012 2416}%
\special{pa 4002 2386}%
\special{pa 3990 2354}%
\special{pa 3970 2294}%
\special{pa 3960 2262}%
\special{pa 3952 2232}%
\special{pa 3942 2202}%
\special{pa 3936 2170}%
\special{pa 3928 2140}%
\special{pa 3922 2110}%
\special{pa 3916 2078}%
\special{pa 3912 2048}%
\special{pa 3908 2016}%
\special{pa 3904 1986}%
\special{pa 3902 1954}%
\special{pa 3900 1924}%
\special{pa 3900 1862}%
\special{pa 3902 1830}%
\special{pa 3904 1800}%
\special{pa 3912 1736}%
\special{pa 3916 1706}%
\special{pa 3928 1642}%
\special{pa 3934 1612}%
\special{pa 3940 1580}%
\special{pa 3948 1548}%
\special{pa 3956 1518}%
\special{pa 3972 1454}%
\special{pa 3982 1422}%
\special{pa 3990 1392}%
\special{pa 4020 1296}%
\special{pa 4030 1266}%
\special{pa 4040 1234}%
\special{pa 4040 1230}%
\special{fp}%
}}%
%
{\color[named]{Black}{%
\special{pn 8}%
\special{pa 4040 1220}%
\special{pa 4100 1080}%
\special{fp}%
\special{sh 1}%
\special{pa 4100 1080}%
\special{pa 4056 1134}%
\special{pa 4080 1130}%
\special{pa 4092 1150}%
\special{pa 4100 1080}%
\special{fp}%
}}%
%
{\color[named]{Black}{%
\special{pn 8}%
\special{pa 2880 1290}%
\special{pa 2790 1130}%
\special{fp}%
\special{sh 1}%
\special{pa 2790 1130}%
\special{pa 2806 1198}%
\special{pa 2816 1176}%
\special{pa 2840 1178}%
\special{pa 2790 1130}%
\special{fp}%
}}%
%
{\color[named]{Black}{%
\special{pn 8}%
\special{pa 2880 1300}%
\special{pa 2896 1324}%
\special{pa 2914 1346}%
\special{pa 2946 1394}%
\special{pa 2960 1418}%
\special{pa 2976 1442}%
\special{pa 2990 1466}%
\special{pa 3018 1518}%
\special{pa 3042 1574}%
\special{pa 3054 1604}%
\special{pa 3062 1634}%
\special{pa 3072 1666}%
\special{pa 3078 1698}%
\special{pa 3084 1732}%
\special{pa 3090 1768}%
\special{pa 3094 1806}%
\special{pa 3094 1844}%
\special{pa 3096 1886}%
\special{pa 3094 1928}%
\special{pa 3086 2020}%
\special{pa 3078 2068}%
\special{pa 3070 2118}%
\special{pa 3060 2170}%
\special{pa 3050 2220}%
\special{pa 3038 2272}%
\special{pa 2996 2428}%
\special{pa 2968 2528}%
\special{pa 2952 2574}%
\special{pa 2938 2620}%
\special{pa 2924 2664}%
\special{pa 2910 2704}%
\special{pa 2898 2742}%
\special{pa 2884 2776}%
\special{pa 2874 2806}%
\special{pa 2864 2834}%
\special{pa 2854 2856}%
\special{pa 2846 2874}%
\special{pa 2840 2888}%
\special{pa 2836 2896}%
\special{pa 2834 2898}%
\special{pa 2834 2894}%
\special{pa 2838 2884}%
\special{pa 2842 2868}%
\special{pa 2860 2814}%
\special{pa 2860 2810}%
\special{fp}%
}}%
\end{picture}%

\caption{}
\end{figure}

Let $\tilde\varphi$ be the lift of $\varphi$ such that 
$\Vert\tilde\varphi-{\rm Id}\Vert_0<2\delta$. 
Now the sequence $(z_0,z_q,z_{2q},\ldots,z_{\xi q})$ is a periodic orbit
of $(\varphi\circ h)^q$. It has a lift 
$\tilde z_0,\tilde z_q,\ldots,\tilde z_{\xi q}$ that is a periodic orbit of
$(\tilde \varphi\circ \tilde h)^q\circ T^{-p_1-1}$, since the dynamical index of
$\gamma$
is $\xi(q,p_1+1)$.
Hence by the Brouwer plane fixed point theorem, there is a fixed point
of $(\tilde \varphi\circ \tilde h)^q\circ T^{-p_1-1}$. This is contrary to
condition (3). The proof is complete now. \qed

\medskip
In the rest we shall construct
a $\delta$-chain  prohibited in Lemma \ref{l1},
by using the condition $x_1\sim x_2$.
The absurdity will show that $F_0\neq\emptyset$, as is required.

Let $\gamma_1$ be a $\delta$-chain from $x_1$ to $x_2$ of
dynamical index $(i_1,j_1)$, and $\gamma_3$ another
from $x_2$ to $x_1$ of dynamical index
$(i_2,j_2)$.
One can assume that $i_1+i_2$ is a multiple of $q$. In
fact, if it is not the case, consider the concatination
$(\gamma_1\cdot\gamma_3)^{q-1}\cdot\gamma_1$ instead of $\gamma_1$,
leaving $\gamma_3$ unchanged. Thus we can set

\begin{enumerate}\addtocounter{enumi}{3}
\item $i_1+i_2=aq$ for some $a\in\N$ and $j_1+j_2=b$ ($b\in\Z$).
\end{enumerate}
Let 
$$\gamma_2=(x_2,h(x_2),\ldots, h^{q-1}(x_2),x_2),
\ \mbox{ and}$$
$$\gamma_4=(x_1,h(x_1),\ldots, h^{q-1}(x_1),x_1).$$
 They are periodic orbits, and hence $\delta$-cycles,
of dynamical indices $(q,p_2)$ and $(q,p_1)$ respectively.
Consider the concatenation
$\gamma_1\cdot \gamma_2^\eta\cdot \gamma_3\cdot\gamma_4^\zeta$ for some
$\eta,\,\zeta\in\N$. It is a $\delta$-cycle at $x_1$ of
dynamical index $(qa+q\zeta+q\eta,b+\zeta p_1+\eta p_2)$.

We shall show there are $\xi$,
$\zeta$ and $\eta$ such that the above
concatenation becomes a $\delta$-cycle of dynamical index
$\xi(q,p_1+1)$ forbidden in Lemma
\ref{l1}.
The equation for it is the following.

\begin{enumerate}\addtocounter{enumi}{4}
\item $a+\zeta+\eta=\xi$.
\item $b+\zeta p_1+\eta p_2=\xi(p_1+1)$.
\end{enumerate}
Now for  any large $\eta>0$, define $\xi$ and $\zeta$
by $$\xi=\eta(p_2-p_1)+(b-p_1a)\ \mbox{ and }\
\zeta=\eta(p_2-p_1-1)+(b-p_1a-a).$$
Then $\eta$, $\xi$ and $\zeta$ are positive integeres
satisfying (5) and (6).
A contradiction shows Proposition
\ref{p2a}. We are done with the proof of Theorem \ref{main2}.

\medskip
\noindent
{\bf 3.4.} Finally let us show Corollary \ref{c1}.
In view of Theorem \ref{main2}, we only need to show the existence of a periodic point $x_\nu\in C_0$
such that $\rot(\tilde h,x_\nu)=\alpha_\nu$. 
For this  we proceed just as in {\bf 3.3}. The assumption that 
$[\alpha_1,\alpha_2]$ is a nondegenerate interval
is necessary for Proposition \ref{l6},
which uses Lemma \ref{l4}. The proof in the present case is
exactly the same except at the last step, {\sc Case 2}. At that point, we need the following
proposition. 

\begin{proposition} \label{p3}
Suppose $C_0$ is a chain transitive class with 
$\rot(\tilde h,
C_0)=[\alpha_1,\alpha_2]$ with $\alpha_1=p/q$, $(p,q)=1$. 
Then the homeomorphism $\tilde h^q\circ T^{-p}$  admits
a fixed point in $\tilde \AAA$.
\end{proposition}

We emphasize that we have only to show the existence of the fixed point
{\em in the whole} $\tilde \AAA$, since we have followed the argument in {\bf 3.3}.
The rest of this paragraph is devoted to the proof of Proposition \ref{p3}.
The assumption $\rot(\tilde h,C_0)=[p/q,\alpha_2]$ 
implies the following.
\begin{lemma}\label{l3.4}
We have $\rot(\tilde h^q,C_0)=[p,q\alpha_2]$. 
\end{lemma}
Here $C_0$ may not be a single chain transitive class for $h^q$. But
the rotation set $\rot(\tilde h^q,C_0)$ is defined, in the same way,
as the set
of the values $\rot(\tilde h,\mu)$, where $\mu$ runs over the space of
the $h^q$-invariant probability measures
supported on $C_0$.

\medskip
\bd Clearly a $h$-invariant probability measure $\mu$ is
$h^q$-invariant and $\rot(\tilde h^q,\mu)=q\cdot\rot(\tilde h,\mu)$.
To show this, notice that
$$
\rot(\tilde h^q,\mu)=\langle\mu,\Pi_1\circ \tilde h^q-\Pi_1\rangle
=\sum_{i=0}^q\langle\mu,\,\,\Pi_1\circ\tilde h^{i+1}-\Pi_1\circ\tilde
h^i\rangle,
$$
$$\mbox{where }\
\langle\mu,\,\,\Pi_1\circ\tilde h^{i+1}-\Pi_1\circ\tilde h^i\rangle
=\langle\mu,\,\,(\Pi_1\circ\tilde h-\Pi_1)\circ h^i\rangle
$$
$$=\langle h^i_*\mu,\,\,\Pi_1\circ\tilde h-\Pi_1\rangle
=\langle \mu,\,\,\Pi_1\circ\tilde h-\Pi_1\rangle
=
\rot(\tilde h,\mu).
$$
Thus we get $$q\cdot\rot(\tilde h,C_0)\subset\rot(\tilde h^q,C_0).$$

On the other hand, given a $h^q$-invariant probability measure
$\nu$, the average
$  \hat\nu=q^{-1}\sum_{i=0}^{q-1}h^i_*\nu$
is $h$-invariant, and we have
$$
\langle\hat\nu,\,\,\Pi_1\circ \tilde h-\Pi_1\rangle
=q^{-1}\sum_{i=0}^{q-1}\langle \tilde h^i_*\nu,\,\,\Pi_1\circ \tilde h-\Pi_1\rangle
=q^{-1}\langle\nu,\,\,\Pi_1\circ\tilde h^q-\Pi\rangle,$$
showing $\rot(\tilde h,\hat\nu)=q^{-1}\cdot\rot(\tilde h^q,\nu)$.
This implies the converse inclusion $$q\cdot\rot(\tilde h,C_0)\supset\rot(\tilde
h^q,C_0).$$
\qed

\medskip
Since $p$ is an extremal point
of the rotation set $[p,q\alpha_2]$,
there is an ergodic $h^q$-invariant probability measure 
$\mu$ supported on $C_0$
such that $\rot(\tilde h^q,\mu)=p$. 
To see this, any $h^q$-invariant measure is a convex
integral of the ergodic components, and since $p$ is
extremal, almost any ergodic component has rotation number $p$. 

We use the following 
version of the Atkinson theorem (\cite{A}), whose 
proof is found at Proposition 12.1 of \cite{FH}.

\begin{proposition} \label{p4}
Suppose $T:X\to X$ is an ergodic automorphism of a probability space
 $(X,\mu)$ and let $\varphi:X\to\R$ be an integrable function with
 $\langle\mu,\varphi\rangle=0$. Let $S(n,x)=\sum_{i=0}^{n-1}\varphi(T^i(x))$.
Then for any $\varepsilon>0$ the set of $x$ such that
 $\abs{S(n,x)}<\varepsilon$ for infinitely many $n$ is a full measure
 subset of $X$.
\end{proposition}

Notice that Proposition \ref{p4} holds only
for $\R$-valued functions, and fails e.\ g.\ for $\mathbb C$-valued functions.
We apply Proposition \ref{p4} for the transformation
$h^q:C_0\to C_0$, an ergodic measure $\mu$ with $\rot(\tilde h^q,\mu)=p$, 
the function $\varphi:C_0\to\R$ defined by
 $$\varphi\circ\pi=\Pi_1\circ \tilde h^q\circ T^{-p}-\Pi_1=\Pi_1\circ\tilde
 h^q-\Pi_1-p,$$
and  $\varepsilon=1$.
Notice that the condition $\rot(\tilde h^q,\mu)=p$ is equivalent to
$\langle\mu,\varphi\rangle=0$.

Since for any $i\in\N$,
$$
\varphi\circ h^{qi}\circ\pi=\varphi\circ\pi\circ(\tilde h^q\circ T^{-p})^i
=\Pi_1\circ(\tilde h^q\circ T^{-p})^{i+1}-\Pi_1\circ(\tilde h^q\circ
T^{-p})^{i},$$
we have
$$
S(n,\cdot)\circ\pi=\Pi_1\circ(\tilde h^q\circ T^{-p})^n-\Pi_1.$$
By Proposition \ref{p4}, there is a point $x\in C_0$ 
such that $\abs{S(n,x)}<1$ for infinitely many $n\in\N$.

Then a lift $\tilde x$ of $x$ satisfies
\begin{equation}\label{3.2}
\abs{\Pi_1((\tilde h^q\circ T^{-p})^n(\tilde x))-\Pi_1(\tilde x)}<1
\end{equation}
for infinitely many $n\in\N$.
Since the orbit of $\tilde x$ is contained in $\pi^{-1}(C_0)$, a subset
in $\tilde\AAA$ bounded from above and below, (\ref{3.2}) implies that
the $\omega$-limit set of $\tilde x$ for the homeomorphism 
$\tilde h^q\circ T^{-p}$ is nonempty. 
Especially the nonwandering set of $\tilde h^q\circ T^{-p}$
is nonempty. This implies the existence of a fixed
poit of $\tilde h^q\circ T^{-p}$ by virtue of (a variant of)
the Brouwer plane fixed point
theorem (\cite{F3}). This completes the proof of Proposition \ref{p3}.

\section{Realization of a rational prime end rotation number}

The purpose of this section is to give a proof of Theorem \ref{main3}.
We assume throughout that  $\rot(\tilde h,\infty)=p/q$ for $h\in\HH$,
and that $H$ is a $C^\infty$ complete Lyapunov function defined on $\AAA$.
Let $F_1=\pi({\rm Fix}(\tilde h^q\circ
T^{-p}))$,
the Nielsen class associated to the lift $\tilde h^q\circ T^{-q}$ of $h^q$.
Our purpose is to show that $F_1$ is nonvoid.

Let $a$ be
a regular and dynamically regular value of $H$ which satisfies the
following condition:
\begin{equation}\label{e4.1}
A^+\cap \Fr(U_\infty)\neq\emptyset, 
\end{equation}
where $A^+$ is the upper subannulus bounded by 
the unique homotopically nontrivial simple closed curve in  $H^{-1}(a)$.
The lower subannulus is denoted by $A^-$.

\begin{figure}
\unitlength 0.1in
\begin{picture}( 41.9400, 29.6400)( 13.6000,-33.4400)
%
{\color[named]{Black}{%
\special{pn 8}%
\special{pa 1360 1620}%
\special{pa 1384 1596}%
\special{pa 1408 1570}%
\special{pa 1430 1544}%
\special{pa 1478 1496}%
\special{pa 1500 1470}%
\special{pa 1524 1448}%
\special{pa 1548 1424}%
\special{pa 1572 1402}%
\special{pa 1598 1380}%
\special{pa 1622 1358}%
\special{pa 1646 1338}%
\special{pa 1672 1320}%
\special{pa 1696 1302}%
\special{pa 1722 1284}%
\special{pa 1748 1268}%
\special{pa 1776 1254}%
\special{pa 1802 1240}%
\special{pa 1830 1228}%
\special{pa 1858 1218}%
\special{pa 1914 1202}%
\special{pa 1942 1196}%
\special{pa 1972 1192}%
\special{pa 2032 1188}%
\special{pa 2064 1188}%
\special{pa 2094 1190}%
\special{pa 2126 1192}%
\special{pa 2158 1196}%
\special{pa 2190 1202}%
\special{pa 2222 1210}%
\special{pa 2256 1218}%
\special{pa 2288 1226}%
\special{pa 2322 1238}%
\special{pa 2356 1248}%
\special{pa 2424 1276}%
\special{pa 2460 1290}%
\special{pa 2494 1306}%
\special{pa 2530 1324}%
\special{pa 2564 1342}%
\special{pa 2598 1362}%
\special{pa 2634 1382}%
\special{pa 2666 1402}%
\special{pa 2700 1424}%
\special{pa 2732 1446}%
\special{pa 2762 1470}%
\special{pa 2818 1518}%
\special{pa 2842 1542}%
\special{pa 2866 1568}%
\special{pa 2886 1594}%
\special{pa 2904 1620}%
\special{pa 2920 1646}%
\special{pa 2932 1674}%
\special{pa 2942 1700}%
\special{pa 2946 1728}%
\special{pa 2948 1754}%
\special{pa 2946 1782}%
\special{pa 2942 1808}%
\special{pa 2932 1836}%
\special{pa 2900 1888}%
\special{pa 2878 1914}%
\special{pa 2854 1940}%
\special{pa 2828 1964}%
\special{pa 2798 1986}%
\special{pa 2766 2008}%
\special{pa 2734 2028}%
\special{pa 2700 2046}%
\special{pa 2664 2062}%
\special{pa 2628 2076}%
\special{pa 2590 2088}%
\special{pa 2554 2098}%
\special{pa 2518 2106}%
\special{pa 2480 2110}%
\special{pa 2446 2110}%
\special{pa 2412 2108}%
\special{pa 2378 2102}%
\special{pa 2348 2094}%
\special{pa 2318 2082}%
\special{pa 2292 2064}%
\special{pa 2268 2044}%
\special{pa 2246 2020}%
\special{pa 2226 1994}%
\special{pa 2206 1966}%
\special{pa 2188 1938}%
\special{pa 2170 1908}%
\special{pa 2152 1880}%
\special{pa 2132 1852}%
\special{pa 2112 1826}%
\special{pa 2092 1802}%
\special{pa 2068 1782}%
\special{pa 2044 1764}%
\special{pa 2016 1752}%
\special{pa 1984 1744}%
\special{pa 1950 1742}%
\special{pa 1916 1744}%
\special{pa 1882 1752}%
\special{pa 1848 1764}%
\special{pa 1816 1778}%
\special{pa 1788 1798}%
\special{pa 1762 1820}%
\special{pa 1742 1846}%
\special{pa 1730 1874}%
\special{pa 1722 1906}%
\special{pa 1720 1938}%
\special{pa 1724 1972}%
\special{pa 1732 2006}%
\special{pa 1744 2040}%
\special{pa 1760 2072}%
\special{pa 1778 2102}%
\special{pa 1800 2128}%
\special{pa 1824 2152}%
\special{pa 1848 2174}%
\special{pa 1876 2190}%
\special{pa 1904 2204}%
\special{pa 1934 2216}%
\special{pa 1966 2226}%
\special{pa 1996 2234}%
\special{pa 2030 2240}%
\special{pa 2062 2246}%
\special{pa 2096 2248}%
\special{pa 2128 2252}%
\special{pa 2162 2256}%
\special{pa 2226 2260}%
\special{pa 2260 2260}%
\special{pa 2292 2262}%
\special{pa 2454 2262}%
\special{pa 2518 2258}%
\special{pa 2552 2256}%
\special{pa 2616 2252}%
\special{pa 2648 2248}%
\special{pa 2678 2246}%
\special{pa 2838 2226}%
\special{pa 2870 2220}%
\special{pa 2900 2216}%
\special{pa 2932 2210}%
\special{pa 2964 2206}%
\special{pa 2996 2200}%
\special{pa 3026 2194}%
\special{pa 3090 2182}%
\special{pa 3120 2174}%
\special{pa 3184 2162}%
\special{pa 3214 2154}%
\special{pa 3246 2148}%
\special{pa 3278 2140}%
\special{pa 3308 2134}%
\special{pa 3340 2126}%
\special{pa 3370 2118}%
\special{pa 3402 2112}%
\special{pa 3432 2104}%
\special{pa 3496 2088}%
\special{pa 3526 2082}%
\special{pa 3558 2074}%
\special{pa 3588 2066}%
\special{pa 3610 2060}%
\special{fp}%
}}%
%
{\color[named]{Black}{%
\special{pn 8}%
\special{pa 3630 1920}%
\special{pa 3570 1884}%
\special{pa 3540 1864}%
\special{pa 3510 1846}%
\special{pa 3482 1828}%
\special{pa 3422 1788}%
\special{pa 3394 1770}%
\special{pa 3282 1690}%
\special{pa 3256 1670}%
\special{pa 3204 1626}%
\special{pa 3132 1560}%
\special{pa 3066 1488}%
\special{pa 3046 1464}%
\special{pa 2992 1386}%
\special{pa 2976 1358}%
\special{pa 2948 1302}%
\special{pa 2936 1272}%
\special{pa 2916 1212}%
\special{pa 2908 1180}%
\special{pa 2900 1146}%
\special{pa 2892 1078}%
\special{pa 2888 1042}%
\special{pa 2888 970}%
\special{pa 2890 932}%
\special{pa 2894 896}%
\special{pa 2898 858}%
\special{pa 2906 820}%
\special{pa 2922 748}%
\special{pa 2932 714}%
\special{pa 2944 680}%
\special{pa 2972 616}%
\special{pa 3004 560}%
\special{pa 3022 536}%
\special{pa 3042 512}%
\special{pa 3062 492}%
\special{pa 3084 476}%
\special{pa 3106 462}%
\special{pa 3130 450}%
\special{pa 3156 442}%
\special{pa 3182 438}%
\special{pa 3208 438}%
\special{pa 3236 440}%
\special{pa 3264 444}%
\special{pa 3292 452}%
\special{pa 3322 462}%
\special{pa 3350 474}%
\special{pa 3410 506}%
\special{pa 3438 526}%
\special{pa 3468 546}%
\special{pa 3552 618}%
\special{pa 3580 646}%
\special{pa 3632 702}%
\special{pa 3656 732}%
\special{pa 3680 764}%
\special{pa 3702 794}%
\special{pa 3724 826}%
\special{pa 3760 890}%
\special{pa 3792 954}%
\special{pa 3816 1018}%
\special{pa 3826 1050}%
\special{pa 3834 1082}%
\special{pa 3852 1178}%
\special{pa 3854 1210}%
\special{pa 3858 1242}%
\special{pa 3858 1338}%
\special{pa 3854 1402}%
\special{pa 3846 1466}%
\special{pa 3822 1594}%
\special{pa 3814 1626}%
\special{pa 3808 1658}%
\special{pa 3792 1722}%
\special{pa 3782 1754}%
\special{pa 3766 1818}%
\special{pa 3756 1850}%
\special{pa 3750 1870}%
\special{fp}%
}}%
%
{\color[named]{Black}{%
\special{pn 8}%
\special{pa 3870 2070}%
\special{pa 3954 2078}%
\special{pa 3998 2082}%
\special{pa 4082 2090}%
\special{pa 4122 2094}%
\special{pa 4164 2096}%
\special{pa 4248 2104}%
\special{pa 4288 2106}%
\special{pa 4330 2108}%
\special{pa 4370 2112}%
\special{pa 4410 2114}%
\special{pa 4450 2114}%
\special{pa 4530 2118}%
\special{pa 4646 2118}%
\special{pa 4684 2116}%
\special{pa 4720 2116}%
\special{pa 4758 2114}%
\special{pa 4794 2110}%
\special{pa 4830 2108}%
\special{pa 4864 2104}%
\special{pa 4900 2100}%
\special{pa 4934 2094}%
\special{pa 4966 2090}%
\special{pa 5000 2084}%
\special{pa 5032 2076}%
\special{pa 5062 2068}%
\special{pa 5094 2060}%
\special{pa 5122 2050}%
\special{pa 5152 2040}%
\special{pa 5208 2016}%
\special{pa 5234 2004}%
\special{pa 5260 1990}%
\special{pa 5284 1976}%
\special{pa 5308 1960}%
\special{pa 5332 1942}%
\special{pa 5354 1924}%
\special{pa 5374 1906}%
\special{pa 5394 1886}%
\special{pa 5414 1864}%
\special{pa 5432 1842}%
\special{pa 5464 1794}%
\special{pa 5478 1768}%
\special{pa 5492 1740}%
\special{pa 5504 1712}%
\special{pa 5524 1652}%
\special{pa 5534 1620}%
\special{pa 5546 1552}%
\special{pa 5554 1480}%
\special{pa 5554 1408}%
\special{pa 5552 1372}%
\special{pa 5546 1336}%
\special{pa 5540 1302}%
\special{pa 5530 1270}%
\special{pa 5518 1238}%
\special{pa 5504 1206}%
\special{pa 5488 1178}%
\special{pa 5470 1152}%
\special{pa 5448 1128}%
\special{pa 5422 1106}%
\special{pa 5394 1088}%
\special{pa 5364 1074}%
\special{pa 5330 1062}%
\special{pa 5296 1054}%
\special{pa 5262 1050}%
\special{pa 5228 1052}%
\special{pa 5196 1060}%
\special{pa 5170 1074}%
\special{pa 5146 1094}%
\special{pa 5128 1120}%
\special{pa 5114 1150}%
\special{pa 5106 1184}%
\special{pa 5102 1216}%
\special{pa 5100 1250}%
\special{pa 5102 1282}%
\special{pa 5106 1314}%
\special{pa 5110 1378}%
\special{pa 5110 1408}%
\special{pa 5106 1440}%
\special{pa 5100 1472}%
\special{pa 5092 1502}%
\special{pa 5082 1534}%
\special{pa 5070 1564}%
\special{pa 5056 1594}%
\special{pa 5038 1624}%
\special{pa 5020 1652}%
\special{pa 5000 1678}%
\special{pa 4976 1704}%
\special{pa 4928 1752}%
\special{pa 4900 1774}%
\special{pa 4872 1794}%
\special{pa 4842 1812}%
\special{pa 4810 1828}%
\special{pa 4746 1856}%
\special{pa 4712 1866}%
\special{pa 4644 1882}%
\special{pa 4576 1890}%
\special{pa 4508 1890}%
\special{pa 4474 1886}%
\special{pa 4410 1874}%
\special{pa 4378 1864}%
\special{pa 4348 1852}%
\special{pa 4320 1838}%
\special{pa 4294 1822}%
\special{pa 4268 1804}%
\special{pa 4246 1784}%
\special{pa 4224 1762}%
\special{pa 4202 1738}%
\special{pa 4184 1712}%
\special{pa 4152 1656}%
\special{pa 4138 1626}%
\special{pa 4124 1594}%
\special{pa 4104 1530}%
\special{pa 4094 1496}%
\special{pa 4076 1394}%
\special{pa 4074 1360}%
\special{pa 4070 1326}%
\special{pa 4070 1224}%
\special{pa 4074 1192}%
\special{pa 4076 1160}%
\special{pa 4080 1128}%
\special{pa 4086 1096}%
\special{pa 4092 1066}%
\special{pa 4100 1034}%
\special{pa 4110 1004}%
\special{pa 4118 972}%
\special{pa 4142 912}%
\special{pa 4154 880}%
\special{pa 4170 850}%
\special{pa 4184 820}%
\special{pa 4200 790}%
\special{pa 4218 760}%
\special{pa 4254 704}%
\special{pa 4274 676}%
\special{pa 4296 650}%
\special{pa 4340 602}%
\special{pa 4388 558}%
\special{pa 4414 540}%
\special{pa 4438 522}%
\special{pa 4466 506}%
\special{pa 4522 482}%
\special{pa 4550 474}%
\special{pa 4580 466}%
\special{pa 4610 462}%
\special{pa 4674 458}%
\special{pa 4706 458}%
\special{pa 4738 460}%
\special{pa 4772 464}%
\special{pa 4808 470}%
\special{pa 4842 476}%
\special{pa 4878 484}%
\special{pa 4950 504}%
\special{pa 4988 514}%
\special{pa 5024 524}%
\special{pa 5062 536}%
\special{pa 5098 546}%
\special{pa 5134 554}%
\special{pa 5170 564}%
\special{pa 5238 576}%
\special{pa 5270 580}%
\special{pa 5332 580}%
\special{pa 5360 576}%
\special{pa 5388 570}%
\special{pa 5412 558}%
\special{pa 5434 544}%
\special{pa 5454 526}%
\special{pa 5472 504}%
\special{pa 5488 478}%
\special{pa 5502 450}%
\special{pa 5516 418}%
\special{pa 5528 386}%
\special{pa 5530 380}%
\special{fp}%
}}%
%
{\color[named]{Black}{%
\special{pn 8}%
\special{ar 3750 2020 1322 1322  3.0810606  3.6560440}%
}}%
%
{\color[named]{Black}{%
\special{pn 8}%
\special{ar 3740 2020 1336 1336  3.9971458  4.6292477}%
}}%
%
{\color[named]{Black}{%
\special{pn 8}%
\special{ar 3750 2020 1342 1342  5.0686086  6.0038621}%
}}%
%
{\color[named]{Black}{%
\special{pn 20}%
\special{ar 3750 2030 1328 1328  3.6702172  4.0061899}%
}}%
%
{\color[named]{Black}{%
\special{pn 20}%
\special{ar 3730 2020 1332 1332  4.6427344  5.0746901}%
}}%
%
{\color[named]{Black}{%
\special{pn 8}%
\special{ar 3750 2020 1324 1324  0.0227234  2.9576492}%
}}%
\put(24.0000,-4.8000){\makebox(0,0){$V$}}%
\put(25.9000,-10.9000){\makebox(0,0){$c_\nu$}}%
\put(31.7000,-19.4000){\makebox(0,0){$E_\nu$}}%
\put(38.8000,-5.2000){\makebox(0,0){$c_{\nu'}$}}%
\put(53.5000,-15.8000){\makebox(0,0){$E_{\nu'}$}}%
%
{\color[named]{Black}{%
\special{pn 8}%
\special{pa 3730 2163}%
\special{pa 3722 2163}%
\special{fp}%
\special{pa 3687 2159}%
\special{pa 3683 2159}%
\special{pa 3683 2158}%
\special{pa 3680 2158}%
\special{fp}%
\special{pa 3652 2147}%
\special{pa 3651 2147}%
\special{pa 3651 2146}%
\special{pa 3649 2146}%
\special{pa 3649 2145}%
\special{pa 3647 2145}%
\special{pa 3647 2145}%
\special{fp}%
\special{pa 3623 2128}%
\special{pa 3623 2128}%
\special{pa 3622 2128}%
\special{pa 3622 2127}%
\special{pa 3621 2127}%
\special{pa 3621 2126}%
\special{pa 3620 2126}%
\special{pa 3619 2125}%
\special{pa 3619 2124}%
\special{fp}%
\special{pa 3600 2104}%
\special{pa 3599 2104}%
\special{pa 3599 2103}%
\special{pa 3598 2103}%
\special{pa 3598 2102}%
\special{pa 3597 2102}%
\special{pa 3597 2101}%
\special{pa 3596 2100}%
\special{pa 3596 2099}%
\special{fp}%
\special{pa 3581 2074}%
\special{pa 3581 2073}%
\special{pa 3580 2073}%
\special{pa 3580 2071}%
\special{pa 3579 2070}%
\special{pa 3579 2068}%
\special{fp}%
\special{pa 3569 2037}%
\special{pa 3569 2037}%
\special{pa 3569 2032}%
\special{pa 3568 2031}%
\special{pa 3568 2030}%
\special{fp}%
\special{pa 3568 1994}%
\special{pa 3568 1989}%
\special{pa 3569 1989}%
\special{pa 3569 1987}%
\special{fp}%
\special{pa 3577 1955}%
\special{pa 3577 1954}%
\special{pa 3578 1954}%
\special{pa 3578 1952}%
\special{pa 3579 1952}%
\special{pa 3579 1950}%
\special{pa 3580 1949}%
\special{fp}%
\special{pa 3594 1924}%
\special{pa 3594 1923}%
\special{pa 3595 1922}%
\special{pa 3595 1921}%
\special{pa 3596 1921}%
\special{pa 3596 1920}%
\special{pa 3597 1920}%
\special{pa 3597 1919}%
\special{pa 3597 1919}%
\special{fp}%
\special{pa 3616 1898}%
\special{pa 3616 1897}%
\special{pa 3617 1897}%
\special{pa 3620 1894}%
\special{pa 3621 1894}%
\special{pa 3621 1893}%
\special{fp}%
\special{pa 3645 1877}%
\special{pa 3645 1877}%
\special{pa 3646 1876}%
\special{pa 3647 1876}%
\special{pa 3647 1875}%
\special{pa 3649 1875}%
\special{pa 3649 1874}%
\special{pa 3650 1874}%
\special{fp}%
\special{pa 3678 1863}%
\special{pa 3679 1863}%
\special{pa 3679 1862}%
\special{pa 3683 1862}%
\special{pa 3683 1861}%
\special{pa 3684 1861}%
\special{fp}%
\special{pa 3719 1857}%
\special{pa 3727 1857}%
\special{fp}%
\special{pa 3761 1863}%
\special{pa 3761 1863}%
\special{pa 3764 1863}%
\special{pa 3764 1864}%
\special{pa 3767 1864}%
\special{fp}%
\special{pa 3793 1876}%
\special{pa 3794 1876}%
\special{pa 3795 1877}%
\special{pa 3796 1877}%
\special{pa 3797 1878}%
\special{pa 3798 1878}%
\special{pa 3798 1879}%
\special{pa 3799 1879}%
\special{fp}%
\special{pa 3823 1896}%
\special{pa 3823 1896}%
\special{pa 3823 1897}%
\special{pa 3824 1898}%
\special{pa 3825 1898}%
\special{pa 3825 1899}%
\special{pa 3826 1899}%
\special{pa 3826 1900}%
\special{pa 3827 1900}%
\special{pa 3827 1900}%
\special{fp}%
\special{pa 3845 1922}%
\special{pa 3845 1922}%
\special{pa 3846 1923}%
\special{pa 3846 1924}%
\special{pa 3847 1924}%
\special{pa 3847 1925}%
\special{pa 3848 1926}%
\special{pa 3848 1927}%
\special{pa 3849 1927}%
\special{fp}%
\special{pa 3862 1952}%
\special{pa 3862 1954}%
\special{pa 3863 1954}%
\special{pa 3863 1957}%
\special{pa 3864 1957}%
\special{pa 3864 1958}%
\special{fp}%
\special{pa 3872 1989}%
\special{pa 3872 1997}%
\special{fp}%
\special{pa 3871 2032}%
\special{pa 3871 2037}%
\special{pa 3870 2038}%
\special{pa 3870 2039}%
\special{fp}%
\special{pa 3861 2069}%
\special{pa 3861 2070}%
\special{pa 3860 2071}%
\special{pa 3860 2073}%
\special{pa 3859 2073}%
\special{pa 3859 2075}%
\special{pa 3858 2075}%
\special{fp}%
\special{pa 3844 2100}%
\special{pa 3844 2100}%
\special{pa 3843 2101}%
\special{pa 3843 2102}%
\special{pa 3842 2102}%
\special{pa 3842 2103}%
\special{pa 3841 2103}%
\special{pa 3841 2104}%
\special{pa 3840 2105}%
\special{fp}%
\special{pa 3820 2125}%
\special{pa 3820 2126}%
\special{pa 3819 2126}%
\special{pa 3819 2127}%
\special{pa 3818 2127}%
\special{pa 3818 2128}%
\special{pa 3817 2128}%
\special{pa 3817 2129}%
\special{pa 3816 2129}%
\special{fp}%
\special{pa 3791 2145}%
\special{pa 3791 2146}%
\special{pa 3789 2146}%
\special{pa 3789 2147}%
\special{pa 3787 2147}%
\special{pa 3787 2148}%
\special{pa 3786 2148}%
\special{fp}%
\special{pa 3757 2158}%
\special{pa 3756 2159}%
\special{pa 3753 2159}%
\special{pa 3752 2160}%
\special{pa 3750 2160}%
\special{fp}%
}}%
\put(36.8000,-26.3000){\makebox(0,0){$A^-$}}%
\end{picture}%

\caption{}
\end{figure}

Let $V$ be the unique unbounded component of $U_\infty\cap A^+$.
See Figure 4. Let 
$$\Cl_{U_\infty}(V)\cap\partial A^-=\coprod_{\nu\in I}c_\nu,$$
where $c_\nu$ are cross cuts of $U_\infty$.
The cross cuts $c_\nu$ are at most countable and oriented according to
the orientation of $V$. 
Let $E_\nu$ be the connected component of $U_\infty\setminus c_\nu$
disjoint from $V$. 
Since $A^-$ is forward invariant, $\Cl(E_\nu)$ is mapped
by $h^q$ into  
some $E_{\nu'}$.

Let $p_\nu$ (resp.\  $q_\nu$) be the 
innitial point (resp.\ terminal point) of $c_\nu$. As in {\bf 2.3},
the cross cut $c_\nu$ with endpoint $p_\nu$ (resp.\ $q_\nu$) defines
a prime end denoted by $\hat p_\nu$ (resp.\ $\hat q_\nu$).
Denote by $\hat c_\nu$ the closed interval in the
set of prime ends $\partial U_\infty^*$ 
bounded by $\hat p_\nu$ and
$\hat q_\nu$. 
In other words,
$$\hat c_\nu=\Cl_{U_\infty^*}(E_\nu)\cap\partial U_\infty^*.$$
Of course they are mutually disjoint, and $\hat c_\nu$ is mapped
by $(h^*_\infty)^q$ into the interior of 
some $\hat c_{\nu'}$. If $\nu\neq\nu'$, then there 
is no fixed point of $(h^*_\infty)^q$ in $\hat c_\nu$. If $\nu=\nu'$, 
$\hat c_\nu$ is mapped into the interior of itself
by $\hat c_\nu$.
On the other hand, there is a fixed point of
$(h^*_\infty)^q$, since $\rot(\tilde h,\infty)=p/q$.
Therefore there must be a fixed point $\xi$ of $(h_\infty^*)^q$ in the set
$$\Xi=\partial U^*_\infty\setminus(\bigcup_\nu \hat
c_\nu).$$
More precisely, any lift $\tilde\xi$ of $\xi$ to $\partial \tilde
U_\infty^*$,
the universal cover of $\partial U_\infty^*$, is fixed 
by $(\tilde h_\infty^*)^q\circ T^{-p}$.

The {\em principal point set} $\Pi(\xi)$ of the prime end
$\xi$ is defined to be the
set of all the limit points of topological chains which represent 
the prime end $\xi$. As is well known (\cite{P}), the principal point
set $\Pi(\xi)$ is closed, connected and invariant by $h^q$. Clearly it
is contained in $\Fr(U_\infty)$. Also since $\xi$ is contained in
$\Xi$,  any cross cut $c_i$ of any topological chain $\{c_i\}$
representing $\xi$ must intersect $\Cl(V)$.
The set $\Pi(\xi)$ is contained in $\Cl(V)$,
since ${\rm diam}(c_i)\to0$.  This implies that
$\Pi(\xi)$ is compact.
Let $\hat\Pi(\xi)$ be the union of $\Pi(\xi)$ with all the bounded connected
components of the complement. The set $\hat\Pi(\xi)$ is also a
$h^q$-invariant continuum, and therefore it
 does not separate two ends of $\AAA$ by the assumption on $h$.
It is also nonseparating, in the sense that its complement is
connected.

The Cartwright-Littlewood
theorem (\cite{CL}) asserts that any planar homeomorphism leaving a
nonseparating continuum invariant has a fixed point in it.
Thus there is a fixed point $y$ of $h^q$ in $\hat\Pi(\xi)$.
In the rest of this section, we shall
 show  $y\in F_1$, i.\ e.\ a lift $\tilde y$ of $y$ is
a fixed point of $\tilde h^q\circ T^{-p}$. But in fact, we shall
find such a point $\tilde y$ at the very end of the proof.

Recall that for a bounded cross cut  $c$ of $U_\infty$,
$V(c)$ denotes the component of $U_\infty\setminus c$ which
is homeomorphic to an open disc.
Likewise we define the component $V(\tilde c)$ for a lift $\tilde c$
of $c$ to be  the lift of $V(c)$ bounded by $\tilde c$.

Given a topological chain $\{c_i\}$ of $U_\infty$, a {\em lift}
$\{\tilde c_i\}$  of $\{c_i\}$ is defined as follows.  For $i=1$,
let $\tilde c_1$ be an arbitrary lift of $c_1$. 
For $i>1$, let $\tilde c_i$ be the unique lift of $c_i$ contained
in $V(\tilde c_{i-1})$. Then we have $V(\tilde c_i)\subset V(\tilde c_{i-1})$
($\forall i>1$), and the lift $\{\tilde c_i\}$ is determined uniquely by the
choice of $\tilde c_1$.

Let $x\in\Pi(\xi)$ be an arbitrary point, and let 
$\{c_i\}$ be a topological chain representing $\xi$
such that $c_i\to x$. Let $\{\tilde c_i\}$ be a lift
of $\{c_i\}$ and $\tilde x$ a lift of $x$.
Then since $c_i\to x$,
there is a sequence of integers $n_i$ such that 
$T^{n_i}(\tilde c_i)\to \tilde x$. Let us show that $n_i$ is
identical for any large $i$. 

Since $\Pi(\xi)$ is compact
and does not separate the two ends of $\AAA$, there is a simple
closed curve $\Gamma$ such that $\Pi(\xi)$ is contained in
the open disc $E$ bounded by $\Gamma$. 
Assume that there are infinitely many $i$ such that $n_{i+1}\neq n_i$.
For any large $i$, the cross cuts $c_i$ and $c_{i+1}$ are contained in
$E$. Consider a simple path $\gamma$ joining $c_i$ to $c_{i+1}$ 
in $V(c_i)\setminus V(c_{i+1})$.
Then $\gamma$, starting and ending in $E$, must wind the annulus $\AAA$
since $n_{i+1}\neq n_i$.
Thus
there is a cross cut $c'_i$ contained in $\Gamma$ which separates
$c_{i+1}$ and $c_i$. 
Passing to a further subsequence, we may assume $\Cl(c_i')$ are
disjoint, since $c_i'$ are disjoint open intervals of a single
curve $\Gamma$. 
We also have ${\rm diam}(c'_i)\to0$. Thus $\{c_i'\}$
is a topological chain contained in $\Gamma$, 
which is clearly equivalent to $\{c_i\}$. 
Thus any accumulation point of $\{c'_i\}$ must be contained in the
principal point set $\Pi(\xi)$.
This contradicts the choice
of $\Gamma$: $\Gamma\cap\Pi(\xi)=\emptyset$.

Now one can assume, changing the lift $\tilde x$ of
$x$ if necessary,  that $\tilde c_i\to\tilde x$
for the lift $\{\tilde c_i\}$. 
By the definition of the topology of the prime end
compactification (Section 2),
the family $\{V(c_i)\}$ forms a fundamental neighbourhood system
of the prime end $\xi\in U_\infty^*$.
Then it follows immediately that
$\{V(\tilde c_i)\}$ forms a fundamental neighbourhood
system of a lift $\tilde \xi$ of $\xi$.
On the other hand, we have 
$(\tilde h^*_\infty)^q\circ T^{-p}(\tilde\xi)=\tilde\xi$. Thus
$\{\tilde h^q\circ T^{-p}(V(\tilde c_i))\}$ 
is also a fundamental neighbourhood system of $\tilde\xi$.
That is, $\{\tilde h^q\circ T^{-p}(\tilde c_i)\}$ and
$\{\tilde c_i\}$ are equivalent in the sense that for any $i$, there is
$j$ such that $\tilde c_j\subset V(\tilde h^q\circ T^{-p}(\tilde c_i))$ and
$\tilde h^q\circ T^{-p}(\tilde c_j)\subset V(\tilde c_i)$.

Let $\tilde \Pi(\xi)$ be the lift of $\Pi(\xi)$ which contains 
the point $\tilde x$.
The set $\tilde \Pi(\xi)$ is characterized as the set of the limit
points of lifts of topological chains which are equivalent
to $\{\tilde c_i\}$. Since
$\{\tilde h^q\circ T^{-p}(\tilde c_i)\}$ is equivalent to
$\{\tilde c_i\}$, and 
$\tilde h^q\circ T^{-p}(\tilde c_i)\to\tilde h^q\circ T^{-p}(\tilde x)$,
we have 
$\tilde h^q\circ T^{-p}(\tilde x)\in
\tilde\Pi(\xi)$.
But then since $h^q(\Pi(\xi))=\Pi(\xi)$, we have
$\tilde h^q\circ T^{-p}(\tilde\Pi(\xi))=\tilde\Pi(\xi)$.

Finally by the Cartwright-Littlewood theorem, there is
a fixed point $\tilde y$ of $\tilde h^q\circ T^{-p}$
in the corresponding lift of $\hat\Pi(\xi)$,
completing the proof of Theorem \ref{main3}.

\section{Accessible case}
 
This section is devoted to the proof of Theorem \ref{main4}.
Let $h\in\HH$ be a homeomorphism satisfying $\rot(\tilde h,\infty)=\alpha$
for some lift $\tilde h$ and $\alpha\in\R$ such that $-\infty$ is
accessible
from $U_\infty$.
By changing the coordinates of $\AAA$, one may assume that  $h$ satifies
$$
h(\theta,t)=(\theta,t-1),\ \ \forall(\theta,t)\in B,
$$
where $B=\{(\theta,t)\in\AAA\mid t\leq 0\}$.
Clearly $B\subset U_{-\infty}$.
Let
$$Z=\AAA\setminus(U_\infty\cup U_{-\infty}).$$
We shall show that $\lim_{i\to\infty}i^{-1}\Pi_1(\tilde h^i(z))=\alpha$
for any $z\in\pi^{-1}(Z)$. Clearly this implies (1) of Theorem
\ref{main4}. 

Let $V$ be the unbounded component of $U_\infty\cap(\AAA\setminus B)$.
{\color{red}Notice that $V\subset h(V)$.}
It is an essential open subannulus of $\AAA$.
Let $\{c_\nu\}$ be the family of cross cuts of $U_\infty$ contained
in $\partial B\cap\Cl(V)$ and let $V_\nu$ be the connected component of
$U_\infty\setminus c_\nu$ which is disjoint from $V$.
The components $V_\nu$ are mutually disjoint
open discs, which may intersect $\AAA\setminus B$.
The cross cut $c_\nu$ is called the {\em gate} of $V_\nu$.
{\color{red}Since $U_\infty=h(U_\infty)$ and $V\subset h(V)$, we have
$h(\cup_\nu V_\nu)\subset \cup_\nu V_\nu$.}

A component $V_\nu$ is said to be {\em accessible} if $-\infty$ is
accessible from $V_\nu$. 
This means that there is a path $\gamma:(-\infty,0]\to V_\nu$ such that
$\Pi_2\circ\gamma(t)\to-\infty$ as $t\to-\infty$, where $\Pi_2:\AAA\to\R$ is
the projection onto the second factor (the hight function).
There is an accessible component by the assumption. 
For any $V_\nu$, there exists $V_{\nu'}$ such that
$h(V_\nu)\subset V_{\nu'}$, and if $V_\nu$ is accessible, so is
$V_{\nu'}$.

Choose a sequence $V_i$ ($i\in\N$) from the family $\{V_\nu\}$ as
follows. Let $V_1$ be any accessible component. For $i>1$, let
$V_i$ be the component such that $h(V_{i-1})\subset V_i$. Then any $V_i$ is
accessible. 
The sequence $\{V_i\}$ may be all distinct or eventually periodic,
that is, there is $p\in\N$ such that $V_{i+p}=V_i$ for any large $i$.

To the gate $c_i$ of $V_i$ is associated a closed interval
$\hat c_i$ in the set of prime ends $\partial U_\infty^*$,
defined by $$\hat c_i=\Cl_{U_\infty^*}(V_i)\cap\partial U_\infty^*.$$
Since $h(V_{i-1})\subset V_i$, we have $h^*_\infty(\hat c_{i-1})\subset
\hat c_i$.
The cyclic orders of the family $\{c_i\}$ in $\partial B$ and
$\{\hat c_i\}$ in $\partial U_\infty^*$ are the same, and there is
a homeomorphism $\varphi:\partial B\to\partial U_\infty^*$ such that
$\varphi(\Cl(c_i))=\hat c_i$ ($\forall i$).

Fix once and for all a lift $\tilde V_i$ of $V_i$ to $\tilde\AAA$ in
the following way. Let $\tilde V_1$ be any lift of $V_1$,
and for $i>1$, $\tilde V_i$  the unique lift of $V_i$
which satisfies $\tilde h(\tilde V_{i-1})\subset \tilde V_i$ 
for the prescribed lift
$\tilde h$.
The gate of $\tilde V_i$ is denoted by $\tilde c_i$, that is,
$\tilde c_i$ is the frontier of $\tilde V_i$ in $\pi^{-1}(U_\infty)$.
It is a lift of $c_i$.
A closed interval $ \tilde \hat c_i$ of 
$\partial \tilde U_\infty^*=\pi^{-1}(\partial  U_\infty^*)$
is defined  by 
$$\tilde\hat{c_i}=\Cl( \tilde V_i)\cap\partial \tilde U_\infty^*.$$
 It is a lift of $\hat
c_i$, and the map $\tilde h_\infty^*$ defined on $\partial \tilde
U_\infty^*$ as an extension of $\tilde h$,
satisfy $\tilde h_\infty^*(\,\tilde\hat c_{i-1})\subset \tilde\hat c_i$.

Denote by $T$ the generator of the covering transformations of
both $\tilde \AAA$ and $\partial \tilde U_\infty^*$. There is a lift 
$$\tilde \varphi:\pi^{-1}(\partial B)\to \partial\tilde U_\infty^*$$
of $\varphi$ such that $\tilde \varphi(T^j(\Cl(\tilde c_i))= T^j(\,\tilde\hat c_i)$ 
($\forall i\in\N,\ \forall j\in\Z$). We identify $\partial\tilde U_\infty^*$
with $\pi^{-1}(\partial B)$ by $\tilde \varphi^{-1}$, and then with
$\R$ by $\Pi_1$. Thus $T$ is the right translation by $1$.

Let us denote the interval $\tilde\hat c_i=[a_i,b_i]$, where $a_i$ and
$b_i$ are real numbers by the above identification.
Recall that $\alpha=\rot(\tilde h,\infty)$ is, by definition, the rotation number
of 
$\tilde h_\infty^*:\partial\tilde U_\infty^*\to\partial\tilde
U_\infty^*$. 
Since  $\tilde h_\infty^*(\,\tilde\hat c_{i-1})\subset \tilde\hat c_i$
and the length of each $\tilde\hat c_i$ is always less than 1,
we have 
\begin{equation}\label{e5.1}
\alpha=\lim_{i\to\infty}i^{-1}a_i.
\end{equation}
Below we consider  $(a_i,b_i)$ to be the interval 
$\tilde c_i\subset\pi^{-1}(\partial B)$ by the above identification.
It is important that (\ref{e5.1}) still holds.

Our aim is to show that $\lim_{i\to\infty}\Pi_1\circ\tilde
h^i(z)=\alpha$
for any $z\in\pi^{-1}(Z)$.
But we shall show only 
$\lim_{i\to\infty}\Pi_1\circ\tilde
h^i(z)\leq\alpha$, the other inequality being shown similarly.

Let us denote by $\Gamma_i$ the set of all the simple curves $l:\R \to
\tilde U_\infty$ such that 
\begin{enumerate}
\item $\Pi_2\circ l(t)\to\pm\infty$ as $t\to\pm\infty$,  and
\item $l(t)\in \tilde V_i$ for all negative $t$.
\end{enumerate}
Since $\tilde V_i$ is the lift of an accessible component, $\Gamma_i$ is
nonempty for any $i\in\N$.

\begin{definition}
Let $z\in\pi^{-1}(Z)$. We say $z\leq \tilde V_i$ if there is $l\in\Gamma_i$
such that $z$ lies on the left side of $l$. 
\end{definition}
See Figure 5.

\begin{figure}
\unitlength 0.1in
\begin{picture}( 44.6000, 23.6100)( 12.4000,-24.9000)
%
{\color[named]{Black}{%
\special{pn 8}%
\special{pa 1900 670}%
\special{pa 1884 704}%
\special{pa 1848 768}%
\special{pa 1832 802}%
\special{pa 1814 834}%
\special{pa 1798 866}%
\special{pa 1780 900}%
\special{pa 1748 964}%
\special{pa 1730 998}%
\special{pa 1716 1030}%
\special{pa 1684 1094}%
\special{pa 1670 1128}%
\special{pa 1654 1160}%
\special{pa 1626 1224}%
\special{pa 1614 1256}%
\special{pa 1600 1288}%
\special{pa 1576 1352}%
\special{pa 1536 1480}%
\special{pa 1512 1576}%
\special{pa 1506 1606}%
\special{pa 1500 1638}%
\special{pa 1496 1670}%
\special{pa 1492 1700}%
\special{pa 1488 1732}%
\special{pa 1486 1762}%
\special{pa 1486 1856}%
\special{pa 1490 1916}%
\special{pa 1498 1976}%
\special{pa 1504 2008}%
\special{pa 1512 2036}%
\special{pa 1520 2066}%
\special{pa 1540 2126}%
\special{pa 1552 2156}%
\special{pa 1566 2184}%
\special{pa 1580 2214}%
\special{pa 1596 2244}%
\special{pa 1632 2300}%
\special{pa 1652 2328}%
\special{pa 1672 2352}%
\special{pa 1696 2376}%
\special{pa 1718 2396}%
\special{pa 1742 2412}%
\special{pa 1768 2424}%
\special{pa 1794 2434}%
\special{pa 1822 2438}%
\special{pa 1850 2440}%
\special{pa 1878 2438}%
\special{pa 1908 2434}%
\special{pa 1938 2426}%
\special{pa 1968 2416}%
\special{pa 1998 2404}%
\special{pa 2058 2372}%
\special{pa 2090 2352}%
\special{pa 2120 2332}%
\special{pa 2150 2310}%
\special{pa 2182 2286}%
\special{pa 2212 2262}%
\special{pa 2240 2236}%
\special{pa 2270 2210}%
\special{pa 2298 2184}%
\special{pa 2326 2156}%
\special{pa 2352 2128}%
\special{pa 2378 2102}%
\special{pa 2404 2074}%
\special{pa 2428 2046}%
\special{pa 2452 2020}%
\special{pa 2476 1992}%
\special{pa 2498 1966}%
\special{pa 2518 1938}%
\special{pa 2540 1912}%
\special{pa 2560 1886}%
\special{pa 2580 1858}%
\special{pa 2634 1780}%
\special{pa 2652 1752}%
\special{pa 2668 1726}%
\special{pa 2686 1700}%
\special{pa 2718 1648}%
\special{pa 2732 1622}%
\special{pa 2748 1596}%
\special{pa 2762 1570}%
\special{pa 2778 1544}%
\special{pa 2834 1440}%
\special{pa 2850 1414}%
\special{pa 2948 1232}%
\special{pa 2964 1206}%
\special{pa 2978 1182}%
\special{pa 3058 1052}%
\special{pa 3076 1026}%
\special{pa 3092 1000}%
\special{pa 3110 974}%
\special{pa 3130 950}%
\special{pa 3148 924}%
\special{pa 3188 872}%
\special{pa 3210 846}%
\special{pa 3230 820}%
\special{pa 3254 794}%
\special{pa 3276 768}%
\special{pa 3324 716}%
\special{pa 3376 664}%
\special{pa 3432 612}%
\special{pa 3460 588}%
\special{pa 3488 570}%
\special{pa 3516 560}%
\special{pa 3540 564}%
\special{pa 3564 578}%
\special{pa 3586 604}%
\special{pa 3606 634}%
\special{pa 3624 668}%
\special{pa 3642 700}%
\special{pa 3678 768}%
\special{pa 3694 800}%
\special{pa 3712 834}%
\special{pa 3726 866}%
\special{pa 3742 898}%
\special{pa 3756 932}%
\special{pa 3798 1028}%
\special{pa 3846 1156}%
\special{pa 3886 1284}%
\special{pa 3894 1316}%
\special{pa 3904 1348}%
\special{pa 3912 1378}%
\special{pa 3920 1410}%
\special{pa 3926 1442}%
\special{pa 3934 1472}%
\special{pa 3942 1504}%
\special{pa 3948 1534}%
\special{pa 3954 1566}%
\special{pa 3960 1596}%
\special{pa 3966 1628}%
\special{pa 3970 1658}%
\special{pa 3976 1688}%
\special{pa 3980 1720}%
\special{pa 3988 1780}%
\special{pa 3992 1812}%
\special{pa 4004 1902}%
\special{pa 4006 1934}%
\special{pa 4010 1964}%
\special{pa 4014 2024}%
\special{pa 4018 2054}%
\special{pa 4020 2084}%
\special{pa 4022 2116}%
\special{pa 4028 2206}%
\special{pa 4028 2236}%
\special{pa 4034 2326}%
\special{pa 4034 2356}%
\special{pa 4038 2416}%
\special{pa 4038 2446}%
\special{pa 4040 2476}%
\special{pa 4040 2490}%
\special{fp}%
}}%
%
{\color[named]{Black}{%
\special{pn 8}%
\special{pa 4220 2490}%
\special{pa 4232 2400}%
\special{pa 4236 2368}%
\special{pa 4240 2338}%
\special{pa 4242 2308}%
\special{pa 4246 2278}%
\special{pa 4250 2246}%
\special{pa 4254 2216}%
\special{pa 4256 2186}%
\special{pa 4260 2154}%
\special{pa 4262 2124}%
\special{pa 4266 2092}%
\special{pa 4268 2062}%
\special{pa 4272 2032}%
\special{pa 4274 2000}%
\special{pa 4276 1970}%
\special{pa 4278 1938}%
\special{pa 4282 1906}%
\special{pa 4282 1876}%
\special{pa 4288 1780}%
\special{pa 4288 1748}%
\special{pa 4290 1716}%
\special{pa 4290 1524}%
\special{pa 4288 1490}%
\special{pa 4288 1458}%
\special{pa 4284 1390}%
\special{pa 4282 1358}%
\special{pa 4280 1324}%
\special{pa 4276 1290}%
\special{pa 4274 1256}%
\special{pa 4266 1188}%
\special{pa 4260 1152}%
\special{pa 4256 1118}%
\special{pa 4250 1084}%
\special{pa 4246 1048}%
\special{pa 4240 1012}%
\special{pa 4234 978}%
\special{pa 4228 942}%
\special{pa 4224 908}%
\special{pa 4220 872}%
\special{pa 4218 838}%
\special{pa 4218 772}%
\special{pa 4226 708}%
\special{pa 4236 678}%
\special{pa 4246 650}%
\special{pa 4260 622}%
\special{pa 4278 596}%
\special{pa 4300 570}%
\special{pa 4326 548}%
\special{pa 4354 526}%
\special{pa 4386 510}%
\special{pa 4418 498}%
\special{pa 4450 496}%
\special{pa 4478 500}%
\special{pa 4506 514}%
\special{pa 4528 536}%
\special{pa 4550 564}%
\special{pa 4568 592}%
\special{pa 4582 620}%
\special{pa 4596 650}%
\special{pa 4608 680}%
\special{pa 4620 708}%
\special{pa 4640 768}%
\special{pa 4650 796}%
\special{pa 4674 886}%
\special{pa 4680 916}%
\special{pa 4688 948}%
\special{pa 4694 978}%
\special{pa 4698 1008}%
\special{pa 4704 1040}%
\special{pa 4712 1100}%
\special{pa 4720 1164}%
\special{pa 4722 1194}%
\special{pa 4726 1258}%
\special{pa 4728 1288}%
\special{pa 4730 1320}%
\special{pa 4730 1480}%
\special{pa 4728 1512}%
\special{pa 4728 1544}%
\special{pa 4726 1576}%
\special{pa 4724 1610}%
\special{pa 4718 1706}%
\special{pa 4714 1740}%
\special{pa 4712 1772}%
\special{pa 4708 1804}%
\special{pa 4704 1838}%
\special{pa 4700 1870}%
\special{pa 4696 1904}%
\special{pa 4692 1936}%
\special{pa 4688 1970}%
\special{pa 4684 2002}%
\special{pa 4680 2036}%
\special{pa 4674 2068}%
\special{pa 4670 2102}%
\special{pa 4664 2136}%
\special{pa 4660 2168}%
\special{pa 4654 2202}%
\special{pa 4650 2234}%
\special{pa 4638 2302}%
\special{pa 4634 2334}%
\special{pa 4622 2402}%
\special{pa 4616 2434}%
\special{pa 4610 2468}%
\special{pa 4610 2470}%
\special{fp}%
}}%
%
{\color[named]{Black}{%
\special{pn 8}%
\special{pa 4790 2450}%
\special{pa 4806 2390}%
\special{pa 4812 2360}%
\special{pa 4820 2330}%
\special{pa 4826 2298}%
\special{pa 4834 2268}%
\special{pa 4840 2238}%
\special{pa 4848 2208}%
\special{pa 4856 2176}%
\special{pa 4868 2116}%
\special{pa 4876 2086}%
\special{pa 4882 2054}%
\special{pa 4890 2024}%
\special{pa 4896 1992}%
\special{pa 4902 1962}%
\special{pa 4908 1930}%
\special{pa 4914 1900}%
\special{pa 4922 1868}%
\special{pa 4928 1836}%
\special{pa 4934 1806}%
\special{pa 4938 1774}%
\special{pa 4956 1678}%
\special{pa 4960 1646}%
\special{pa 4966 1614}%
\special{pa 4974 1550}%
\special{pa 4980 1516}%
\special{pa 4984 1484}%
\special{pa 4988 1450}%
\special{pa 4992 1418}%
\special{pa 4996 1384}%
\special{pa 4998 1352}%
\special{pa 5002 1318}%
\special{pa 5004 1284}%
\special{pa 5008 1250}%
\special{pa 5014 1148}%
\special{pa 5014 1012}%
\special{pa 5010 944}%
\special{pa 5008 912}%
\special{pa 5004 880}%
\special{pa 4998 848}%
\special{pa 4994 816}%
\special{pa 4986 784}%
\special{pa 4980 754}%
\special{pa 4950 664}%
\special{pa 4938 636}%
\special{pa 4924 608}%
\special{pa 4910 582}%
\special{pa 4894 556}%
\special{pa 4876 532}%
\special{pa 4858 506}%
\special{pa 4818 462}%
\special{pa 4794 440}%
\special{pa 4772 418}%
\special{pa 4746 398}%
\special{pa 4722 380}%
\special{pa 4694 360}%
\special{pa 4668 342}%
\special{pa 4640 326}%
\special{pa 4610 310}%
\special{pa 4582 294}%
\special{pa 4550 280}%
\special{pa 4520 266}%
\special{pa 4488 252}%
\special{pa 4458 240}%
\special{pa 4426 228}%
\special{pa 4392 216}%
\special{pa 4328 196}%
\special{pa 4294 188}%
\special{pa 4262 180}%
\special{pa 4194 164}%
\special{pa 4162 158}%
\special{pa 4094 146}%
\special{pa 4060 142}%
\special{pa 4028 138}%
\special{pa 3994 136}%
\special{pa 3960 132}%
\special{pa 3926 130}%
\special{pa 3794 130}%
\special{pa 3730 134}%
\special{pa 3602 150}%
\special{pa 3572 156}%
\special{pa 3540 164}%
\special{pa 3510 170}%
\special{pa 3450 186}%
\special{pa 3420 196}%
\special{pa 3392 206}%
\special{pa 3364 218}%
\special{pa 3336 228}%
\special{pa 3308 240}%
\special{pa 3280 254}%
\special{pa 3254 266}%
\special{pa 3228 282}%
\special{pa 3202 296}%
\special{pa 3178 312}%
\special{pa 3152 328}%
\special{pa 3128 344}%
\special{pa 3104 362}%
\special{pa 3082 380}%
\special{pa 3058 398}%
\special{pa 2992 458}%
\special{pa 2928 522}%
\special{pa 2908 546}%
\special{pa 2888 568}%
\special{pa 2868 592}%
\special{pa 2850 616}%
\special{pa 2830 640}%
\special{pa 2812 666}%
\special{pa 2792 690}%
\special{pa 2774 716}%
\special{pa 2756 744}%
\special{pa 2740 770}%
\special{pa 2722 798}%
\special{pa 2704 824}%
\special{pa 2688 852}%
\special{pa 2672 882}%
\special{pa 2654 910}%
\special{pa 2638 938}%
\special{pa 2622 968}%
\special{pa 2608 998}%
\special{pa 2576 1058}%
\special{pa 2560 1090}%
\special{pa 2546 1120}%
\special{pa 2532 1152}%
\special{pa 2516 1184}%
\special{pa 2446 1344}%
\special{pa 2418 1412}%
\special{pa 2404 1444}%
\special{pa 2362 1546}%
\special{pa 2350 1580}%
\special{pa 2308 1682}%
\special{pa 2296 1716}%
\special{pa 2282 1750}%
\special{pa 2268 1786}%
\special{pa 2256 1820}%
\special{pa 2228 1888}%
\special{pa 2196 1952}%
\special{pa 2180 1982}%
\special{pa 2162 2008}%
\special{pa 2142 2034}%
\special{pa 2120 2056}%
\special{pa 2098 2076}%
\special{pa 2072 2092}%
\special{pa 2044 2104}%
\special{pa 2014 2112}%
\special{pa 1982 2116}%
\special{pa 1946 2116}%
\special{pa 1912 2112}%
\special{pa 1878 2104}%
\special{pa 1844 2094}%
\special{pa 1812 2080}%
\special{pa 1784 2062}%
\special{pa 1760 2042}%
\special{pa 1740 2020}%
\special{pa 1726 1996}%
\special{pa 1718 1970}%
\special{pa 1714 1940}%
\special{pa 1716 1910}%
\special{pa 1720 1880}%
\special{pa 1728 1848}%
\special{pa 1752 1780}%
\special{pa 1766 1746}%
\special{pa 1782 1710}%
\special{pa 1846 1574}%
\special{pa 1862 1542}%
\special{pa 1878 1508}%
\special{pa 1892 1476}%
\special{pa 1908 1444}%
\special{pa 1922 1414}%
\special{pa 1936 1382}%
\special{pa 1952 1352}%
\special{pa 1980 1292}%
\special{pa 1996 1262}%
\special{pa 2010 1232}%
\special{pa 2024 1204}%
\special{pa 2040 1176}%
\special{pa 2054 1148}%
\special{pa 2070 1120}%
\special{pa 2086 1094}%
\special{pa 2100 1066}%
\special{pa 2116 1040}%
\special{pa 2134 1014}%
\special{pa 2150 990}%
\special{pa 2166 964}%
\special{pa 2184 940}%
\special{pa 2202 914}%
\special{pa 2238 866}%
\special{pa 2256 844}%
\special{pa 2276 820}%
\special{pa 2296 798}%
\special{pa 2318 776}%
\special{pa 2338 754}%
\special{pa 2360 732}%
\special{pa 2384 712}%
\special{pa 2406 690}%
\special{pa 2430 670}%
\special{pa 2456 650}%
\special{pa 2480 630}%
\special{fp}%
}}%
%
{\color[named]{Black}{%
\special{pn 8}%
\special{pa 2910 1310}%
\special{pa 3900 1310}%
\special{fp}%
}}%
%
{\color[named]{Black}{%
\special{pn 8}%
\special{pa 4280 1310}%
\special{pa 4730 1310}%
\special{fp}%
}}%
%
{\color[named]{Black}{%
\special{pn 8}%
\special{pa 1240 1310}%
\special{pa 1580 1310}%
\special{fp}%
}}%
%
{\color[named]{Black}{%
\special{pn 8}%
\special{pa 1970 1310}%
\special{pa 2470 1310}%
\special{fp}%
}}%
%
{\color[named]{Black}{%
\special{pn 8}%
\special{pa 5010 1310}%
\special{pa 5700 1310}%
\special{fp}%
}}%
%
{\color[named]{Black}{%
\special{pn 20}%
\special{pa 1580 1310}%
\special{pa 1970 1300}%
\special{fp}%
}}%
\put(38.2000,-4.0000){\makebox(0,0){$V_i$}}%
\put(21.8000,-4.7000){\makebox(0,0){$V$}}%
%
{\color[named]{Black}{%
\special{pn 4}%
\special{sh 1}%
\special{ar 2220 1170 8 8 0  6.28318530717959E+0000}%
\special{sh 1}%
\special{ar 3400 1170 8 8 0  6.28318530717959E+0000}%
\special{sh 1}%
\special{ar 4480 1180 8 8 0  6.28318530717959E+0000}%
\special{sh 1}%
\special{ar 5530 1630 8 8 0  6.28318530717959E+0000}%
}}%
\put(23.5000,-9.9000){\makebox(0,0){$z_3$}}%
\put(44.6000,-9.7000){\makebox(0,0){$z_2$}}%
\put(34.1000,-9.9000){\makebox(0,0){$z_1$}}%
\end{picture}%

\caption{$z_1\leq V_i$, $z_2\leq V_i$, $z_3\not\leq V_i$}
\end{figure}
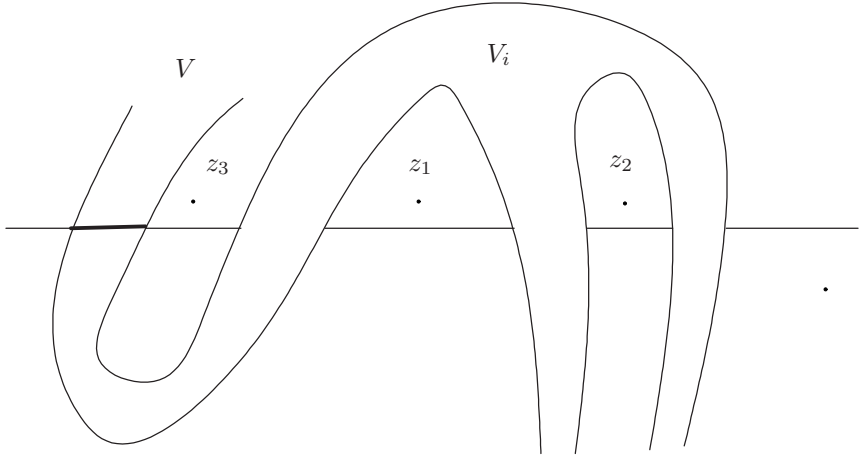

\begin{lemma}\label{l5.1}
If $z\leq \tilde V_{i-1}$ for $z\in\pi^{-1}(Z)$ and $i>1$, then
$\tilde h(z)\leq \tilde V_i$.
\end{lemma}

\bd If $l\in\Gamma_{i-1}$, then $\tilde h(l)\in\Gamma_i$. The lemma
follows from this. \qed

\begin{lemma}\label{l5.2}
There is $M>0$ such that if $z\leq \tilde V_i$
{\em ($z\in\pi^{-1}(Z)$)}, then $\Pi_1(z)\leq a_i+M$.
\end{lemma}

\bd 
We shall show the following.
\begin{enumerate}
\item There is $M>0$ such that if $z\leq\tilde V_1$ ($z\in\pi^{-1}(Z)$),
then $\Pi_1(z)\leq a_1+M-1$.
\end{enumerate}
Let us explain why this is sufficient.
Considering the action of covering transformations, (1)
implies the following.
\begin{enumerate}\addtocounter{enumi}{1}
\item  If $z\leq T^n(\tilde V_1)$ ($z\in\pi^{-1}(Z)$, $n\in\Z$)
under the similar definition, 
then $\Pi_1(z)\leq a_1+n+M-1$.
\end{enumerate}
To deduce the lemma from (2),
let $n$ be the integer such that $a_1+n-1\leq a_i< a_1+n$. 
The last inequality means
that the interval $T^n(\tilde c_1)$ lies on the right of $\tilde c_i$
in $\pi^{-1}(\partial B)$, and therefore $z\leq \tilde V_i$
implies that $z\leq T^n(\tilde V_1)$.
Then by (2), we have  $$\Pi_1(z)\leq a_1+n+M-1\leq a_i+M.$$

Let us start the proof of (1).
Let $\delta$ be a simple curve in $V$ 
joining  $\pi(a_1)$
to $\pi(b_1)$ which is not homotopic to $\pi([a_1,b_1])$, and let
$\gamma=\pi([b_1,a_1+1])\subset\partial B$.
 Choose $\delta$ so that the concatenation $\delta\cdot\gamma$ is a
simple closed curve which
bounds a closed disc $D$ containing $Z$ in its interior. This is
possible
because $Z$ is a compactum not separating both ends of $\AAA$.
There is a lift $\tilde D$ of $D$ which is bounded by 
the concatenation $\tilde \delta\cdot\tilde\gamma$,
where $\tilde \delta$ is a lift of $\delta$
and $\tilde\gamma=[b_1,a_1+1]$.
Let $Z_0=\pi^{-1}(Z)\cap \tilde D$.
Then we have $\pi^{-1}(Z)=\coprod_{i\in\Z}T^i(Z_0)$.

We shall show that the point $z\in \pi^{-1}(Z)$ satisfying $z\leq \tilde V_1$
is contained in $T^i(Z_0)$ for some $i\leq0$. Clearly this is sufficient 
for our purpose since $Z_0$ is
compact. Assume the contrary, say, $z\in T(Z_0)$. 
Since $z\leq \tilde V_1$, there is a curve $l$ in $\Gamma_1$ which
contains $z\in T(Z_0)$ on its left side. Let $t_0$ be the smallest
value such that $l(t_0)\in(a_1,b_1)$. 
The curve $l$ is homotopic in the family $\Gamma_1$ to a curve, still denoted by
$l$,
such that $l(t_0,\infty)$ is contained in $\pi^{-1}(V)$.
It can further be homotoped so that $l(t_0,\infty)$ does not intersect 
the disc
$T(\tilde D)$, since $\pi^{-1}(V)$ is simply connected.

The other half of the curve, $l((-\infty,t_0))$, is
contained in $\tilde V_1$. It must intersect $[b_1+1,a_1+2]$,
the lower boundary of $T(\tilde D)$,
since a point $z\in T(D)$ is still on the left side of the new curve $l$.

Consider the curve $T\circ l$.  The two curves
$l(-\infty, t_0)$ and $T\circ l(-\infty,t_0)$ must intersect. See Figure 6.
But the former is contained in $\tilde V_1$ while the latter in
$T(\tilde V_1)$. Since $\tilde V_1\cap T(\tilde V_1)\neq\emptyset$,
this is impossible. \qed

\begin{figure}
\unitlength 0.1in
\begin{picture}( 45.8000, 30.1500)( 16.7000,-37.3000)
%
{\color[named]{Black}{%
\special{pn 8}%
\special{pa 1770 2260}%
\special{pa 6250 2260}%
\special{fp}%
}}%
%
{\color[named]{Black}{%
\special{pn 8}%
\special{pa 1670 1090}%
\special{pa 1694 1114}%
\special{pa 1718 1136}%
\special{pa 1742 1160}%
\special{pa 1766 1182}%
\special{pa 1788 1206}%
\special{pa 1812 1230}%
\special{pa 1878 1302}%
\special{pa 1918 1350}%
\special{pa 1938 1376}%
\special{pa 1956 1400}%
\special{pa 1992 1452}%
\special{pa 2008 1480}%
\special{pa 2024 1506}%
\special{pa 2036 1534}%
\special{pa 2050 1562}%
\special{pa 2062 1590}%
\special{pa 2072 1620}%
\special{pa 2082 1648}%
\special{pa 2092 1678}%
\special{pa 2100 1708}%
\special{pa 2106 1740}%
\special{pa 2112 1770}%
\special{pa 2124 1834}%
\special{pa 2136 1930}%
\special{pa 2138 1962}%
\special{pa 2140 1996}%
\special{pa 2142 2028}%
\special{pa 2144 2062}%
\special{pa 2146 2094}%
\special{pa 2148 2128}%
\special{pa 2148 2162}%
\special{pa 2150 2196}%
\special{pa 2150 2250}%
\special{fp}%
}}%
%
{\color[named]{Black}{%
\special{pn 8}%
\special{pa 2550 2270}%
\special{pa 2546 2232}%
\special{pa 2540 2192}%
\special{pa 2534 2154}%
\special{pa 2530 2114}%
\special{pa 2522 2038}%
\special{pa 2520 2002}%
\special{pa 2518 1964}%
\special{pa 2518 1894}%
\special{pa 2520 1860}%
\special{pa 2524 1826}%
\special{pa 2528 1794}%
\special{pa 2534 1764}%
\special{pa 2544 1734}%
\special{pa 2554 1706}%
\special{pa 2566 1680}%
\special{pa 2580 1656}%
\special{pa 2598 1632}%
\special{pa 2618 1610}%
\special{pa 2640 1592}%
\special{pa 2666 1574}%
\special{pa 2694 1558}%
\special{pa 2724 1546}%
\special{pa 2756 1534}%
\special{pa 2788 1526}%
\special{pa 2822 1520}%
\special{pa 2854 1518}%
\special{pa 2886 1518}%
\special{pa 2914 1522}%
\special{pa 2942 1528}%
\special{pa 2966 1538}%
\special{pa 2988 1550}%
\special{pa 3008 1566}%
\special{pa 3026 1584}%
\special{pa 3042 1606}%
\special{pa 3054 1630}%
\special{pa 3066 1656}%
\special{pa 3076 1684}%
\special{pa 3086 1714}%
\special{pa 3092 1746}%
\special{pa 3104 1814}%
\special{pa 3108 1890}%
\special{pa 3110 1930}%
\special{pa 3110 2054}%
\special{pa 3108 2096}%
\special{pa 3102 2228}%
\special{pa 3100 2260}%
\special{fp}%
}}%
%
{\color[named]{Black}{%
\special{pn 8}%
\special{pa 3710 2260}%
\special{pa 3726 2232}%
\special{pa 3740 2204}%
\special{pa 3820 2064}%
\special{pa 3836 2038}%
\special{pa 3854 2010}%
\special{pa 3870 1984}%
\special{pa 3888 1958}%
\special{pa 3948 1880}%
\special{pa 3968 1856}%
\special{pa 4012 1808}%
\special{pa 4036 1784}%
\special{pa 4060 1762}%
\special{pa 4138 1702}%
\special{pa 4194 1670}%
\special{pa 4224 1656}%
\special{pa 4254 1644}%
\special{pa 4284 1634}%
\special{pa 4314 1628}%
\special{pa 4346 1622}%
\special{pa 4376 1620}%
\special{pa 4440 1620}%
\special{pa 4472 1624}%
\special{pa 4506 1628}%
\special{pa 4570 1640}%
\special{pa 4602 1648}%
\special{pa 4634 1658}%
\special{pa 4664 1666}%
\special{pa 4696 1678}%
\special{pa 4726 1688}%
\special{pa 4756 1700}%
\special{pa 4846 1742}%
\special{pa 4902 1774}%
\special{pa 4928 1790}%
\special{pa 4956 1808}%
\special{pa 4982 1826}%
\special{pa 5006 1846}%
\special{pa 5032 1866}%
\special{pa 5056 1886}%
\special{pa 5078 1908}%
\special{pa 5102 1930}%
\special{pa 5146 1978}%
\special{pa 5166 2002}%
\special{pa 5188 2026}%
\special{pa 5208 2050}%
\special{pa 5252 2098}%
\special{pa 5272 2122}%
\special{pa 5294 2146}%
\special{pa 5316 2168}%
\special{pa 5340 2190}%
\special{pa 5362 2214}%
\special{pa 5386 2236}%
\special{pa 5400 2250}%
\special{fp}%
}}%
%
{\color[named]{Black}{%
\special{pn 20}%
\special{pa 3720 2250}%
\special{pa 5400 2250}%
\special{fp}%
}}%
%
{\color[named]{Black}{%
\special{pn 20}%
\special{pa 3710 2240}%
\special{pa 3702 2210}%
\special{pa 3694 2178}%
\special{pa 3688 2146}%
\special{pa 3680 2114}%
\special{pa 3674 2082}%
\special{pa 3668 2052}%
\special{pa 3662 2020}%
\special{pa 3658 1988}%
\special{pa 3654 1924}%
\special{pa 3654 1892}%
\special{pa 3658 1828}%
\special{pa 3662 1796}%
\special{pa 3670 1762}%
\special{pa 3678 1730}%
\special{pa 3688 1700}%
\special{pa 3700 1668}%
\special{pa 3714 1638}%
\special{pa 3730 1610}%
\special{pa 3746 1584}%
\special{pa 3766 1558}%
\special{pa 3786 1536}%
\special{pa 3808 1514}%
\special{pa 3830 1496}%
\special{pa 3856 1480}%
\special{pa 3908 1452}%
\special{pa 3938 1440}%
\special{pa 3968 1430}%
\special{pa 3998 1422}%
\special{pa 4030 1416}%
\special{pa 4062 1408}%
\special{pa 4096 1404}%
\special{pa 4130 1398}%
\special{pa 4164 1394}%
\special{pa 4200 1390}%
\special{pa 4236 1388}%
\special{pa 4308 1380}%
\special{pa 4346 1378}%
\special{pa 4382 1374}%
\special{pa 4420 1370}%
\special{pa 4458 1368}%
\special{pa 4494 1364}%
\special{pa 4532 1362}%
\special{pa 4570 1358}%
\special{pa 4606 1356}%
\special{pa 4644 1354}%
\special{pa 4716 1350}%
\special{pa 4858 1350}%
\special{pa 4926 1354}%
\special{pa 5022 1366}%
\special{pa 5052 1372}%
\special{pa 5136 1396}%
\special{pa 5162 1408}%
\special{pa 5186 1420}%
\special{pa 5208 1432}%
\special{pa 5230 1446}%
\special{pa 5252 1462}%
\special{pa 5288 1498}%
\special{pa 5304 1520}%
\special{pa 5332 1564}%
\special{pa 5344 1590}%
\special{pa 5364 1642}%
\special{pa 5372 1670}%
\special{pa 5384 1730}%
\special{pa 5390 1762}%
\special{pa 5394 1794}%
\special{pa 5398 1828}%
\special{pa 5404 1930}%
\special{pa 5406 1966}%
\special{pa 5406 2076}%
\special{pa 5404 2114}%
\special{pa 5404 2152}%
\special{pa 5402 2190}%
\special{pa 5402 2228}%
\special{pa 5400 2250}%
\special{fp}%
}}%
%
{\color[named]{Black}{%
\special{pn 13}%
\special{pa 2350 1310}%
\special{pa 2348 1342}%
\special{pa 2346 1372}%
\special{pa 2342 1402}%
\special{pa 2340 1432}%
\special{pa 2336 1464}%
\special{pa 2334 1494}%
\special{pa 2332 1526}%
\special{pa 2330 1556}%
\special{pa 2326 1586}%
\special{pa 2324 1618}%
\special{pa 2322 1648}%
\special{pa 2320 1680}%
\special{pa 2318 1710}%
\special{pa 2316 1742}%
\special{pa 2314 1772}%
\special{pa 2314 1804}%
\special{pa 2312 1836}%
\special{pa 2310 1866}%
\special{pa 2310 1898}%
\special{pa 2308 1930}%
\special{pa 2308 1962}%
\special{pa 2306 1994}%
\special{pa 2306 2122}%
\special{pa 2308 2154}%
\special{pa 2308 2186}%
\special{pa 2310 2220}%
\special{pa 2310 2252}%
\special{pa 2312 2286}%
\special{pa 2314 2318}%
\special{pa 2316 2352}%
\special{pa 2320 2386}%
\special{pa 2322 2420}%
\special{pa 2338 2556}%
\special{pa 2342 2592}%
\special{pa 2348 2626}%
\special{pa 2356 2658}%
\special{pa 2366 2688}%
\special{pa 2382 2712}%
\special{pa 2402 2730}%
\special{pa 2428 2744}%
\special{pa 2458 2752}%
\special{pa 2492 2756}%
\special{pa 2528 2760}%
\special{pa 2790 2760}%
\special{pa 2828 2758}%
\special{pa 2864 2758}%
\special{pa 2902 2756}%
\special{pa 2938 2756}%
\special{pa 2976 2754}%
\special{pa 3048 2750}%
\special{pa 3084 2746}%
\special{pa 3120 2744}%
\special{pa 3154 2740}%
\special{pa 3190 2736}%
\special{pa 3258 2728}%
\special{pa 3360 2710}%
\special{pa 3392 2704}%
\special{pa 3456 2688}%
\special{pa 3486 2680}%
\special{pa 3518 2670}%
\special{pa 3548 2660}%
\special{pa 3576 2650}%
\special{pa 3606 2640}%
\special{pa 3634 2628}%
\special{pa 3662 2614}%
\special{pa 3688 2602}%
\special{pa 3716 2588}%
\special{pa 3740 2572}%
\special{pa 3766 2558}%
\special{pa 3790 2540}%
\special{pa 3814 2524}%
\special{pa 3836 2504}%
\special{pa 3858 2486}%
\special{pa 3878 2466}%
\special{pa 3918 2422}%
\special{pa 3936 2400}%
\special{pa 3954 2376}%
\special{pa 3970 2350}%
\special{pa 3986 2326}%
\special{pa 4002 2298}%
\special{pa 4030 2242}%
\special{pa 4044 2212}%
\special{pa 4056 2182}%
\special{pa 4070 2152}%
\special{pa 4082 2120}%
\special{pa 4096 2088}%
\special{pa 4110 2054}%
\special{pa 4124 2022}%
\special{pa 4138 1992}%
\special{pa 4154 1962}%
\special{pa 4174 1936}%
\special{pa 4194 1912}%
\special{pa 4216 1894}%
\special{pa 4242 1880}%
\special{pa 4272 1870}%
\special{pa 4304 1866}%
\special{pa 4336 1866}%
\special{pa 4370 1870}%
\special{pa 4404 1878}%
\special{pa 4434 1888}%
\special{pa 4494 1916}%
\special{pa 4546 1952}%
\special{pa 4572 1972}%
\special{pa 4594 1992}%
\special{pa 4618 2012}%
\special{pa 4638 2034}%
\special{pa 4660 2056}%
\special{pa 4678 2078}%
\special{pa 4698 2100}%
\special{pa 4730 2148}%
\special{pa 4746 2174}%
\special{pa 4788 2252}%
\special{pa 4800 2278}%
\special{pa 4820 2334}%
\special{pa 4830 2364}%
\special{pa 4838 2392}%
\special{pa 4854 2452}%
\special{pa 4866 2512}%
\special{pa 4872 2544}%
\special{pa 4880 2608}%
\special{pa 4882 2640}%
\special{pa 4886 2672}%
\special{pa 4888 2706}%
\special{pa 4888 2740}%
\special{pa 4890 2772}%
\special{pa 4890 2876}%
\special{pa 4888 2910}%
\special{pa 4886 2946}%
\special{pa 4886 2982}%
\special{pa 4882 3018}%
\special{pa 4880 3052}%
\special{pa 4878 3090}%
\special{pa 4866 3198}%
\special{pa 4862 3236}%
\special{pa 4858 3272}%
\special{pa 4854 3310}%
\special{pa 4848 3346}%
\special{pa 4844 3384}%
\special{pa 4838 3422}%
\special{pa 4834 3458}%
\special{pa 4816 3572}%
\special{pa 4812 3610}%
\special{pa 4806 3648}%
\special{pa 4800 3680}%
\special{fp}%
}}%
%
{\color[named]{Black}{%
\special{pn 13}%
\special{pa 2350 1310}%
\special{pa 2400 1050}%
\special{fp}%
\special{sh 1}%
\special{pa 2400 1050}%
\special{pa 2368 1112}%
\special{pa 2390 1102}%
\special{pa 2408 1120}%
\special{pa 2400 1050}%
\special{fp}%
}}%
%
{\color[named]{Black}{%
\special{pn 13}%
\special{pa 3440 1250}%
\special{pa 3438 1282}%
\special{pa 3438 1314}%
\special{pa 3432 1410}%
\special{pa 3432 1442}%
\special{pa 3428 1506}%
\special{pa 3428 1538}%
\special{pa 3422 1634}%
\special{pa 3422 1666}%
\special{pa 3418 1730}%
\special{pa 3418 1762}%
\special{pa 3414 1826}%
\special{pa 3414 1858}%
\special{pa 3412 1890}%
\special{pa 3412 1922}%
\special{pa 3410 1954}%
\special{pa 3410 1986}%
\special{pa 3408 2018}%
\special{pa 3408 2050}%
\special{pa 3406 2082}%
\special{pa 3406 2114}%
\special{pa 3404 2146}%
\special{pa 3404 2178}%
\special{pa 3402 2210}%
\special{pa 3402 2274}%
\special{pa 3400 2306}%
\special{pa 3400 2434}%
\special{pa 3398 2466}%
\special{pa 3398 2690}%
\special{pa 3400 2722}%
\special{pa 3400 2818}%
\special{pa 3402 2850}%
\special{pa 3402 2914}%
\special{pa 3404 2946}%
\special{pa 3404 2978}%
\special{pa 3406 3010}%
\special{pa 3406 3042}%
\special{pa 3408 3074}%
\special{pa 3408 3106}%
\special{pa 3412 3170}%
\special{pa 3412 3202}%
\special{pa 3418 3298}%
\special{pa 3418 3330}%
\special{pa 3428 3490}%
\special{pa 3428 3522}%
\special{pa 3440 3714}%
\special{pa 3440 3730}%
\special{fp}%
}}%
%
{\color[named]{Black}{%
\special{pn 13}%
\special{pa 3430 1250}%
\special{pa 3480 1040}%
\special{fp}%
\special{sh 1}%
\special{pa 3480 1040}%
\special{pa 3446 1100}%
\special{pa 3468 1092}%
\special{pa 3484 1110}%
\special{pa 3480 1040}%
\special{fp}%
}}%
%
{\color[named]{Black}{%
\special{pn 8}%
\special{pa 5840 2260}%
\special{pa 5836 2228}%
\special{pa 5830 2196}%
\special{pa 5826 2162}%
\special{pa 5822 2130}%
\special{pa 5816 2096}%
\special{pa 5804 2000}%
\special{pa 5802 1966}%
\special{pa 5798 1934}%
\special{pa 5794 1870}%
\special{pa 5794 1838}%
\special{pa 5792 1806}%
\special{pa 5792 1774}%
\special{pa 5794 1742}%
\special{pa 5796 1712}%
\special{pa 5798 1680}%
\special{pa 5802 1648}%
\special{pa 5806 1618}%
\special{pa 5810 1586}%
\special{pa 5816 1556}%
\special{pa 5832 1496}%
\special{pa 5852 1436}%
\special{pa 5864 1406}%
\special{pa 5874 1376}%
\special{pa 5888 1346}%
\special{pa 5900 1316}%
\special{pa 5914 1288}%
\special{pa 5942 1228}%
\special{pa 5956 1200}%
\special{pa 5970 1170}%
\special{pa 5980 1150}%
\special{fp}%
}}%
\put(42.6000,-8.5000){\makebox(0,0){$\pi^{-1}(V)$}}%
\put(50.0000,-16.1000){\makebox(0,0){$T(\tilde D)$}}%
%
{\color[named]{Black}{%
\special{pn 4}%
\special{sh 1}%
\special{ar 4350 2180 8 8 0  6.28318530717959E+0000}%
\special{sh 1}%
\special{ar 4290 2600 8 8 0  6.28318530717959E+0000}%
}}%
\put(43.4000,-20.4000){\makebox(0,0){$z$}}%
\put(20.2000,-24.6000){\makebox(0,0){$a_1$}}%
\put(25.0000,-24.7000){\makebox(0,0){$b_1$}}%
\put(30.5000,-24.6000){\makebox(0,0){$a_1+1$}}%
\put(38.8000,-29.3000){\makebox(0,0){$b_1+1$}}%
%
{\color[named]{Black}{%
\special{pn 8}%
\special{pa 3800 2790}%
\special{pa 3750 2350}%
\special{fp}%
\special{sh 1}%
\special{pa 3750 2350}%
\special{pa 3738 2418}%
\special{pa 3756 2404}%
\special{pa 3778 2414}%
\special{pa 3750 2350}%
\special{fp}%
}}%
\put(24.3000,-7.8000){\makebox(0,0){$l$}}%
\put(34.6000,-7.9000){\makebox(0,0){$T(l)$}}%
\end{picture}%

\caption{}
\end{figure}

\medskip
To finish, let $z\in\pi^{-1}(Z)$. One may assume $z\leq \tilde V_1$
by replacing $z$ by $T^{-n}(z)$ if necessary. Then by successive use
of Lemma \ref{l5.1},
we have $\tilde h^i(z)\leq\tilde V_i$ for any $i\in\N$.
Then by Lemma \ref{l5.2}, we have $$\Pi_1(\tilde h^i(z))\leq a_i+M,$$
showing that
$$
\lim_{i\to\infty}i^{-1}\Pi_1(\tilde h^i(z))\leq\lim i^{-1}a_i=\alpha,$$
completing the proof of Theorem \ref{main4} (1).

To show (2), just consider a lift of a  point in $Z$ accessible
from $U_{-\infty}$. Details are left to the reader.

\section{Appendix: $C^\infty$ complete Lyapunov functions}

We fix $h\in\HH$. Here is a criterion of the chain recurrent set $C$
and a chain transitive class in terms of attractors and repellors (\cite{F}).
A subset $A_i$ in $S^2$ is called an {\em attractor} if there is an open
neighbourhood $V_i$ of $A_i$ such that $h(\Cl(V_i))\subset V_i$ and
$ \bigcap_{j\geq 0}f^j(\Cl(V_i))=A_i$. 
The set $V_i$ is called an {\em isolating block} of $A_i$,
and the set
$  A_i^*=\bigcap_{j\geq0}f^{-j}(S^2\setminus V_i)$ 
the {\em dual repellor} of $A$.
The totality of attractors is at most countable, and we denote it
by $\{A_i\}_{i\in I}$. Then we have (\cite{F})
$$
C=\bigcap_{i\in I}(A_i\cup A_i^*).$$
For $x,\, y\in C$, we also have 
$$
x\sim y\, \Longleftrightarrow\, \forall i\in I, \ \mbox{ either }\
x,\,y\in A_i\ \mbox{ or}\ x,\,y\in A_i^*.$$

 We begin with the following well known
fact due to H. Whitney.

\begin{lemma} \label{Whitney}
For any closed subset $P$ in $S^2$, there is a $C^\infty$ function
$\varphi_P:S^2\to[0,1]$ such that $\varphi_P^{-1}(0)=P$. \qed
\end{lemma}

\begin{lemma}
For any disjoint closed subets $P$ and $Q$ of $S^1$, there is
a $C^\infty$ function
$\psi:S^2\to[0,1]$ such that $\psi^{-1}(0)=P$ and
 $\psi^{-1}(1)=Q$.
\end{lemma}

\bd The function $\varphi_P$ in Lemma \ref{Whitney} can easily
be modified so as to satisfy $Q\subset\varphi_P^{-1}(1)$.
Define a function $\varphi_Q$ replacing the roles of $P$ and $Q$, and set 
$$
\psi=2^{-1}(\varphi_P+1-\varphi_Q).$$
\qed

Recall that $\{A_i\}_{i\in I}$ is the family of the attractors of $h$.

\begin{lemma} \label{la}
For each $i\in I$, there is a $C^\infty$ function $H_i:S^2\to[0,1]$
such that 
\begin{enumerate}
\item
$H_i^{-1}(0)=A_i$ and $H_i^{-1}(1)=A^*_i$. 
\item
For any $x\in S^2\setminus(A_i\cup A_i^*)$, we have $H_i(h(x))<H_i(x)$.
\end{enumerate}
\end{lemma}

\bd Let $V_i$ be an isolationg block of $A_i$.
Then for any $j\in\Z$, there is a $C^\infty$ function $\psi_j:S^2\to [0,1]$
such that $\psi_j^{-1}(0)=f^j(\Cl(V_i))$ and
$\psi_j^{-1}(1)=S^2\setminus
f^{j-1}(V_i))$.
Choose a sequence $c_j>0$ such that
$$\sum_{j\in\N}c_j\Vert
\psi_j\Vert_{\abs{j}}<\infty,$$ where $\Vert\cdot\Vert_{\abs{j}}$ denotes the $C^{\abs{j}}$ norm.
Then the function
$  \sum_{j\in\Z}c_j\psi_j$ is a $C^\infty$ function, and after normalized
it satisfies the conditions of Lemma \ref{la}. 
\qed

\medskip
\noindent
{\bf Proof of Proposition \ref{p1}.} By an appropriate choice of
 positive numbers
$a_i$, 
the function $  H=\sum_{i\in I}a_i H_i$
is a $C^\infty$ function satisfying (1) and (2) of Definition
\ref{d2}. 
If the indexing set $I$ is infinite, set
$I=\N$ and choose $a_i$ such that $a_{i+1}<3^{-1}a_i$ ($\forall i$).
Then we obtain that $H(C)$ is closed and that the Lebesgue measure
of $H(C)$ is zero.

\end{document}